\documentclass[12pt]{article}
\usepackage[utf8]{inputenc}
\usepackage[russian,english]{babel}
\usepackage{amssymb}
\usepackage{bm}
\usepackage{booktabs,caption}
\usepackage{multirow}
\usepackage{graphicx}
\usepackage{amsmath,amsthm}
\usepackage{a4wide}
\usepackage{version}
\usepackage{parskip}
\usepackage{tikz}
\usepackage{mathtools}
\usetikzlibrary{arrows,decorations.markings, matrix}
\usepackage{hyperref}
\usepackage{comment}

\usepackage[all]{xy}

\numberwithin{equation}{section}

\newcommand{\BB}{\mathcal{B}}
\newcommand{\CC}{\mathbb{C}}

\newcommand{\eE}{\mathcal{E}}

\newcommand{\QQ}{\mathbb{Q}}
\newcommand{\FF}{\mathbb{F}}
\newcommand{\GG}{\mathbb{G}}
\newcommand{\HH}{\mathbb{H}}
\newcommand{\LL}{\mathcal{L}}
\newcommand{\RR}{\mathbb{R}}
\newcommand{\OO}{\mathcal{O}}
\newcommand{\PP}{\mathbb{P}}

\newcommand{\XX}{\mathcal{X}}
\newcommand{\YY}{\mathcal{Y}}
\newcommand{\ZZ}{\mathbb{Z}}

\newcommand{\OP}[1]{\mathrm{#1}}
\newcommand{\HP}{\HH\PP}
\newcommand{\Utrop}{U^{\mathrm{trop}}}
\newcommand{\Vtrop}{V^{\mathrm{trop}}}
\newcommand{\Trop}{\OP{Trop}^\vee}
\newcommand{\bin}{\mathfrak{b}}
\newcommand{\sslash}{\mathbin{/\mkern-6mu/}}
\newcommand{\sm}[1]{#1^\circ}
  
\tikzset{->-/.style={decoration={
              markings,
              mark=at position .5 with {\arrow{>}}},postaction={decorate}}}

\tikzset{-<-/.style={decoration={
              markings,
              mark=at position .5 with {\arrow{<}}},postaction={decorate}}}

\title{Tropical methods for stable octic double planes}

\author{Jonny Evans\hspace{2cm} Angelica Simonetti\hspace{2cm} Giancarlo Urz\'{u}a}

\newcommand{\pg}{\paragraph{\hspace{-0.6cm}}}

\usepackage{caption}
\usepackage{subcaption}
\usepackage{accents}
\usetikzlibrary{shapes.multipart}
\usetikzlibrary{decorations.pathmorphing,shapes}
\usetikzlibrary{decorations.pathreplacing}
\DeclareRobustCommand{\dbtilde}[1]{\accentset{\approx}{#1}}

\begin{document}

\maketitle

\begin{abstract}
  This paper has been written to illustrate the power of
  techniques from tropical geometry and mirror symmetry for
  studying the KSBA moduli space of surfaces on or near the
  Noether line. We focus on the moduli space of octic double
  planes (\(K^2=2\), \(p_g=3\)) and use methods from tropical
  and toric geometry to classify the strata corresponding to
  normal KSBA-stable surfaces, focusing on the non-Gorenstein
  case.
\end{abstract}

\section{Introduction}

\pg The moduli space of smooth complex surfaces of general type with
fixed Chern numbers is not compact, but it admits a
compactification suggested by Koll\'{a}r and Shepherd-Barron
\cite{KSB}. This compactified moduli space was shown to be
separated by Koll\'{a}r and Shepherd-Barron in their original
paper, and was shown to be compact by Alexeev \cite{Alexeev}, so
it is usually called the {\em
  Koll\'{a}r--Shepherd-Barron--Alexeev (KSBA) moduli space}. See
the recent comprehensive monograph by Koll\'{a}r
\cite{KollarFamilies} for the state of the art. The boundary
points of the KSBA moduli space correspond to certain singular
surfaces called {\em KSBA-stable surfaces}. Whilst the KSBA
boundary is uncharted territory in general, an explicit
understanding of KSBA limits for specific choices of Chern
numbers has been achieved in the last decade; progress has been
made for \(K^2=5\), \(p_g=4\) (quintics) \cite{RanaQuintics},
for \(K^2=1\), \(p_g=2\) \cite{CFPRR,FPRR,FPR,GPSZ,RU} and for
\(K^2=6\), \(p_g=0\) (Burniat and Campedelli surfaces with
certain fundamental groups)
\cite{AlexeevPardiniBurniatCampedelli}.

\pg In this paper, we will focus on another specific case, when
\(K^2=2\) and \(p_g=3\), in order to show how ideas from
tropical geometry and mirror symmetry can help to organise the
case analysis. The smooth surfaces in this moduli space are {\em
  octic double planes}, that is double covers of \(\PP^2\)
branched over a smooth octic curve
{\cite[\selectlanguage{russian} Глава VI \S 3 Теорема
  6]{Shafarevich}}. One way to produce degenerations of double
branched covers \(X\to \PP^2\) is to allow the branch curve to
degenerate; this can only produce Gorenstein singularities. The
{\em Gorenstein} KSBA-stable limits of octic double planes were
classified by Anthes \cite{Anthes}, who also gave the first
examples of non-Gorenstein limits {\cite[Example
  5.5]{Anthes}}. Still more examples can be found if we also
allow \(\PP^2\) to degenerate as well. The degenerations of
\(\PP^2\) are very well-understood thanks to work of B\u{a}descu
\cite{Badescu1,Badescu2}, Manetti \cite{Manetti1,ManettiThesis}
and Hacking and Prokhorov \cite{HP2,HP1}; this will allow us to
use techniques from {\em tropical} geometry from the work of
Gross, Hacking, Keel and Kontsevich \cite{GHK2,GHK1,GHKK} to
understand the degeneration of the double branched cover. We
will be able to understand which non-Gorenstein singularities
appear for normal surfaces. The integral affine geometry of some
associated tropicalisations is the key ingredient which allows
us to tame the case analysis. These tropicalisations play a role
analogous to the Newton polygon in toric geometry but for
certain non-toric degenerations of \(\PP^2\). Tropical geometry
allows us to read off information easily about discrepancies of
singularities (see \ref{pg:necessary_condition}) and the
operation of {\em mutation} which connects the tropicalisations
of different degenerations of \(\PP^2\) gives us an inductive
technique for ruling out infinitely many cases in Section
\ref{sct:reducing_to_finite}.

\pg We now state our main theorem. We always work over
\(\CC\). We denote by \(\HP(5)\) the partial smoothing of
\(\PP(1,4,25)\) which keeps the \(\frac{1}{25}(1,4)\)
singularity and smooths the \(\frac{1}{4}(1,1)\)
singularity. This surface can be realised as the hypersurface
\(\{x_1x_4+x_2^{13}+x_3^2=0\}\subset\PP(1,2,13,25)\); we will
discuss this and related surfaces in more detail later. See also
Appendix \ref{sct:almost_toric} for a detailed description of
\(\HP(5)\).

\paragraph{Theorem.}\label{thm:main_thm} {\em Let \(X\) be a
  normal KSBA-stable surface which admits a \(\QQ\)-Gorenstein
  smoothing whose general fibre is a smooth octic double
  plane. Then \(X\) is one of the following:
  \begin{itemize}
  \item[1.] A double cover of \(\PP^2\) branched over a singular
    octic. In this case \(X\) has at worst Gorenstein
    singularities as classified by Anthes \cite{Anthes}.
  \item[2.] A double cover of \(Y=\PP(1,1,4)\) branched over a curve
    of weighted degree \(16\). In addition to any Gorenstein
    singularities coming from singularities of the branch curve,
    \(X\) may have precisely one of the following baskets of
    singularities:
    \begin{itemize}
    \item[I.] Two \(\frac{1}{4}(1,1)\) singularities. These
      surfaces were constructed by Anthes {\cite[Example
        5.5]{Anthes}}.
    \item[II.] A singularity whose minimal resolution has
      one of the dual graphs shown in Figure
      \ref{fig:graphs}(a-c). These are \(\ZZ/2\)-quotients of
      (a) simple elliptic or (b--c) cusp singularities. 
    \end{itemize}
  \item[3.] A double cover of \(Y=\PP(1,4,25)\) branched over a
    curve of weighted degree \(80\). In addition any Gorenstein singularities coming from singularities of the branch curve, \(X\) has a \(\frac{1}{50}(1,29)\)
    singularity and precisely one of the baskets listed for
    \(\PP(1,1,4)\).
  \item[4.] A double cover of \(Y=\HP(5)\). In addition any Gorenstein singularities, \(X\) has a \(\frac{1}{50}(1,29)\)
    singularity.
  \end{itemize}}

\paragraph{Remark.} We have not studied which combinations of
Gorenstein and non-Gorenstein singularities can appear.

\begin{figure}[htb]
  \begin{center}
    \begin{tikzpicture}
      \node (a) at (0,0) {\(\bullet\)};
      \node at (a) [right] {\footnotesize \(-4\)};
      \node (b) at (45:1) {\(\bullet\)};
      \node at (b) [right] {\footnotesize \(-2\)};
      \node (c) at (-45:1) {\(\bullet\)};
      \node at (c) [right] {\footnotesize \(-2\)};
      \node (d) at (135:1) {\(\bullet\)};
      \node at (d) [left] {\footnotesize \(-2\)};
      \node (e) at (-135:1) {\(\bullet\)};
      \node at (e) [left] {\footnotesize \(-2\)};
      \draw (a) -- (b);
      \draw (a) -- (c);
      \draw (a) -- (d);
      \draw (a) -- (e);
      \node at (-2,0) {(a)};
      \begin{scope}[shift={(5,0)}]
        \node (a) at (0,0) {\(\bullet\)};
        \node (b) at (1,0) {\(\cdots\)};
        \node (c) at (2,0) {\(\bullet\)};
        \node (d) at (3,0) {\(\bullet\)};
        \node at (d) [right] {\footnotesize \(-4\)};
        \node at (a) [above] {\footnotesize \(-2\)};
        \node at (c) [above] {\footnotesize \(-2\)};
        \node (e) at ({3+cos(-45)},{sin(-45)}) {\(\bullet\)};
        \node at (e) [right] {\footnotesize \(-2\)};
        \node (f) at ({3+cos(45)},{sin(45)}) {\(\bullet\)};
        \node at (f) [right] {\footnotesize \(-2\)};
        \node (g) at (135:1) {\(\bullet\)};
        \node at (g) [left] {\footnotesize \(-2\)};
        \node (h) at (-135:1) {\(\bullet\)};
        \node at (h) [left] {\footnotesize \(-2\)};
        \draw (a) -- (b);
        \draw (b) -- (c);
        \draw (c) -- (d);
        \draw (d) -- (f);
        \draw (d) -- (e);
        \draw (a) -- (g);
        \draw (a) -- (h);
        \node at (-2,0) {(b)};
        \draw [decorate,decoration={brace,amplitude=5pt,mirror}] (0,-0.5) -- (2,-0.5) node[midway,yshift=-1em]{\(t\)};
      \end{scope}
      \begin{scope}[shift={(0,-3)}]
        \node (a) at (0,0) {\(\bullet\)};
        \node (b) at (1,0) {\(\cdots\)};
        \node (c) at (2,0) {\(\bullet\)};
        \node (d) at (3,0) {\(\bullet\)};
        \node (i) at (4,0) {\(\bullet\)};
        \node (j) at (5,0) {\(\cdots\)};
        \node (k) at (6,0) {\(\bullet\)};
        \node at (d) [above] {\footnotesize \(-4\)};
        \node at (a) [above] {\footnotesize \(-2\)};
        \node at (c) [above] {\footnotesize \(-2\)};
        \node at (i) [above] {\footnotesize \(-2\)};
        \node at (k) [above] {\footnotesize \(-2\)};
        \node (e) at ({6+cos(-45)},{sin(-45)}) {\(\bullet\)};
        \node at (e) [right] {\footnotesize \(-2\)};
        \node (f) at ({6+cos(45)},{sin(45)}) {\(\bullet\)};
        \node at (f) [right] {\footnotesize \(-2\)};
        \node (g) at (135:1) {\(\bullet\)};
        \node at (g) [left] {\footnotesize \(-2\)};
        \node (h) at (-135:1) {\(\bullet\)};
        \node at (h) [left] {\footnotesize \(-2\)};
        \draw (a) -- (b);
        \draw (b) -- (c);
        \draw (c) -- (d);
        \draw (d) -- (i);
        \draw (i) -- (j);
        \draw (j) -- (k);
        \draw (a) -- (g);
        \draw (a) -- (h);
        \draw (k) -- (e);
        \draw (k) -- (f);
        \node at (-2,0) {(c)};
        \draw [decorate,decoration={brace,amplitude=5pt,mirror}] (0,-0.5) -- (2,-0.5) node[midway,yshift=-1em]{\(t_1\)};
        \draw [decorate,decoration={brace,amplitude=5pt,mirror}] (4,-0.5) -- (6,-0.5) node[midway,yshift=-1em]{\(t_2\)};
  \end{scope}
    \end{tikzpicture}
    \caption{The dual graphs for minimal resolutions of the
      singularities appearing in Case II of Theorem
      \ref{thm:main_thm}.}
    \label{fig:graphs}
  \end{center}
\end{figure}

\paragraph{Remark.} In an earlier version of this paper, we
included the graphs from Figure \ref{fig:graphs}(b) and (c) with
\(t=1\) and \(t_1=t_2=1\), but missed the possibility that
\(t,t_1,t_2\) could be greater than \(1\). We are indebted to
H. Akaike, M. Enokizono, M. Hattori and Y. Koto for pointing out
this gap. In \ref{pg:bound_on_ts}, we show that \(t\leq 19\) for
graph (b) and that this equality is sharp, and that
\(t_1+t_2\leq 38\) for graph (c), however the most extreme
example we have found has \(t_1+t_2\leq 21\). We conjecture that
this stronger bound should hold.

\paragraph{Remark.} In principle, the same techniques should
apply more generally to double covers of other rational
surfaces, but the classification of degenerations would be
substantially more complicated. Our original motivation for
studying this case was to understand the more general case of
{\em Horikawa surfaces}: surfaces which lie on or near the
Noether line \(2p_g = K^2+4\) and which were studied
exhaustively by Horikawa in his seminal sequence of papers
\cite{Hor1,Hor2,Hor3,Hor4,Hor5}. These all arise as double
covers of rational surfaces. For some interesting recent
progress on KSBA-stable Horikawa surfaces, see the papers of
Rana and Rollenske \cite{RanaRollenske} and Monreal, Negrete and
Urz\'{u}a \cite{MNU}.

\paragraph{Remark.} The problem of studying KSBA limits of
branch covers of \(\PP^2\) is closely related to studying KSBA
limits of pairs \((\PP^2,kB)\) for \(B\) a plane curve. Such
KSBA moduli spaces of pairs have been introduced by Hacking
\cite{HackingPlaneCurves} and were the basis for Anthes's
approach in \cite{Anthes} and for that of Alexeev and Pardini in
\cite{AlexeevPardiniBurniatCampedelli}. DeVleming and Stapleton
\cite{DeVlemingStapleton} used Hacking's ideas to study the
failure or otherwise of planarity of a curve to persist under
taking limits; their work is very much relevant to this paper
and we will use some of their results to show our stable
surfaces are smoothable. A tropical perspective on moduli of log
Calabi-Yau pairs \((Y,B)\) was recently introduced by Alexeev,
Arg\"{u}z and Bousseau \cite{AAB}.

\paragraph{Remark.} Monreal, Negrete and Urz\'{u}a \cite{MNU}
give a list of all possible combinations of cyclic quotient
singularities that could conceivably occur on a normal
KSBA-stable surface on the Noether line. They achieve this by
observing that the minimal model of the minimal resolution of
such a surface is an elliptic surface, and then proving
restrictions on the possible Hirzebruch-Jung chains of rational
curves one could find on a blow-up of an elliptic surface. They
construct stable surfaces realising all these combinations of
singularities, however, the sheaf cohomology group which
measures local-to-global obstructions to smoothing their
surfaces is not zero, and indeed most of these surfaces are not
smoothable. In the language of {\cite[Theorem 4.3, Case
  \(p_g=3\)]{MNU}}, only the Lee--Park case (i) and the family
S1F.4 with \(n=0,1\) are smoothable. Even though the other
KSBA-stable surfaces with these Chern numbers they find are not
smoothable, one can rationally blow-down their minimal
resolutions to get smooth 4-manifolds, and it would be
interesting to know if these smooth 4-manifolds are
diffeomorphic to smooth octic double surfaces. Their work raises
similar questions about other combinations of Chern numbers.

\paragraph{Overview.} We begin by reviewing some background from
algebraic and almost toric geometry in Section
\ref{sct:background}. Then in Section \ref{sct:strategy}, we
show that any normal KSBA limit of an octic double plane is a
double cover of a degeneration of \(\PP^2\) and reduce to the
case where the limit of \(\PP^2\) has at worst quotient
singularities (the {\em Manetti surfaces}). In Section
\ref{sct:manetti}, we review what is known about Manetti
surfaces and their mirror tropicalisations; in Section
\ref{sct:integer_points} we enumerate all of the integer points
in the mirror tropicalisations. In Section
\ref{sct:geom_min_res} we explain how to read off some key
geometric information about the minimal resolution from the
mirror tropicalisation. We use this in Section
\ref{sct:reducing_to_finite} to show that most Manetti surfaces
necessarily give rise to non-normal double covers. Section
\ref{sct:discrep} gives a necessary numerical criterion for the
double cover to be log canonical. Section \ref{sct:classif}
gives the full classification. We conclude in Section
\ref{sct:conclusion} with some comments about how these strata
fit together in the moduli space. Appendix
\ref{sct:almost_toric} gives a worked example which should help
the interested reader better understand how the mirror
tropicalisations are just Symington's almost toric fibrations by
another name.

\paragraph{Acknowledgements.}
JE and AS would like to thank GU and the Pontificia Universidad
Cat\'{o}lica de Chile for their hospitality during a research
visit focused on this work. We would like to thank Marcos
Canedo, Fabrizio Catanese, Kristin DeVleming, Jan Grabowski,
Paul Hacking, Vicente Monreal, Jaime Negrete, Julie Rana,
S\"{o}nke Rollenske, Jonathan Wahl, and Juan Pablo
Z\'{u}\~{n}iga for helpful conversations and
communications. Special thanks are due to the referees of the
earlier versions of this paper; their comments have helped us
improve both the results and the exposition substantially. We
are extremely grateful to Hiroto Akaike, Makoto Enokizono,
Masafumi Hattori and Yuki Koto for pointing out a gap in our
case analysis and for sharing the preliminary version of their
paper \cite{AEHK}. JE and AS are supported by EPSRC Standard
Grant EP/W015749/1. GU is supported by FONDECYT Regular Grant
1230065 and a Marie Curie FRIAS Fellowship. For the purpose of
open access, the authors have applied a Creative Commons
Attribution (CC BY) licence to any Author Accepted Manuscript
version arising.

\section{Background}
\label{sct:background}

\subsection*{Singularities and double covers}

\paragraph{Double covers.}\label{pg:ap_double_covers}
Let \(Y\) be a normal surface and \(\sm{Y}\) its smooth
locus. By normality of \(Y\), we have
\(\OP{codim}_Y(Y\setminus \sm{Y})=2\). Let \(B\subset Y\) be a
reduced Weil divisor defined by the vanishing of some section
\(\sigma\) of the sheaf \(\OO(B)\). Let \(\sm{B}=B\cap \sm{Y}\);
since \(\sm{Y}\) is smooth, \(\sm{B}\) is a Cartier divisor in
\(\sm{Y}\). Suppose that there is a line bundle
\(\phi\colon L\to \sm{Y}\) with
\(L^{\otimes 2}\cong \OO_{\sm{Y}}(\sm{B})\). Let
\(\sm{X}=\{x\in L\,:\,x^{\otimes 2}=\sigma(\phi(x))\}\) be the
double cover
\(f_0\coloneqq \phi|_{\sm{X}}\colon \sm{X}\to \sm{Y}\) branched
along \(\sm{B}\). Using a construction of Alexeev and Pardini
{\cite[Lemma 1.2]{AlexeevPardini}}, we find an \(S_2\) surface
\(X\) and a \(\ZZ_2\)-cover \(f\colon X\to Y\) extending
\(f_0\); such an extension is unique up to isomorphism. Note
that we still have \(\OP{codim}_X(X\setminus \sm{X})=2\), so
\(X\) is \(R_1\), and by Serre's condition for normality \(X\)
is itself a normal surface.

\paragraph{Discrepancies of pairs.} See {\cite[Section 2.3]{KM}}
for a thorough introduction to discrepancies. Let \(Y\) be a
surface and \(\Delta=\sum\delta_i\Delta_i\) be a Weil
\(\QQ\)-divisor, and assume that \(K_Y+\Delta\) is
\(\QQ\)-Cartier. Let \(\pi\colon Z\to Y\) be a proper birational
morphism. We can write
\begin{equation}
  \label{eq:discrep}K_Z+\pi^{-1}_*\Delta = \pi^*(K_Y+\Delta)+\sum_Ea(E,Y,\Delta)E
\end{equation}  
for some coefficients \(a(E,Y,\Delta)\), where the sum is over the
exceptional divisors \(E\) of \(\pi\); the coefficient \(a(E,Y,\Delta)\)
is called the discrepancy of \(E\) with respect to \((Y,\Delta)\). The
discrepancy of \((Y,\Delta)\) is defined to be
\[\OP{discrep}(Y,\Delta)=\inf_E(a(E,Y,\Delta))\] where
the infimum is taken over {\em all exceptional divisors of all
  possible birational maps \(\pi\)}. A pair \((Y,\Delta)\)
(respectively \(Y\) itself) is called {\em terminal}, {\em
  canonical}, {\em pure log terminal (plt)}, or {\em log
  canonical (lc)} if \(\OP{discrep}(Y,\Delta)\) (respectively
\(\OP{discrep}(Y,0)\)) is positive, nonnegative, \(>1\), or
\(\geq 1\) respectively. If \(\Delta=0\) we usually omit the
``pure'' from plt.

KSBA-stable surfaces are allowed to have {\em semi log canonical}
singularities (an extension of the log canonical condition to
non-normal surfaces). If we restrict to {\em normal} surfaces then
\(X\) is KSBA-stable if it is log canonical and has ample canonical
class.


\pg\label{pg:useful_properties_lt_lc} We note five useful
properties of log terminal and log canonical surface
singularities:
\begin{itemize}
\item[(1)] Log terminal is equivalent (for surface
  singularities) to having an isolated quotient singularity
  {\cite[Proposition 4.18]{KM}}.
\item[(2)] Both the log terminal {\cite[Theorem
    2.5]{EsnaultViehweg}} and log canonical conditions
  \cite{Ishii1,Ishii2} are preserved under deformations.
\item[(3)] The log terminal condition is preserved under taking
  quotients {\cite[Proposition 5.20(4)]{KM}}.
\item[(4)] If \(X\) is the double cover of a normal surface
  \(Y\) branched along a divisor \(B\) then \(X\) is log
  canonical if and only if the pair
  \(\left(Y,\frac{1}{2}B\right)\) is log canonical
  {\cite[Proposition 2.5]{AlexeevPardini}}.
\item[(5)] If \((Y,\Delta)\) is log canonical and
  \(p\in \Delta\) is a point then the germ of \((Y,0)\) at \(p\)
  is log terminal {\cite[\S 3.29.4]{KollarSingsofMMP}}.
\end{itemize}

\paragraph{Classification of log canonical surface
  singularities.} Here we list the log canonical surface
singularities together with some properties that we will use
(see {\cite[Theorem 9.6]{Kawamata}} or {\cite[Theorem
  4.7]{KM}}).
\begin{itemize}
\item {\bf Irrational cases:}
  \begin{itemize}
  \item Simple elliptic singularities (not plt). The exceptional
    curve of the minimal resolution is a smooth elliptic curve.
  \item Cusp singularities (not plt). The exceptional locus of
    the minimal resolution is either a nodal elliptic curve or a
    cycle of two or more rational curves.
  \end{itemize}
\item {\bf Rational cases:}
  \begin{itemize}
  \item Isolated quotient singularities (plt), classified by
    Brieskorn \cite{Bries}. The only ones which admit
    \(\QQ\)-Gorenstein smoothings with Milnor number zero are
    the {\em Wahl} singularities, cyclic quotients of the form
    \(\frac{1}{p^2}(1,pq-1)\).
  \item Certain \(\ZZ/2\), \(\ZZ/3\), \(\ZZ/4\) and \(\ZZ/6\)
    quotients of simple elliptic singularities (not plt). There are three
    which admit smoothings of Milnor number zero: in these three
    cases, the dual graph of the exceptional locus of the
    minimal resolution is one of the three shown in Figure
    \ref{fig:simple_elliptic_quot}.
  \item Certain \(\ZZ/2\) quotients of cusp singularities (not
    plt). None of these admit smoothings of Milnor number zero.
  \end{itemize}
\end{itemize}

\begin{figure}[htb]
  \begin{center}
    \begin{tikzpicture}
      \node (a) at (0,0) {\(\bullet\)};
      \node at (a) [above right] {\footnotesize \(-4\)};
      \node (b) at (90:1) {\(\bullet\)};
      \node at (b) [left] {\footnotesize \(-3\)};
      \node (c) at (-30:1.5) {\(\bullet\)};
      \node at (c) [above] {\footnotesize \(-3\)};
      \node (d) at (210:1.5) {\(\bullet\)};
      \node at (d) [above] {\footnotesize \(-3\)};
      \draw (a) -- (b);
      \draw (a) -- (c);
      \draw (a) -- (d);
      \begin{scope}[shift={(4,0)}]
        \node (a2) at (0,0) {\(\bullet\)};
        \node at (a2) [above right] {\footnotesize \(-3\)};
        \node (b2) at (90:1) {\(\bullet\)};
        \node at (b2) [left] {\footnotesize \(-4\)};
        \node (c2) at (-30:1.5) {\(\bullet\)};
        \node at (c2) [above] {\footnotesize \(-2\)};
        \node (d2) at (210:1.5) {\(\bullet\)};
        \node at (d2) [above] {\footnotesize \(-2\)};
        \draw (a2) -- (b2);
        \draw (a2) -- (c2);
        \draw (a2) -- (d2);
      \end{scope}
      \begin{scope}[shift={(8,0)}]
        \node (a3) at (0,0) {\(\bullet\)};
        \node at (a3) [above right] {\footnotesize \(-2\)};
        \node (b3) at (90:1)  {\(\bullet\)};
        \node at (b3) [left] {\footnotesize \(-6\)};
        \node (c3) at (-30:1.5) {\(\bullet\)};
        \node at (c3) [above] {\footnotesize \(-2\)};
        \node (d3) at (210:1.5) {\(\bullet\)};
        \node at (d3) [above] {\footnotesize \(-3\)};
        \draw (a3) -- (b3);
        \draw (a3) -- (c3);
        \draw (a3) -- (d3);
      \end{scope}
    \end{tikzpicture}
    \caption{The dual graphs of the minimal resolutions of the
      elliptic quotient singularities admitting smoothings with Milnor
      number zero.}
    \label{fig:simple_elliptic_quot}
  \end{center}
\end{figure}

\subsection*{Toric and almost toric geometry}

\paragraph{Toric geometry.} We take a moment to recap some basic
toric geometry and then explain how it generalises following the
work of Gross, Hacking, Keel and Kontsevich
\cite{GHK2,GHK1,GHKK} to a wider class of varieties which we
will call {\em almost toric surfaces}. See Cox, Little and
Schenck \cite{CoxLittleSchenck}, Danilov \cite{Danilov} or
Fulton \cite{Fulton} for a proper introduction to toric
geometry. We begin by fixing a lattice \(N\cong \ZZ^n\) and its
dual lattice \(M\). A toric variety \(Y\) is associated to a fan
\(\Sigma\) in \(N_\RR=N\otimes_\ZZ\RR\). Each cone
\(\sigma\in\Sigma\) has a dual cone \(\sigma^\vee\subset M_\RR\)
and the monoid \(\sigma^\vee\cap M\) of integer points in the
dual cone gives us a ring \(\CC[\sigma^\vee\cap M]\) and hence
an affine variety. The variety \(Y\) is obtained by gluing
together these affine pieces using transition functions
determined by the fan. The integer points of \(M\) therefore
have interpretations as local functions on \(Y\). The correct
global interpretation is as follows. Each ray
\(\rho\in\Sigma(1)\) is associated with a torus-invariant
divisor \(D_\rho\subset Y\), and if \(B\subset Y\) is a Weil
\(\QQ\)-divisor \(B\in\left|\sum b_\rho D_\rho\right|\), then we
get a polytope
\begin{equation}\label{eq:toric_tropicalisation}
  P_B\coloneqq\left\{u\in M\otimes_\ZZ\RR\,:\,\langle
    u,\rho\rangle \geq -b_\rho\,\forall\,\rho\in\Sigma(1)\right\}\subset M_\RR.
\end{equation}
The integer points of \(P_B\) correspond to a basis of the
space of sections of the sheaf \(\OO(B)\) (see
{\cite[Proposition 4.3.3, Example
  5.4.5]{CoxLittleSchenck}}). The divisor
\(D=\sum_{\rho\in\Sigma(1)} D_\rho\) is called the {\em toric
  boundary} of \(Y\).

\paragraph{Almost toric geometry.}\label{pg:almost_toric}
Let \(Y\) be a rational surface and \(D\subset Y\) be a cycle of
rational curves supporting an anticanonical divisor. If
\(p\in D\) is an intersection point between two components of
the cycle then the {\em toric blow-up} of \((Y,D)\) at \(p\) is
the pair \((\OP{Bl}_p(Y),\pi^*(D))\). If \(p\in D\) is any other
point then the {\em non-toric blow-up} of \((Y,D)\) at \(p\) is
the pair \((\OP{Bl}_p(Y),\pi^{-1}_*(D))\). In both cases,
\(\pi\colon\OP{Bl}_p(Y)\to Y\) denotes the blow-up map. We will
call a surface \(Y\) an {\em almost toric surface} with {\em
  almost toric boundary \(D\)} if there exists a toric surface
\(\bar{Y}\) with toric boundary \(\bar{D}\) and a sequence of
toric and non-toric blow-ups yielding a pair
\((\tilde{Y},\tilde{D})\) that dominates \((Y,D)\), via a map
\(\tilde{Y}\to Y\) that contracts some subchains of the cycle
\(\tilde{D}\). Gross, Hacking, Keel and Kontsevich show that the
geometry of almost toric surfaces is governed by {\em tropical
  manifolds}, just as toric varieties are governed by
polytopes. A concise exposition of their work which highlights
the parallels with toric geometry can be found in the work of
Mandel \cite{Mandel}; we now summarise the important points.

\paragraph{Tropicalisation.}\label{pg:tropicalisation}
The almost toric analogue of the lattices \(N\) and \(M\) are
the {\em tropicalisations} \(\Vtrop(\ZZ)\) and \(\Utrop(\ZZ)\)
of the log Calabi-Yau surface \(V=Y\setminus D\) and of its
mirror \(U\); each of these is a constellation of integer points
in an integral affine manifold\footnote{Technically this is not
  a manifold: the origin is a singular point.}  (respectively
\(\Vtrop\) and \(\Utrop\)) and there is a canonical ``pairing''
\(\langle,\rangle\colon \Utrop\times\Vtrop\to\RR\) which is
integral when restricted to the integer points. Despite the
forbidding terminology, these integral affine manifolds are
constructed in an elementary way from the cycle \(D\), see
{\cite[Section 1.1]{GHK1}} or {\cite[Sections 2.2,
  4.1]{Mandel}}. The manifold \(\Vtrop\) contains a fan
\(\Sigma\) consisting of rays \(\rho_1,\ldots,\rho_n\)
corresponding to the curves \(D_i\) in the almost toric
boundary. Given a Weil \(\QQ\)-divisor \(B\subset Y\) with
\(B\in\left|\sum b_iD_i\right|\), we get a {\em strongly convex
  polygon} \(\Trop(Y,B)\subset \Utrop\) defined by
\begin{equation}\label{eq:almost_toric_tropicalisation}
  \Trop(Y,B)\coloneqq \left\{u\in \Utrop\,:\,\langle
    u,\rho_i\rangle \geq -b_i\,\forall i\right\}\subset \Utrop,
\end{equation}
see {\cite[Definition 5.14]{Mandel}}. We will call
\(\Trop(Y,B)\) the {\em mirror tropicalisation of the pair
  \((Y,B)\)} and will also adopt this notation for the polygon
\(P_B\) from Equation \eqref{eq:toric_tropicalisation}. If \(B\)
is integral, the integer points \(\Trop(Y,B)\cap \Utrop(\ZZ)\)
correspond with a basis for the space of sections of \(\OO(B)\)
{\cite[Proposition 5.15]{Mandel}}. The basis elements are the
canonical GHK theta functions.

\paragraph{Remark.}\label{pg:rmk_symplectic}
For symplectically-oriented readers, who are used to Symington's
almost toric geometry \cite{Symington}, there is a close
connection. If \(B\) is Poincar\'{e}-dual to the cohomology
class of a symplectic form on \(Y\) then the integral affine
polygon \(\Trop(Y,B)\) is simply the almost toric base diagram
for an almost toric fibration on \(Y\), except with all the
nodes slid along their eigenlines to the barycentre\footnote{One
  cannot, in practice, realise this as an almost toric base
  diagram if the nodes collide in this way, but it still makes
  sense as an integral affine manifold. If you want a precise
  statement, you should think of \(\Utrop\setminus\{(0,0)\}\) as
  the almost toric base diagram for an almost toric fibration on
  the symplectisation of the ideal contact boundary of
  \(Y\setminus D\).}. The rays in \(\Sigma\) are the
eigenrays. We will not use any symplectic almost toric geometry
in this paper, but we will use the terminology of branch cuts,
nodal trades/slides and mutations for manipulating integral
affine manifolds, see for example {\cite[Chapter 8]{EvansBook}}
for a detailed introduction to this subject.

\paragraph{Wahl vertices.}
Let \(\Pi\) be a convex rational polygon and let \(\bm{p}\) be a
vertex of \(\Pi\) where two edges \(e_1,e_2\) meet, ordered so that
\(\Pi\) lies anticlockwise of \(e_1\) and clockwise of \(e_2\). Let
\(\bm{u}_1,\bm{u}_2\) be the primitive integral vectors pointing along
\(e_1,e_2\), oriented away from \(\bm{p}\). We say that \(\bm{p}\) is
a {\em Wahl vertex} if \(\bm{u}_1\wedge \bm{u}_2=c^2\) and
\(\bm{u}_1+\bm{u}_2=c\bm{w}\) for some integer \(c\) and some
primitive integral vector \(\bm{w}\). Here, \(\bm{u}_1\wedge \bm{u}_2\)
is the determinant of the matrix whose columns are
\(\bm{u}_1,\bm{u}_2\). We call \(c\) the {\em index} and \(\bm{w}\)
the {\em eigendirection} at \(\bm{p}\); indeed, \(\bm{w}\) spans the
unique eigenspace of the linear map \(M_{\bm{p}}(\bm{u})\coloneqq
\bm{u}-(\bm{u}\wedge \bm{w})\bm{w}\). The {\em eigenline} at
\(\bm{p}\) is the ray emanating from \(\bm{p}\) in the \(\bm{w}\)
direction.

\begin{figure}[htb]
  \begin{center}
    \begin{tikzpicture}
      \draw[gray,opacity=0.5] (-1,0) grid (4,3);
      \node (a) at (0,0) {\(\bullet\)};
      \node at (a) [below] {\(\bm{p}\)};
      \draw[very thick,->] (a.center) -- (-1,2) node [above right] {\(\bm{u}_2\)};
      \draw[very thick,->] (a.center) -- (4,1) node [above left] {\(\bm{u}_1\)};
      \draw[dotted,->] (a.center) -- (3,3) node [above] {\(\bm{u}_1+\bm{u}_2\)};
      \draw[very thick,dashed,->] (a.center) -- (1,1) node [above] {\(\bm{w}\)};
      \draw (a.center) -- (-2,4) node [above] {\(e_2\)};
      \draw (a.center) -- (5,1.25) node [right] {\(e_1\)};
      \draw[dotted] (-1,2) -- (3,3) -- (4,1);
    \end{tikzpicture}
    \caption{A Wahl vertex with index \(3\) and eigendirection \((1,1)\).}
    \label{fig:wahl_vertex}
  \end{center}
\end{figure}

\paragraph{Mutations.} Given a Wahl vertex \(\bm{p}\) of
\(\Pi\), the eigenline at \(\bm{p}\) separates \(\Pi\) into two
pieces, \(\Pi_1\) and \(\Pi_2\), such that \(\bm{u}_k\) points
into \(\Pi_k\) for \(k=1,2\). We define the {\em clockwise
  (respectively anticlockwise) mutation}\footnote{Mutations
  appear naturally in the context of integral affine manifolds
  and almost toric fibrations. The word
  ``clockwise''/``anticlockwise'' refers to the way we rotate
  the branch cut in an almost toric base diagram to pass between
  the two polygons. See {\cite[Section 8.4]{EvansBook}}} of
\(\Pi\) to be the polygon
\(\overset{\curvearrowright}{\mu}_{\bm{p}}\Pi\coloneqq \Pi_1\cup
M_{\bm{p}}^{-1}(\Pi_2)\) (respectively
\(\overset{\curvearrowleft}{\mu}_{\bm{p}}\Pi=M_{\bm{p}}(\Pi_1)\cup\Pi_2\)). In
all the examples we will consider, \(\Pi\) will be a triangle
contained in the upper half-plane, with a horizontal edge
running along the \(x\)-axis connecting two Wahl vertices:
\(\bm{p}_0\) on the left and \(\bm{p}_1\) on the right. In this
situation, for brevity, we will write
\(\mu_0\coloneqq\overset{\curvearrowright}{\mu}_{\bm{p}_0}\) and
\(\mu_1\coloneqq\overset{\curvearrowleft}{\mu}_{\bm{p}_1}\).
More generally, given a binary sequence
\(\bm{\bin}=(\bin_1,\bin_2,\ldots,\bin_m)\), we will write
\(\mu_{\bm{\bin}}\Pi\coloneqq
\mu_{\bin_m}\mu_{\bin_{m-1}}\cdots\mu_{\bin_1}\Pi\).

\begin{figure}[htb]
  \begin{center}
    \begin{tikzpicture}
      \draw (0,0) -- (4,0) -- (1,3) -- cycle;
      \draw[dashed] (0,0) -- (2,2);
      \draw[->,thick] (0,0) -- (0.5,0.5) node [right] {\(\bm{w}_0\)};
      \node at (0,0) [below] {\(\bm{p}_0\)};
      \node at (4,0) [below] {\(\bm{p}_1\)};
      \node (a) at (0.2,2) {\(\Pi_2\)};
      \node at (2,1) {\(\Pi_1\)};
      \draw (2,2) -- (-2,0) -- (0,0);
      \node (b) at (-1.2,1) {\(M_{\bm{p}_0}^{-1}(\Pi_2)\)};
      \begin{scope}[shift={(7,0)}]
        \draw (0,0) -- (4,0) -- (1,3) -- cycle;
        \draw[dashed] (4,0) -- (2/3,2);
        \draw[->,thick] (4,0) --++ (-1-1/9,2/3) node [below] {\(\bm{w}_1\)};
        \node at (0,0) [below] {\(\bm{p}_0\)};
        \node at (4,0) [below] {\(\bm{p}_1\)};
        \node (a) at (2,2.5) {\(\Pi_1\)};
        \node at (1.5,0.5) {\(\Pi_2\)};
        \node at (6,0.5) {\(M_{\bm{p}_1}(\Pi_1)\)};
        \draw (4,0) -- (6,0) -- (2/3,2);
      \end{scope}
    \end{tikzpicture}
    \caption{The mutations \(\mu_0\) (left) and \(\mu_1\) (right).}
    \label{fig:mutations}
  \end{center}
\end{figure}

\section{Strategy}
\label{sct:strategy}

\pg\label{pg:reduce_to_p2} Suppose that
\(\rho\colon\mathcal{X}\to \Delta\) is a KSBA-stable
degeneration over the disc with fibres
\(X_t\colon=\rho^{-1}(t)\). Suppose that the general fibre
\(X_t\), \(t\neq 0\), is a double branched cover of \(\PP^2\)
and assume that the stable limit \(X_0\) is {\em normal}. By
{\cite[Proposition 2.6]{FPRR}}, after a base-change
\(c\colon\Delta'\to\Delta\), the fibrewise covering involution
\(i_t\colon X_t\to X_t\) defined for \(t\neq 0\) extends over
\(X_0\) to give an involution \(i_0\colon X_0\to X_0\). Write
\(Y_0=X/i_0\) and \(\YY\) for the fibrewise quotient of
\(c^*\XX\). Since \(Y_0\) is a finite quotient of a normal
variety, it is normal, so \(\varrho\colon\YY\to \Delta'\) is a
normal degeneration of \(\PP^2\).

\paragraph{Normal degenerations of
  \(\PP^2\).}\label{pg:normal_deg_of_p2}
Normal degenerations of \(\PP^2\) were studied by B\u{a}descu
\cite{Badescu1,Badescu2} and Manetti
\cite{Manetti1,ManettiThesis} and Hacking and Prokhorov
\cite{HP2,HP1}. We summarise the principal facts we shall
need. If \(Y\) is a normal projective degeneration of \(\PP^2\)
then:
\begin{itemize}
\item[(1)] \(Y\) can have at most one {\em irrational}
  singularity. In this case, the minimal resolution
  \(\tilde{Y}\) is (a blow-up of) an irrational ruled surface of
  irregularity \(q>0\). The exceptional locus of the irrational
  singularity comprises exactly one section together with
  possibly some rational curves from the (blown-up) fibres
  {\cite[Theorem 1, Case (B)]{Badescu1}}.
\item[(2)] If \(y\) is a rational singularity of \(Y\) then the
  germ \((Y,y)\) admits a \(\QQ\)-Gorenstein smoothing with
  Milnor number zero {\cite[IV, Proposition
      2.5]{ManettiThesis}}.
\item[(3)] Any normal degeneration of \(\PP^2\) is a
  \(\QQ\)-Gorenstein smoothing; in particular, \(K^2_Y=9\)
                     {\cite[Corollary 5]{Manetti1}}.
\item[(3)] If \(Y\) has at worst log terminal (quotient)
  singularities then we call \(Y\) a {\em Manetti surface}. The
  Manetti surfaces were completely classified by Hacking and
  Prokhorov (see {\cite[Corollary 1.2]{HP1}} or {\cite[Theorem
    1.1]{HP2}}). Some of the Manetti surfaces are toric: they
  are weighted projective planes \(Y=\PP(a^2,b^2,c^2)\)
  associated to Markov triples \((a,b,c)\). These have (up to)
  three quotient singularities modelled on
  \[\frac{1}{a^2}(b^2,c^2),\quad \frac{1}{b^2}(c^2,a^2)\quad
    \mbox{and}\quad\frac{1}{c^2}(a^2,b^2).\] The remaining
  Manetti surfaces are partial smoothings of the toric Manetti
  surfaces, keeping a subset of these singularities and
  smoothing the rest. They are usually not toric, but they are
  {\em almost toric} in the sense of \ref{pg:almost_toric}. We
  will discuss them in more detail in Section
  \ref{sct:manetti}.
\end{itemize}

\paragraph{Log canonical degenerations of \(\PP^2\).} The
observations in \ref{pg:normal_deg_of_p2} allow us to cut down
the possible {\em log canonical} normal degenerations of
\(\PP^2\):
\begin{itemize}
\item Cusp singularities cannot appear: they are irrational, but
  the exceptional locus contains only rational curves, which is
  incompatible with \ref{pg:normal_deg_of_p2}(1).
\item Of the rational singularities, only Wahl singularities or
  the three elliptic quotients in Figure
  \ref{fig:simple_elliptic_quot} can appear (none of the others
  admit smoothings with Milnor number zero).
\end{itemize}
Manetti {\cite[p.70--71]{ManettiThesis}} shows that it is also
possible to rule out the \([3,3,3;4]\) elliptic quotient case. In the
next lemma, we will show that no elliptic quotient singularities can
occur, which means that the only possibilities are simple elliptic
singularities and Wahl singularities.

\paragraph{Lemma.}\label{pg:elliptic_quotient}
{\em A normal log canonical degeneration \(Y_0\) of \(\PP^2\) cannot admit
  an elliptic quotient singularity.}
\begin{proof}
  Hacking and Prokhorov {\cite[Proposition 3.1]{HP2}} guarantee that
  we can smooth away any of the other singularities, so we only need
  to consider the case where there is precisely one singularity. Since
  they are the only ones admitting smoothings with Milnor number zero,
  we only need to consider the cases in Figure
  \ref{fig:simple_elliptic_quot}. Manetti {\cite[Theorem 4]{Manetti1}}
  shows that \(Y_0\) has Picard rank \(1\) and B\u{a}descu
  {\cite[Theorem 1]{Badescu1}} tells us that the minimal resolution is
  a blow-up of a Hirzebruch surface: the exceptional locus is the
  negative section together with possibly trees of curves contained in
  fibres. In our case, there is a \(3\)-valent exceptional curve
  (which must be the section) and three chains of length \(1\), so
  three fibres \(F_1,F_2,F_3\) are blown up. Suppose that \(F_i\)
  experiences \(n_i\) blow-ups. When we smooth (rationally blow down
  the exceptional locus) we collapse four curves and, since the Picard
  rank of the Hirzebruch surface is \(2\), the Picard rank of the
  smoothing (\(\PP^2\)) is
  \[2+\sum n_i-4=1,\] so \(n_1=n_2=n_3=1\). But if we blow-up a
  fibre (square zero) only once then we only obtain
  \(-1\)-curves, rather than the \(-2,-3,-4,-6\) curves needed
  in the exceptional locus, so these graphs cannot appear.
\end{proof}

\paragraph{Simple elliptic case.}\label{pg:simple_elliptic}
The simple elliptic case does arise as a degeneration of
\(\PP^2\). Suppose \(\PP^2\) is embedded in \(\PP^9\) via its
anticanonical linear system and let \(H\) be a hyperplane which
intersects the image of \(\PP^2\) transversely in a cubic curve
\(C\). Then \(\PP^2\) can be deformed to the cone on
\(C\subset H\). The minimal resolution of this is a ruled
surface of irregularity \(1\) and \(K^2=9\).

\paragraph{Lemma.}\label{lma:milnor_fibre} {\em Suppose that
  \(\YY\to\Delta\) is a normal degeneration of \(\PP^2\) such
  that \(Y_0\) has a simple elliptic singularity (and possibly
  several Wahl singularities). Let \(U_t\subset Y_t\) be the
  Milnor fibre of the simple elliptic singularity embedded in
  \(Y_t\cong\PP^2\). Let \(i_t\colon Y_t\setminus U_t\to Y_t\)
  be the inclusion map. The image of
  \((i_t)_*\colon H_2(Y_t\setminus U_t;\ZZ)\to
  H_2(Y_t;\ZZ)=\ZZ\) is \(3\ZZ\).}
\begin{proof}
  The obstruction to partial smoothing vanishes {\cite[Proposition
      3.1]{HP1}}. Let \(Y'\) be the result of a partial smoothing
  which smooths the Wahl singularities but keeps a simple elliptic
  singularity \(y'\). B\u{a}descu {\cite[Theorem 4.3]{Badescu1}} then
  tells us that \(Y'\) is the cone over an elliptic curve of square
  \(9\). By symplectic parallel transport in the partial smoothing,
  the complement \(Y_t\setminus U_t\) is diffeomorphic to the
  complement \(Y'\setminus\{y'\}\), so its homology is generated by a
  class \(C\) of square \(9\). Therefore \((i_t)_*C=3H\in
  H_2(\PP^2;\ZZ)\).
\end{proof}

\paragraph{Corollary.}\label{cor:branch_curve_through_point}
{\em Suppose that \(\YY\to \Delta\) is a degeneration of
  \(\PP^2\) such that \(Y_0\) has a simple elliptic
  singularity. Suppose that \(\BB\) is a divisor in \(\YY\)
  which is flat over \(\Delta\). Write
  \(B_t\coloneqq\BB\cap Y_t\); by flatness, the degree
  \([B_t]\in H_2(Y_t;\ZZ)\) is constant in \(t\) {\cite[III,
    Theorem 9.9]{Hartshorne}}. If \([B_t]\) is not divisible by
  \(3\) then \(B_0\) must pass through the simple elliptic
  singularity.}
\begin{proof}
  Let \(U_0\) be a small closed Euclidean neighbourhood of the
  simple elliptic singularity \(y\in Y_0\) and let
  \(U_t\subset Y_t\) be the Milnor fibre which maps to \(U_0\)
  under symplectic parallel transport.\footnote{We are
    transporting \(Y_t\) here, not \(B_t\).} If \(B_t\) is not
  divisible by \(3\) then \(B_t\) must intersect \(U_t\) by
  Lemma \ref{lma:milnor_fibre}. Therefore \(B_0\) intersects
  \(U_0\). Since \(U_0\) can be taken arbitrarily small, we see
  that \(B_0\) must pass through \(y\) itself.
\end{proof}

\paragraph{Corollary.} {\em Suppose that \(\XX\to\Delta\) is a
  normal KSBA degeneration of octic double planes and let
  \(\YY\to\Delta\) be the corresponding degeneration of
  \(\PP^2\) coming from \ref{pg:reduce_to_p2}. Then \(Y_0\)
  cannot have a simple elliptic singularity.}
\begin{proof}
  First note that the branch locus \(\BB\) of \(\XX\to\YY\) is
  flat over \(\Delta\); this follows from {\cite[III,
    Proposition 9.7]{Hartshorne}} since normality of \(\XX\)
  implies that \(\BB\) is reduced. If \(Y_0\) had a simple
  elliptic singularity then the branch curve \(B_0=\BB\cap Y_0\)
  would need to pass through it by Corollary
  \ref{cor:branch_curve_through_point} (because \(3\) does not
  divide \(8\)). But since a simple elliptic singularity is not
  log terminal, the pair \((Y_0,\frac{1}{2}B_0)\) fails to be
  log canonical by \ref{pg:useful_properties_lt_lc}(5), and so
  by \ref{pg:useful_properties_lt_lc}(4), the double cover
  \(X_0\) fails to be log canonical.
\end{proof}

\pg The conclusion of our arguments so far is that, in the
context of \ref{pg:reduce_to_p2}, a normal KSBA degeneration of
an octic double plane is a double branched cover \(X_0\to Y_0\)
of a Manetti surface. Note that it is possible for the branch
curve \(B_0\) to develop singularities away from the singular
locus of the Manetti surface; the allowable singularities were
classified by Anthes {\cite[Table 2]{Anthes}} and yield
Gorenstein singularities of the double cover. We will ignore
these in what follows, and focus on the non-Gorenstein
singularities which appear when \(B_0\) passes through the
singularities of \(Y_0\). Our final goal in this section is to
better understand the branch locus \(B_0\). This will be
achieved in Lemma \ref{lma:octic_branch} below, after some
preparation.

\pg\label{pg:octic_topology} In the setting of
\ref{pg:reduce_to_p2}, let \(\mathcal{Y}\to \Delta\) be the
smoothing of the Manetti surface \(Y_0\) over the disc obtained
by quotienting \(\XX\) by its fibrewise involution, so that the
general fibre \(Y_t\) is \(\PP^2\). For each
singularity\footnote{Each singularity like
  \(\frac{1}{a^2}(b^2,c^2)\) can be put into this form with
  \(p=a\) and \(q=3c/b\mod a\) because then
  \(c^2/b^2=pq-1\mod a^2\).} \(\frac{1}{p_i^2}(1,p_iq_i-1)\),
\(i=1,\ldots,m\), of \(Y_0\), we find an embedded copy of the
Milnor fibre \(B_{p_i,q_i}\) inside \(\PP^2\) (note that
\(m\leq 3\)). We can assume that the \(B_{p_i,q_i}\) for
different singularities are pairwise disjoint since the
singularities are distinct. We have
\(H_2(Y_0;\ZZ)\cong H_2\left(\PP^2,\bigcup
  B_{p_i,q_i};\ZZ\right)\). The long exact sequence for this
relative homology group yields
\[\cdots\to H_2(\PP^2;\ZZ)\to H_2(Y_0;\ZZ)\to
  H_1\left(\bigcup B_{p_i,q_i};\ZZ\right)\to H_1(\PP^2;\ZZ)=0.\]
We have
\(H_1\left(\bigcup B_{p_i,q_i};\ZZ\right)=\prod \ZZ/p_i\) and \(H_*\left(\bigcup B_{p_i,q_i};\ZZ\right)=0\) for \(*>1\), so
the exact sequence becomes\footnote{Note that in a Manetti
  surface, the \(p_i\) form part of a Markov triple, and hence
  are pairwise coprime, which is why
  \(\prod\ZZ/p_i\cong\ZZ/\left(\prod p_i\right)\).}
\[0\to \ZZ\to \ZZ\to \ZZ/\left(\prod p_i\right)\to 0.\]
Therefore the map \(\ZZ=H_2(\PP^2;\ZZ)\to H_2(Y_0;\ZZ)=\ZZ\)
takes an octic curve \([B]=8\in \ZZ= H_2(\PP^2;\ZZ)\) to a curve
in the class \(8\prod p_i\in \ZZ=H_2(Y_0;\ZZ)\).

\paragraph{Definition.}\label{dfn:octic}
We define a Weil divisor on \(Y_0\) to be an {\em octic} if it
lives in the homology class
\(8\prod_{i=1}^m p_i\in\ZZ=H_2(Y_0;\ZZ)\). For example, if
\(Y=\PP(a^2,b^2,c^2)\) then we take a divisor in the linear
system \(|\OO(8abc)|\), i.e. a divisor cut out by a polynomial
of weighted degree \(8abc\) in the weighted homogeneous
coordinates.

\paragraph{Lemma.}\label{lma:octic_branch} {\em In the context
  of \ref{pg:reduce_to_p2}, the branch curve of the double cover
  \(X_0\to Y_0\) is an octic curve in \(Y_0\) in the sense of
  Definition \ref{dfn:octic}.}
\begin{proof}
  Since the degeneration is KSBA-stable the central fibre and the
  general fibre have the same value of the characteristic number
  \(K^2\). For \(t\in\Delta\), let \(f_t\colon X_t\to Y_t\) be the
  double cover with branch locus \(B_t\subset Y_t\) and ramification
  locus \(R_t\subset X_t\). Then
  \(K_{X_t}=f_t^*(K_{Y_t}+B_t/2)\). When \(t\neq 0\), we have
  \(Y_t=\PP^2\), \(K_{Y_t}=-3h\) and \([B_t]=8h\) where \(h\in
  H_2(Y_t;\ZZ)\) is a generator. So \(K_{X_t}^2 = f_t^*(h^2)=2\)
  because \(h^2\) is a single point and \(f_t\) is a double
  cover. When \(t=0\) the generator \(h_0\in H_2(Y_0;\ZZ)\) has
  self-intersection \(\prod_{i=1}^mp^{-2}_i\) and
  \(K_{Y_0}=-3h_0\prod_{i=1}^mp_i\). Write \(B=\beta h_0\) for some
  \(\beta\in\ZZ\). The only way to achieve \(K_{X_0}^2=2\) is to take
  \(\beta=8\prod_{i=1}^mp_i\), which implies that \(B\) is an octic in
  the sense of Definition \ref{dfn:octic}.
\end{proof}

\paragraph{Octic double Manetti
  surfaces.}\label{pg:octic_double_manetti}
To summarise, we see that a normal projective surface arising as
a KSBA-stable limit of octic double planes is an {\em octic
  double Manetti surface}, that is a double cover of a Manetti
surface branched over an octic curve in the sense of Definition
\ref{dfn:octic}. The rest of the paper will be dedicated to
understanding which octic double Manetti surfaces are normal
with at worst quotient singularities.

\section{Toric and almost toric Manetti surfaces}
\label{sct:manetti}

\pg In this section, we review the integral affine geometry of
the mirror tropicalisations \(\Trop(Y,B)\) where \(Y\) is a
Manetti surface and \(B\) is an octic Weil divisor.

\paragraph{Markov triples.} A triple of positive integers
\(a,b,c\) is a {\em Markov triple} if
\[a^2+b^2+c^2=3abc.\] These triples can all be obtained from
\(1,1,1\) by a sequence of {\em mutations}, in which \(a,b,c\)
is replaced by \(a,b,c'=3ab-c\). Drawing a tree whose vertices
are (unordered) Markov triples and whose edges correspond to
mutations allows us to organise the triples in the {\em Markov
  topograph} (see Figure \ref{fig:markov_tree}): each region
into which the tree separates the plane is labelled by a Markov
number which appears in every triple around its boundary; each
edge corresponding to a mutation \(a,b,c\leadsto a,b,c'\) is
labelled with the pair \((a,b)\).

\begin{figure}[htb]
  \begin{center}
    \begin{tikzpicture}[scale=0.7, every node/.style={scale=0.7}]
      \node (a) at (0,0) {\((1,1,1)\)};
      \node (b) at (0,2) {\((1,1,2)\)};
      \node (c) at (0,4) {\((1,2,5)\)};
      \node[shape=ellipse,draw=black] (d) at (-4,5) {\((1,5,13)\)};
      \node[fill=lightgray] (e) at (4,5) {\((2,5,29)\)};
      \node[shape=ellipse,draw=black] (f) at (-6,3) {\((1,13,34)\)};
      \node[draw=black] (g) at (-5,7) {\((5,13,194)\)};
      \node[fill=lightgray] (h) at (5,7) {\((5,29,433)\)};
      \node[fill=lightgray] (i) at (6,3) {\((2,29,169)\)};
      \node[shape=ellipse,draw=black] (j) at (-8,1) {\((1,34,89)\)};
      \node[fill=lightgray] (k) at (8,1) {\((2,169,985)\)};
      \node[draw=black] (l) at (-8,5) {\((13,34,1325)\)};
      \node[fill=lightgray] (m) at (8,5) {\((29,169,14701)\)};
      \draw[->] (a) -- (b) node[midway,sloped,above] {\((1,1)\)};
      \draw[->] (b) -- (c) node[midway,sloped,above] {\((1,2)\)};
      \draw[->] (c) -- (d) node[midway,sloped,below] {\((1,5)\)};
      \draw[->] (c) -- (e) node[midway,sloped,below] {\((2,5)\)};
      \draw[->] (d) -- (f) node[midway,sloped,below] {\((1,13)\)};
      \draw[->] (d) -- (g) node[midway,sloped,below] {\((5,13)\)};
      \draw[->] (e) -- (h) node[midway,sloped,above] {\((2,29)\)};
      \draw[->] (e) -- (i) node[midway,sloped,below] {\((5,29)\)};
      \draw[->] (f) -- (j) node[midway,sloped,below] {\((1,34)\)};
      \draw[->] (i) -- (k) node[midway,sloped,below] {\((2,169)\)};
      \draw[->] (i) -- (m) node[midway,sloped,below] {\small \((29,169)\)};
      \draw[->] (f) -- (l) node[midway,sloped,below] {\((13,34)\)};
      \node at (-3,3) {\(\bm{1}\)};
      \node at (3,3) {\(\bm{2}\)};
      \node at (0,6) {\(\bm{5}\)};
      \node at (-6,5) {\(13\)};
      \node at (6,5) {\(29\)};
      \node at (-8,2.5) {\(34\)};
      \node at (8,2.5) {\(169\)};
      \node at (-5,7.7) {\(194\)};
      \node at (-8,5.7) {\(1325\)};
      \node at (5,7.7) {\(433\)};
      \node at (8,5.7) {\(14701\)};
      \draw[dotted] (g) --++ (45:1);
      \draw[dotted] (g) --++ (135:1);
      \draw[dotted] (h) --++ (45:1);
      \draw[dotted] (h) --++ (135:1);
      \draw[dotted] (l) --++ (45:1);
      \draw[dotted] (l) --++ (135:1);
      \draw[dotted] (m) --++ (45:1);
      \draw[dotted] (m) --++ (135:1);
      \draw[dotted] (j) --++ (135:1);
      \draw[dotted] (j) --++ (-135:1);
      \draw[dotted] (k) --++ (45:1);
      \draw[dotted] (k) --++ (-45:1);
    \end{tikzpicture}
    \caption{The Markov topograph. Shapes and shading of the vertices will be explained in \ref{dfn:Z_polygons}--\ref{pg:integer_points}. In the end, we will see that the only normal KSBA-stable octic double Manetti surfaces are double covers of \(\HP(1)\), \(\HP(2)\) and \(\HP(5)\), shown in bold here.}
    \label{fig:markov_tree}
  \end{center}
\end{figure}

The Markov topograph is a ``map'' of the set of Manetti surfaces
in the following sense. To every vertex of the Markov topograph
there is an associated toric Manetti surface, the weighted
projective space \(\PP(a^2,b^2,c^2)\). To the (unique) edge
labelled \((a,b)\) there is a non-toric Manetti surface which we
will call \(\HP(a,b)\), obtained from \(\PP(a^2,b^2,c^2)\) by
smoothing the \(\frac{1}{c^2}(a^2,b^2)\) singularity (leaving
the singularities \(\frac{1}{b^2}(c^2,a^2)\) and
\(\frac{1}{a^2}(b^2,c^2)\)). To every region labelled \(c\)
there is a non-toric Manetti surface\footnote{Our notation is a
  little sloppy: it is a longstanding conjecture \cite{Aigner}
  that each Markov number labels only one region. In the end,
  the only non-toric Manetti surface that will arise for us will
  be \(\HP(5)\), and there is a unique region labelled \(5\).}
\(\HP(c)\) obtained from \(\PP(a^2,b^2,c^2)\) by smoothing all
singularities except \(\frac{1}{c^2}(a^2,b^2)\). There are some
redundancies here: \(\HP(1)=\HP(1,1)=\PP^2\),
\(\HP(2)=\PP(1,1,4)\), and \(\HP(a,b)=\PP(a^2,b^2,1)\) whenever
\((a,b,1)\) is a Markov triple.

\paragraph{Toric Manetti surfaces and Markov
  triangles.}\label{pg:toric_manetti_surfaces}
Let \(Y\) be the toric Manetti surface \(\PP(a^2,b^2,c^2)\) and
\(B\subset Y\) be an octic Weil divisor, that is
\(\OO(B)\cong\OO(8abc)\), where \(\OO(8abc)\) is the sheaf whose
sections are weighted homogeneous polynomials of total degree
\(8abc\). The moment triangle\footnote{i.e. \(\Trop(Y,B)\).}
\(\Pi(a,b,c)\) can be constructed inductively as follows.

Let \(\Pi(1,1,1)\) be the triangle with vertices at
\(\bm{p}_0=(0,0)\), \(\bm{p}_1=(8,0)\) and
\(\bm{p}_2=(0,8)\). All three vertices are Wahl (with index
\(1\)). Define \(\Pi(1,1,2)=\mu_1\Pi(1,1,1)\) and
\(\Pi(1,2,5)=\mu_1\Pi(1,1,2)\). Given a binary sequence
\(\bm{\bin}\), we obtain a polygon
\(\Pi(a,b,c)\coloneqq \mu_{\bm{\bin}}\Pi(1,2,5)\). Here, the
label \((a,b,c)\) is a Markov triple constructed as follows:
start with \((1,2,5)\) and perform a sequence of mutations to
the triple: if \(\bin_i=0\), the \(i\)th mutation replaces \(a\)
by \(3bc-a\); if \(\bin_i=1\), the \(i\)th mutation replaces
\(b\) by \(3ac-b\). We call \(\Pi(a,b,c)\) the {\em Markov
  triangle} associated with the Markov triple \((a,b,c)\); see
Figure \ref{fig:markov_triangle} for an indicative sketch of a
general \(\Pi(a,b,c)\), and Figure \ref{fig:Z_points_1} later
for some concrete examples. We continue to write \(\bm{p}_i\)
for the vertices; these vertices are all Wahl type, with
eigendirections \(\bm{w}_i\) and eigenlines \(W_i\). Note that
the eigenlines all meet at the barycentre \(\bm{c}=(8/3,8/3)\)
for \(\Pi(1,1,1)\), and since mutations fix eigenlines, the
eigenlines must all meet at \(\bm{c}\) for any
\(\Pi(a,b,c)\). Label the edges of \(\Pi(a,b,c)\) as \(e_c\)
(for the horizontal edge), \(e_b\) and \(e_a\) (for the edges
emanating from the lefthand and righthand vertices \(\bm{p}_0\)
and \(\bm{p}_1\)). Writing \(|e|\) for the affine length of
\(e\), we have
\[|e_a|=\frac{8a}{bc},\quad |e_b|=\frac{8b}{ac},\quad
  |e_c|=\frac{8c}{ab}.\] Let \(I_0=\{i\,:\,\bin_i=0\}\) and
\(I_1=\{i\,:\,\bin_i=1\}\), and let
\(\Pi(a_i,b_i,c_i)=\mu_{(\bin_1,\ldots,\bin_i)}\Pi(1,2,5)\). By
construction, we have \(\bm{p}_0=(-\alpha,0)\) and
\(\bm{p}_1=(20+\beta,0)\), where
\begin{equation}\label{eq:affine_lengths}\alpha = \sum_{i\in
    I_0}\frac{8b_i}{a_ic_i}\qquad \beta = \sum_{i\in
    I_1}\frac{8a_i}{b_ic_i}.\end{equation} Note that \(a,b\leq c\) but
it depends on the triple as to whether \(a\leq b\) or \(b\leq a\).

\begin{figure}[htb]
  \centering
  \begin{tikzpicture}
    \filldraw[fill=lightgray,opacity=0.5,draw=none] (0,0) -- (10,0) -- (4,2) -- cycle;
    \node at (4,1) {\(\OP{inn}(\Pi)\)};
    \draw (0,0) -- (10,0) node[midway,below] {\(e_c\)} -- (3,4)
    node [midway,above] {\(e_a\)} -- (0,0) node [midway,above left] {\(e_b\)};
    \draw[dotted] (0,0) -- (4,2) node [midway,above] {\(W_0\)};
    \draw[dotted] (10,0) -- (4,2) node [midway,below] {\(W_1\)};
    \draw [dotted] (3,4) -- (4,2) node [midway,right] {\(W_2\)};
    \node at (0,0) [left] {\(\bm{p}_0\)};
    \node at (10,0) [right] {\(\bm{p}_1\)};
    \node at (3,4) [left] {\(\bm{p}_2\)};
    \draw[->] (0,0) --++ (1,0.5) node[above] {\(\bm{w}_0\)};
    \draw[->] (10,0) --++ (-1,1/3) node[above left] {\(\bm{w}_1\)};
    \draw[->] (3,4) --++ (0.5,-1) node[left] {\(\bm{w}_2\)};
    \node at (4,2) {\(\times\)};
    \node at (4,2) [above right] {\(\bm{c}\)};
  \end{tikzpicture}
  \caption{A Markov triangle \(\Pi\) (not to scale), showing the
    vertices \(\bm{p}_i\), the eigendirections \(\bm{w}_i\),
    eigenlines \(W_i\), and the barycentre
    \(\bm{c}=(8/3,8/3)\). We have also shaded the subtriangle
    \(\OP{inn}(\Pi)\) discussed in the proof of Lemma
    \ref{pg:integer_points}.}
  \label{fig:markov_triangle}
\end{figure}

\paragraph{Almost toric Manetti
  surfaces.}\label{pg:almost_toric_manetti} Fix a toric Manetti
surface \(Y=\PP(a^2,b^2,c^2)\). Pick a proper subset
\(S\subset\{a,b,c\}\) of the Markov triple and write \(\HP(S)\) for
the Manetti surface obtained from \(\PP(a^2,b^2,c^2)\) by smoothing
precisely the subset \(\{p_m,\:\,m\not\in S\}\) of singularities. If
you look carefully at the proofs in {\cite[Section 8]{HP1}} or
{\cite[Section 6]{HP2}}, you will see that Hacking and Prokhorov give
an explicit description of \(\HP(S)\): it is obtained by taking a
(toric) Hirzebruch surface, performing some toric blow-ups, performing
one or two non-toric blow-ups (blowing up points on the interior of
the toric boundary), and then contracting \(|S|\) Wahl chains in the
strict transform of the toric boundary. In other words, \(\HP(S)\) is
{\em almost toric} in the sense of \ref{pg:almost_toric}. In Appendix
\ref{sct:almost_toric} we include an example to aid the unversed
reader in figuring out exactly how to carry this procedure out in a
non-trivial example, and connect this with Symington's almost toric
geometry. The number of non-toric blow-ups required is equal to
\(3-|S|\), and they occur on different components of the toric
boundary. In particular, since all the \(T\)-singularities of
\(\PP(a^2,b^2,c^2)\) have Milnor number zero, we never need to make
multiple non-toric blow-ups on the same boundary component.

\paragraph{The mirror tropicalisation.} Let \(Y\) be an almost
toric Manetti surface and \(B\) an octic Weil divisor. We now
describe the mirror tropicalisation \(\Trop(Y,B)\).
\begin{itemize}
\item For \(Y=\HP(a,b)\), we start with the triangle \(\Pi(a,b,c)\)
  perform a single {\em nodal trade} at the corner \(\bm{p}_2\)
  corresponding to the \(\frac{1}{c^2}(a^2,b^2)\) singularity, and
  slide the node to the barycentre. In other words, we make a branch
  cut along the eigenline \(W_2\) connecting \(\bm{p}_2\) to
  \((8/3,8/3)\) and twist the integral affine structure by the
  monodromy \(M_{\bm{p}_2}\) when we cross this branch cut in a
  clockwise direction.
\item For \(\HP(c)\), we start with the moment triangle of
  \(\PP(a^2,b^2,c^2)\), perform two nodal trades at the corners
  \(\bm{p}_0\) and \(\bm{p}_1\) corresponding to the
  \(\frac{1}{a^2}(b^2,c^2)\) and \(\frac{1}{b^2}(c^2,a^2)\)
  singularities, and slide both nodes to the barycentre. In other
  words, we make branch cuts connecting \(\bm{p}_0\) and \(\bm{p}_1\)
  to \((8/3,8/3)\) and twist the integral affine structure by the
  monodromy \(M_{\bm{p}_0}\) or \(M_{\bm{p}_1}\) when we cross the
  corresponding branch cut clockwise.
\end{itemize}
Note that if any of the ``singularities'' in \(S\) are actually
smooth points (e.g. \(a=1\)), we omit these nodal trades. In
each of the last two cases, the full tropicalisation \(\Utrop\)
is obtained from the resulting picture by erasing the edges of
\(\Trop(Y,B)\) and extending the branch-cuts out to infinity.

We write \(\Pi(S)\) for the mirror tropicalisation of
\((\HP(S),B)\). Figure \ref{fig:tropicalisations} shows the mirror
tropicalisations for two of the simplest non-toric Manetti surfaces
\(\Pi(5)\) and \(\Pi(29)\).

\begin{figure}[htb]
  \centering
  \begin{tikzpicture}[scale=0.5]
    \begin{scope}[shift={(0,-18)}]
      \node at (0,0) {\(\circ\)};
      \foreach \x in {1,2,...,20} {\node at (\x,0) {\(\bullet\)};}
      \foreach \x in {0,1,...,13} {\node at (\x,1) {\(\cdot\)};}
      \foreach \x in {0,1,...,7} {\node at (\x,2) {\(\cdot\)};}
      \foreach \x in {1} {\node at (\x,3) {\(\cdot\)};}
      \node at (0,3) {\(\bullet\)};
      \foreach \y in {1,2,3} {\node at (0,\y) {\(\bullet\)};}
      \draw (0,0) -- (20,0) node[below] {\((20,0)\)} -- (0,16/5) node [left] {\((0,16/5)\)} -- cycle;
      \draw[dashed] (20,0) -- (8/3,8/3) node [midway,above,sloped] {\tiny branch cut \((13,-2)\)};
      \node at (8/3,8/3) {\(\times\)};
      \node at (-2,1) {\(\Pi(5)\)};
    \end{scope}

    \begin{scope}[shift={(0,-24)}]
      \foreach \x in {-3,-2,-1} {\node at (\x,0) {\(\bullet\)};}
      \node at (0,0) {\(\circ\)};
      \foreach \x in {1,2,...,20} {\node at (\x,0) {\(\bullet\)};}
      \foreach \x in {-1,0,...,13} {\node at (\x,1) {\(\cdot\)};}
      \foreach \x in {2,3,...,7} {\node at (\x,2) {\(\cdot\)};}
      \draw (-16/5,0) -- (20,0) node[below] {\((20,0)\)} -- (80/29,80/29) node [above] {\((80/29,80/29)\)} -- (-16/5,0) node [below] {\((-16/5,0)\)};
      \draw[dashed] (20,0) -- (8/3,8/3) node [midway,above,sloped] {\tiny branch cut \((13,-2)\)};
      \draw[dashed] (-16/5,0) -- (8/3,8/3) node [midway,above,sloped] {\tiny branch cut \((11,5)\)};
      \node at (8/3,8/3) {\(\times\)};
      \node at (-2,3) {\(\Pi(29)\)};
    \end{scope}
  \end{tikzpicture}
  \caption{The integral affine manifolds \(\Pi(5)\)
    and \(\Pi(29)\). The branch cuts always meet at the point
    \((8/3,8/3)\) (marked \(\times\)) since this point is
    invariant under mutation. The origin is marked \(\circ\),
    integer points are marked \(\cdot\), and integer points on
    the boundary are marked \(\bullet\).}
  \label{fig:tropicalisations}
\end{figure}

\section{Integer points in mirror tropicalisations}
\label{sct:integer_points}

\pg Let \(Y\) be a Manetti surface and \(B\) an octic Weil
divisor. Other octic divisors are given by the vanishing of a
section of the sheaf \(\OO(B)\), and a basis for this space of
sections is given by the GHK theta functions, which are in
bijection with integer points of \(\Trop(Y,B)\). In this
section, we will enumerate the integer points of \(\Trop(Y,B)\)
in all cases. For \(\Pi(1,1,1)\) and \(\Pi(1,1,2)\), it is easy
to read off the integer points from a picture (see Figure
\ref{fig:Z_points_1}).

\begin{figure}[htb]
  \centering
  \begin{tikzpicture}[scale=0.4]
    \node at (0,0) {\(\circ\)};
    \foreach \x in {1,2,...,8} {\node at (\x,0) {\(\bullet\)};}
    \foreach \x in {0,1,...,7} {\node at (\x,1) {\(\cdot\)};}
    \foreach \x in {0,1,...,6} {\node at (\x,2) {\(\cdot\)};}
    \foreach \x in {0,1,...,5} {\node at (\x,3) {\(\cdot\)};}
    \foreach \x in {0,1,...,4} {\node at (\x,4) {\(\cdot\)};}
    \foreach \x in {0,1,...,3} {\node at (\x,5) {\(\cdot\)};}
    \foreach \x in {0,1,2} {\node at (\x,6) {\(\cdot\)};}
    \foreach \x in {0,1} {\node at (\x,7) {\(\bullet\)};}
    \foreach \x in {0} {\node at (\x,8) {\(\bullet\)};}
    \foreach \y in {1,2,...,8} {\node at (0,\y) {\(\bullet\)};}
    \foreach \y in {1,2,...,8} {\node at (8-\y,\y) {\(\bullet\)};}
    \draw (0,0) -- (8,0) node[right] {\footnotesize \((8,0)\)} -- (0,8) node [left] {\footnotesize \((0,8)\)} -- cycle;
    \node at (7,5) {\(\Pi(1,1,1)\)};
    \begin{scope}[shift={(14,0)}]
      \node at (0,0) {\(\circ\)};
      \foreach \x in {1,2,...,16} {\node at (\x,0) {\(\bullet\)};}
      \foreach \x in {0,1,...,12} {\node at (\x,1) {\(\cdot\)};}
      \foreach \x in {0,1,...,8} {\node at (\x,2) {\(\cdot\)};}
      \foreach \x in {0,1,...,4} {\node at (\x,3) {\(\cdot\)};}
      \foreach \y in {1,...,4} {\node at (0,\y) {\(\bullet\)};}
      \foreach \y in {0,1,...,4} {\node at ({16-4*\y},\y) {\(\bullet\)};}
      \draw (0,0) -- (16,0) node[below right] {\footnotesize \((16,0)\)} -- (0,4) node [left] {\footnotesize \((0,4)\)} -- cycle;
      \node at (11,3) {\(\Pi(1,1,2)\)};
    \end{scope}
  \end{tikzpicture}
  \caption{The integer points in the polygons \(\Pi(1,1,1)\),
    \(\Pi(1,1,2)\).}
  \label{fig:Z_points_1}
\end{figure}

\paragraph{Definition.}\label{dfn:Z_polygons} 
(See Figure \ref{fig:Z_polygons}). Let
\begin{itemize}
\item \(\Pi_5\) be the quadrilateral with vertices at \((0,0)\),
  \((0,3)\), \((1,3)\), \((20,0)\),
\item \(\Pi_A\) be the quadrilateral with vertices at
  \((0,0)\), \((0,3)\), \((14,1)\), and \((20,0)\),
\item \(\Pi_B\) be the pentagon with vertices at \((-3,0)\),
  \((2,2)\), \((7,2)\), \((14,1)\), and \((20,0)\),
\item \(\Pi_C\) be the pentagon with vertices at \((-3,0)\),
  \((-1,1)\), \((2,2)\), \((7,2)\), \((20,0)\).
\end{itemize}

\begin{figure}[htb]
  \centering
  \begin{tikzpicture}[scale=0.5]
    \begin{scope}
      \clip (0,0) -- (0,3) -- (1,3) -- (20,0) -- cycle;
      \foreach \x in {0,1,...,20}{
        \foreach \y in {0,1,2,3} {\node at (\x,\y) {\(\cdot\)};}}
    \end{scope}
    \node at (-1,1.5) {\(\Pi_5\)};
    \draw (0,0) -- (0,3) -- (1,3) -- (20,0) -- cycle;
    \foreach \x in {0,1,...,20} {\node at (\x,0) {\(\bullet\)};}
    \foreach \y in {0,1,2,3} {\node at (0,\y) {\(\bullet\)};}
    \node at (1,3) {\(\bullet\)};
    \begin{scope}[shift={(0,-4.5)}]
      \begin{scope}
        \clip (0,0) -- (0,3) -- (14,1) -- (20,0) -- cycle;
        \foreach \x in {0,1,...,20}{
          \foreach \y in {0,1,2,3} {\node at (\x,\y) {\(\cdot\)};}}
      \end{scope}
      \node at (-1,2) {\(\Pi_A\)};
      \draw (0,0) -- (0,3) -- (14,1) -- (20,0) -- cycle;
      \foreach \x in {0,1,...,20} {\node at (\x,0) {\(\bullet\)};}
      \foreach \y in {0,1,2,3} {\node at (0,\y) {\(\bullet\)};}
      \node at (14,1) {\(\bullet\)};
      \node at (7,2) {\(\bullet\)};
    \end{scope}
    \begin{scope}[shift={(0,-9)}]
      \begin{scope}
        \clip (-3,0) -- (2,2) -- (7,2) -- (14,1) -- (20,0) -- cycle;
        \foreach \x in {-3,-2,...,20}{
          \foreach \y in {0,1,2} {\node at (\x,\y) {\(\cdot\)};}}
      \end{scope}
      \node at (-1,2) {\(\Pi_B\)};
      \draw (-3,0) -- (2,2) -- (7,2) -- (14,1) -- (20,0) -- cycle;
      \foreach \x in {-3,-2,...,20} {\node at (\x,0) {\(\bullet\)};}
      \foreach \x in {2,3,...,7} {\node at (\x,2) {\(\bullet\)};}
      \node at (14,1) {\(\bullet\)};
    \end{scope}
    \begin{scope}[shift={(0,-13.5)}]
      \node at (-1,2) {\(\Pi_C\)};
      \begin{scope}
        \clip (-3,0) -- (-1,1) -- (2,2) -- (7,2) -- (20,0) -- cycle;
        \foreach \x in {-3,-2,...,20}{
          \foreach \y in {0,1,2} {\node at (\x,\y) {\(\cdot\)};}}
      \end{scope}
      \draw (-3,0) -- (-1,1) -- (2,2) -- (7,2) -- (20,0) -- cycle;
      \foreach \x in {-3,-2,...,20} {\node at (\x,0) {\(\bullet\)};}
      \foreach \x in {2,3,...,7} {\node at (\x,2) {\(\bullet\)};}
      \node at (-1,1) {\(\bullet\)};
    \end{scope}
  \end{tikzpicture}
  \caption{The polygons \(\Pi_5\), \(\Pi_A\),
    \(\Pi_B\), and \(\Pi_C\).}
  \label{fig:Z_polygons}
\end{figure}

\paragraph{Lemma.}\label{pg:integer_points}
{\em Write \(\ZZ\Pi\) for the integer points of \(\Pi\).
  Let \(\Pi=\mu_{\bm{\bin}}\Pi(1,2,5)\).
  \begin{itemize}
  \item \(\ZZ\Pi(1,2,5)=\ZZ\Pi_5\).
  \item If \(\bm{\bin}=(1,\ldots,1)\) of length \(n\geq 3\) then
    \(\ZZ\Pi=\ZZ\Pi_A\).
  \item If \(\bin_1=1\) but at least one \(\bin_i=0\) then \(\ZZ\Pi=\ZZ\Pi_B\).
  \item If \(\bin_1=0\) then \(\ZZ\Pi=\ZZ\Pi_C\).
  \end{itemize}
  In each case, we say that \(\Pi\) is respectively of Type A, B
  or C. In Figure \ref{fig:markov_tree}, triples of Type A are circled, triples of Type B are boxed and triples of Type C are shaded grey.}

\begin{proof}
  First, it is easy to check that \(\Pi(1,2,5)\) has vertices at
  \((0,0)\), \((20,0)\) and \((0,16/5)\) and from this that
  \(\ZZ\Pi(1,2,5)=\ZZ\Pi_5\).

  Let \(\OP{inn}(\Pi)\) be the {\em inner triangle} bounded by the
  horizontal edge of \(\Pi\) and the two eigenlines \(W_0\) and
  \(W_1\), see Figure \ref{fig:markov_triangle}. This has vertices
  at \(\bm{p}_0=(-\alpha,0)\), \(\bm{p}_1=(20+\beta,0)\) and
  \(\bm{c}=(8/3,8/3)\) (see Equation \eqref{eq:affine_lengths} for
  \(\alpha\) and \(\beta\)). Note that \(\OP{inn}(\Pi)\) is
  precisely the region which is fixed by both of the mutations
  \(\mu_0\) and \(\mu_1\).

  For \(\Pi(5,2,29)=\mu_0\Pi(1,2,5)\), all integral points are
  contained in \(\OP{inn}(\Pi(5,2,29))\), and since
  \(\OP{inn}(\Pi)\) is unaffected by either \(\mu_0\) or \(\mu_1\)
  mutations, we have established the lemma for triangles of Type
  C, which can all be obtained by mutating \(\Pi(5,2,29)\).

  For \(\Pi(1,5,13)=\mu_1(\Pi(1,2,5))\), we find that
  \(\ZZ\OP{inn}(\Pi(1,5,13))\) consists of\footnote{These are
    precisely the integer points common to both \(\Pi_A\) and
    \(\Pi_B\).}:
  \begin{align*}
    &\quad\quad(2,2),\cdots,(7,2)\\
    &\quad(1,1),\cdots,\cdots,\cdots,(14,1)\\
    &(0,0),\cdots,\cdots,\cdots,\cdots,\cdots,(20,0).
  \end{align*}
  Therefore these integral points are contained in all mutations
  of \(\Pi(1,5,13)\), that is all the triangles of Type A or
  B.

  If we perform only \(\mu_1\)-mutations then the integral points
  \((0,1)\), \((0,2)\), \((1,2)\), and \((0,3)\) stay in the
  polygon because they lie below the eigenline \(W_1\). This
  proves the lemma for Type A triangles.

  Suppose that \(\Pi\) is a Type B triangle obtained from a Type A
  triangle \(\Pi'=\mu_1^n\Pi(1,2,5)\) by performing a sequence of
  mutations starting with \(\mu_0\). Note that \(\Pi'\) has
  vertices \(\bm{p}_0=(0,0)\), \(\bm{p}_1=(8F_{2n+5}/F_{2n+3},0)\)
  and \(\bm{p}_2=(0,8F_{2n+3}/F_{2n+5})\) where
  \[F_1,F_2,F_3,F_4,F_5,\ldots = 1,1,2,3,5,\ldots\] is the
  Fibonacci sequence. The eigendirection \(\bm{w}_0\) for \(\Pi'\)
  is \((1,1)\), so the mutation \(\mu_0\) sends the integral point
  \((1,2)\) to \((0,1)\) and the integral points \((0,k)\) to
  \((-k,0)\). These latter points are absorbed into the bottom
  edge of \(\mu_0\Pi'\) and are certainly contained in
  \(\OP{inn}(\mu_0\Pi')\); it remains to show that \((0,1)\) lies
  in \(\OP{inn}(\mu_0\Pi')\).

  The bottom-left vertex of \(\mu_0\Pi'\) lies at
  \((-8F_{2n+3}/F_{2n+5},0)\) and the eigenline connects this to
  \((8/3,8/3)\). The point \((0,1)\) lies on or below this
  eigenline (and hence in \(\OP{inn}(\mu_0\Pi')\)) provided that
  \(8F_{2n+3}/F_{2n+5}\geq 8/5\). But the sequence
  \(8F_{2n+3}/F_{2n+5}\) converges to
  \(12-4\sqrt{5}\approx 3.0557\) from above, so the inequality
  \(8F_{2n+3}/F_{2n+5}\geq 8/5\) is certainly satisfied.
\end{proof}

\section{Geometry of the minimal resolution}
\label{sct:geom_min_res}

\pg Let \(Y=\HP(S)\) be a Manetti surface and \(B\subset Y\) an
octic curve. It will be important to understand the minimal
resolution \(\pi\colon\tilde{Y}\to Y\) and the total and proper
transforms (respectively \(\pi^*B\) and \(\pi^{-1}_*B\)) of
\(B\).

\paragraph{Mirror tropicalisation of the minimal resolution.}
\label{pg:mirror_trop_min_res} The minimal resolution
\(\tilde{Y}\) is again almost toric and the mirror
tropicalisation \(\Trop(\tilde{Y},\pi^*B)\) is again \(\Pi(S)\):
all of the components \(D_i\) of the exceptional locus gives us
a new edge \(e_i\) with inward normal \(\rho_i\) but having zero
length, residing at one of the vertices of \(\Pi(S)\) (whichever
vertex corresponds to the singularity being resolved). We will
now compute \(\Trop(\tilde{Y},\pi^{-1}_*B)\). Recall that we
defined Types A, B and C for Markov triangles in
\ref{pg:integer_points} and the polygons \(\Pi_K\) for
\(K=5,A,B,C\) in \ref{dfn:Z_polygons}.

\paragraph{Lemma.}\label{lma:trop_strict}
\begin{itemize}
\item {\em If \(Y=\PP(1,4,25)\) or \(\HP(5)\) then
    \(\Trop(\tilde{Y},\pi^{-1}_*B)\) is contained in
    \(\Pi_5\).}
\item {\em If \(\Trop(Y,B)\) is of Type K (for some
    \(K\in\{A,B,C\}\)) then
    \(\Trop(\tilde{Y},\pi^{-1}_*B)\subset\Pi_K\).}
\end{itemize}
\begin{proof}
  By Equation \eqref{eq:almost_toric_tropicalisation}, to find
  \(\Trop(\tilde{Y},\pi^{-1}_*B)\), we simply move each
  of these new edges \(e_i\) parallelly inwards by
  \(\varepsilon_i\rho_i\) where \(\varepsilon_i\) is defined by
  \[\pi^*B=\pi^{-1}_*B+\sum\varepsilon_iD_i.\]
  Since \(\pi^{-1}_*B\) is an integral Cartier divisor on a
  smooth variety, {\cite[Proposition 5.15]{Mandel}} implies that
  \(\Trop(\tilde{Y},\pi^*B)\) has vertices at integral
  points. By Lemma \ref{pg:integer_points}, this means that
  \(\Trop(\tilde{Y},\pi^{-1}_*B)\) is contained in one of the
  three polygons from Definition \ref{dfn:Z_polygons} according
  to its type.
\end{proof}

\pg\label{pg:min_res_notation} In what follows, we will be most
concerned with the Manetti surfaces \(Y=\HP(c)\). Recall that
\(Y\) has a cyclic quotient singularity of type
\(\frac{1}{c^2}(a^2,b^2)\), or equivalently
\(\frac{1}{c^2}(1,cq-1)\) for some \(1\leq q\leq c\) coprime to
\(c\), where \(b^2=(cq-1)a^2\mod c^2\). The exceptional locus of
the minimal resolution \(\pi\colon\tilde{Y}\to Y\) is a chain of
curves \(D_1,\ldots,D_n\) with self-intersections given by
\(D_i^2=-d_i\) where \([d_1,\ldots,d_n]\) is the continued
fraction expansion of \(\frac{c^2}{cq-1}\). The almost toric
boundary is the cycle \(D=\sum_{i=0}^nD_i\) of rational curves
where \(D_1,\ldots,D_n\) are the exceptional curves of the
minimal resolution and \(D_0\) is the proper transform of the
almost toric boundary of \(Y\). Note that \(D_0,D_1,\ldots,D_n\)
is a basis for the homology of \(\tilde{Y}\).

\pg\label{pg:ees} Let us write \(\eE_0\) and \(\eE_1\) for the
non-toric exceptional \(-1\)-curves in \(\tilde{Y}\) (we just
write \(\eE_1\) if there is only one). Define \(i_0\) and
\(i_1\) to be the indices such that
\[\eE_0\cdot D_i=\delta_{ii_0},\qquad \eE_1\cdot D_i=\delta_{ii_1}.\]
One can read off these indices from the tropical picture as
follows. Recall that the full mirror tropicalisation \(\Utrop\)
of \(\tilde{Y}\) has one or two branch cuts emanating from the
point \((8/3,8/3)\), parallel to \(\bm{w}_0\) and \(\bm{w}_1\)
(or just \(\bm{w}_1\) if there is only one). There are \(n+1\)
edges \(e_i\) (possibly of length zero) in any tropicalisation,
with inward normals \(\rho_i\) corresponding to the irreducible
components of \(D\). Since the eigendirection of the node
introduced by a non-toric blow-up is parallel to the edge you
blow up (see for example {\cite[Section 9.1]{EvansBook}}), the
edge \(e_{i_0}\) is parallel to \(\bm{w}_0\) and the edge
\(e_{i_1}\) is parallel to \(\bm{w}_1\).

\paragraph{Geometry from the mirror
  tropicalisation.}\label{pg:geom_mirror_trop}
We can read off much crucial information from the integral
affine manifold \(\Trop(\tilde{Y},\pi^{-1}_*B)\). Let us write
the homology class \(\pi^{-1}_*B\) in terms of the basis
\(D_0,\ldots,D_n\) as \(\sum_{i=0}^n\tilde{b}_iD_i\) so that
\(\Trop(\tilde{Y},\pi^{-1}_*B)\) is defined by the inequalities
\eqref{eq:almost_toric_tropicalisation}:
\(\langle u,\rho_{i}\rangle\geq -\tilde{b}_{i}\). Here we are
expressing \(u\) as a vector with respect to an origin at the
barycentre \((8/3,8/3)\). We have
\begin{itemize}
\item[(i)] The intersection number \(\pi^{-1}_*B\cdot D_i\) is
  equal to the degree of the line bundle \(\OO(\pi^{-1}_*B)\)
  restricted to \(D_i\), which is equal to the affine length of
  the edge corresponding to \(D_i\) in
  \(\Trop(\tilde{Y},\pi^{-1}_*B)\) by {\cite[Proposition
    5.18]{Mandel}}.
\item[(ii)] The intersection number \(\eE_k\cdot\pi^{-1}_*B\) is
  equal to
  \(\eE_k\cdot\left(\sum\tilde{b}_iD_i\right)=\tilde{b}_{i_k}\).
\item[(iii)] The multiplicity \(\varepsilon_i\) of \(D_i\) in
  \(\pi^*B=\pi^{-1}_*B+\sum \varepsilon_iD_i\) can be computed
  by comparing the integral affine manifolds
  \(\Trop(\tilde{Y},\pi^*B)\) and
  \(\Trop(\tilde{Y},\pi^{-1}_*B)\). The edges \(e_i\)
  corresponding to \(D_i\) in these two manifolds are parallel,
  but in the latter it has been translated an affine distance
  \(\varepsilon_i\) in the direction of its inward normal
  \(\rho_i\).
\end{itemize}
This will give us powerful geometric control over the pair
\((Y,B)\).

\pg\label{pg:moving_edges} Finally, we give a method to identify
\(\Trop(\tilde{Y},\pi^{-1}_*B)\). Suppose we have a candidate
\(\tilde{B}\) for \(\pi^{-1}_*B\) but it might be that \(\tilde{B}\)
contains irreducible components supported on the exceptional curves
\(D_i\). Let \(\sigma\) be a section of \(\OO(\tilde{B})\) for which
\(\tilde{B}=\sigma^{-1}(0)\). The space of sections of
\(\OO(\tilde{B})\) has a basis given by the GHK theta functions
corresponding to integer points of \(\Trop(\tilde{Y},\tilde{B})\). Let
\(\OP{supp}(\sigma)\subset\ZZ\Trop(\tilde{Y},\tilde{B})\) be the
minimal set of basis elements whose span contains \(\sigma\). Let
\(e\) be an edge of \(\OP{Trop}(\tilde{Y},\tilde{B})\) and let \(E\)
be the corresponding almost toric boundary curve. Suppose that \(e\)
is disjoint from the convex hull of \(\OP{supp}(\sigma)\). Then the
restriction of \(\sigma\) to \(E\) is zero (the nonvanishing sections
of \(\OO(\tilde{B})|_E\) correspond with integer points on
\(e\)). Therefore \(\tilde{B}=\tilde{B}'+E\) and we can shrink the
mirror tropicalisation by moving \(e\) normally inwards. Conversely,
if a point in \(\OP{supp}(\sigma)\) appears on \(e\) then
\(\sigma|_E\) is not identically zero, so \(E\) does not appear as a
component of \(\tilde{B}\) and so \(e\) is part of the boundary of
\(\Trop(\tilde{Y},\pi^{-1}_*B)\).

\section{Reducing to finitely many cases}
\label{sct:reducing_to_finite}

\pg In this section, we will freely use the notation established
in \ref{pg:min_res_notation}--\ref{pg:geom_mirror_trop}. Our
goal is to establish the following result, which reduces the
number of cases we need to consider to a finite number of
Manetti surfaces.

\paragraph{Theorem.}\label{pg:non_normal} {\em Let \(Y\) be a
  Manetti surface and \(B\) an octic Weil divisor. If the double
  cover of \(Y\) branched along \(B\) is normal then \(Y\) is
  one of the following possibilities:}
\begin{gather*}
  \PP(1,1,2^2),\quad \PP(1,2^2,5^2),\quad \PP(1,5^2,13^2),\quad
  \PP(5^2,2^2,29^2),\quad  \PP(1,13^2,34^2),\\
  \HP(5,29),\quad \HP(2,29),\quad \HP(5),\quad \HP(13),\quad \HP(29),\quad \HP(34).
\end{gather*}

\pg If \(B\) has a component \(B'\) with multiplicity
\(n\geq 2\) then the double cover fails to be normal: in local
analytic coordinates \((x,y)\) centred at a smooth point of
\(B'\) with \(B'=\{y=0\}\) component, the double cover is given
by \(u^2=y^nf(x,y)\). The cover then has non-normal
singularities along the ramification locus \(u=y=0\). So the
theorem will follow if we can show that \(B\) contains an
irreducible component with multiplicity \(\geq 2\). In fact, it
suffices to show that an octic \(B\subset \HP(c)\) has such a
multiple component if \(c\) is not one of \(1,2,5,13,29,34\),
since an octic in \(\HP(b,c)\) or \(\PP(a^2,b^2,c^2)\) is a
degeneration of an octic in \(\HP(c)\), and having a component
with multiplicity \(\geq 2\) is a closed condition. So Theorem
\ref{pg:non_normal} will follow from:

\paragraph{Proposition.}\label{prp:high_mult} {\em If
  \(c\neq 1,2,5,13,29,34\) then an octic curve in \(Y=\HP(c)\)
  has \(\pi(\eE_0)\) or \(\pi(\eE_1)\) as an irreducible
  component with multiplicity \(\geq 2\).}

Here, \(\eE_k\) are the curves in the minimal resolution
\(\pi\colon\tilde{Y}\to Y\) coming from the non-toric blow-ups
in the almost toric construction of \(\tilde{Y}\). Recall that
\(D_1,\ldots,D_n\) are the exceptional curves of \(\pi\), which
together with \(D_0\) form the almost toric boundary of
\(\tilde{Y}\). Write the homology class of \(\pi^{-1}_*B\) as
\(\sum_{i=0}^n\tilde{b}_iD_i\). The first ingredient in the
proof of Proposition \ref{prp:high_mult} is the following lemma.

\paragraph{Lemma.} {\em If \(\tilde{b}_{i_k}\leq -2\) then \(B\)
  contains \(\pi(\eE_k)\) as an irreducible component with
  multiplicity \(\geq 2\).}
\begin{proof}
  By \ref{pg:geom_mirror_trop}(ii), we have
  \(\eE_k\cdot\pi^{-1}_*B\leq -2\). Therefore \(\pi^{-1}_*B\)
  must contain \(\eE_k\) as an irreducible component with
  multiplicity at least \(2\), and hence \(B\) must contain
  \(\pi(\eE_k)\) as a component with multiplicity \(\geq 2\).
\end{proof}

Therefore Proposition \ref{prp:high_mult} follows from the next
proposition.

\paragraph{Proposition.}\label{prp:bound_on_displacement}
{\em If \((a, b, c)\) is not any of the following triples:
\begin{gather}\label{eq:list_of_triples}
  (1, 1, 1), \quad (1, 1, 2),\quad (1, 2, 5),\\
  \nonumber (1, 5, 13),\quad (1, 13, 34),\quad (2, 5,
  29),
\end{gather}
then either \(\tilde{b}_{i_0}\leq -2\) or
\(\tilde{b}_{i_1}\leq -2\).}

The proof of Proposition \ref{prp:bound_on_displacement}
occupies the remainder of this section.

\paragraph{Lemma.} {\em For a Markov triangle \(\Pi(a,b,c)\), we have}
\[\bm{w}_1\wedge\bm{w}_2= 3a,\quad\bm{w}_2\wedge\bm{w}_0= 3b,
\quad\bm{w}_0\wedge\bm{w}_1= 3c\]
\begin{proof}
  One can construct a quiver from a Markov triangle whose vertices
  correspond to vertices of the triangle and where the number of edges
  between vertices \(i\) and \(j\) is \(\bm{w}_i\wedge
  \bm{w}_j\). Mutation of the polygon induces mutation of this quiver
  (see {\cite[Proposition 3.17]{KNP}}). For \(\Pi(1,1,1)\) we get
  \(3,3,3\) and the quiver mutation rule for these arrows is the same
  as the mutation rule for trebled Markov triples \((3a,3b,3c)\).
\end{proof}

\paragraph{Definition.}\label{dfn:ps_and_ws}
Let \(\Pi\coloneqq\Pi(a,b,c)\) be a Markov triangle with vertices
\(\bm{p}_0,\bm{p}_1,\bm{p}_2\) and let \(\bm{n}\in\Pi\) be an integral
point. Let \(J\) be a clockwise \(90^\circ\) rotation, and let
\(\bm{w}_k^\perp=J\bm{w}_k\) for \(k=0,1,2\). Define
\[n_k\coloneqq (\bm{n}-\bm{p}_k)\cdot \bm{w}^\perp_k.\]
Note that if \(\bm{n}\in\OP{inn}(\Pi)\) then
\[n_0\geq 0,\qquad n_1\leq 0.\]

\begin{figure}[htb]
  \centering
  \begin{tikzpicture}
    \draw (0,0) -- (10,0) -- (3,4) -- cycle;
    \draw[dotted] (0,0) -- (4,2);
    \draw[dotted] (10,0) -- (4,2);
    \draw [dotted] (3,4) -- (4,2);
    \node at (0,0) [left] {\(\bm{p}_0\)};
    \node at (10,0) [right] {\(\bm{p}_1\)};
    \node at (3,4) [left] {\(\bm{p}_2\)};
    \draw[->] (3,1.5) --++ (0.5,-1) node[right] {\(\bm{w}^\perp_0\)};
    \draw[->] (4.5,11/6) --++ (1/3,1) node[left] {\(\bm{w}^\perp_1\)};
    \draw[->] (15/4,2.5) --++ (-1,-0.5) node[left] {\(\bm{w}^\perp_2\)};
    \node (n) at (2,0.5) {\(\bullet\)};
    \node at (n) [right] {\(\bm{n}\)};
    \node at (4,2) {\(\times\)};
  \end{tikzpicture}
  \caption{The points \(\bm{n}\), \(\bm{p}_i\) and vectors
    \(\bm{w}_i^\perp\) from Definition \ref{dfn:ps_and_ws}. (Not to
    scale for any Markov triple.)}
  \label{fig:ps_and_ws}
\end{figure}

\paragraph{Lemma.}\label{lma:mutation_n_formula}
{\em If \(\Pi'=\mu_0\Pi\) (and all the quantities/vectors for \(\Pi'\) are
just written as for \(\Pi\) but with primes) then
\[n'_0 = n_2+3bn_0,\qquad n'_1=n_1,\qquad n'_2=-n_0.\]
Similarly, if \(\Pi''=\mu_1\Pi\) then}
\[n''_0=n_0,\qquad n''_1 = n_2+3an_1,\qquad n''_2 = -n_1.\]
\begin{proof}
  We will just prove the formula for \(\mu_0\): the one for \(\mu_1\)
  is similar. We work in coordinates where the origin is at
  \(\bm{p}_0\). In particular, we have
  \(n_0=\bm{n}\cdot\bm{w}_0^\perp\). We have
  \begin{align*}
    \bm{p}'_0&\coloneqq \mu_0\bm{p}_2 = \bm{p}_2 +
               \left(\bm{p}_2\wedge\bm{w}_0\right)\bm{w}_0\\
    \bm{w}'_0&\coloneqq \mu_0\bm{w}_2 = \bm{w}_2 +
               \left(\bm{w}_2\wedge \bm{w}_0\right)\bm{w}_0\\
             &=\bm{w}_2 + 3b\bm{w}_0\\
    \left(\bm{w}'_0\right)^\perp&=J\bm{w}'_0=\bm{w}^\perp_2 +
                                  3b\bm{w}^\perp_0
  \end{align*}  
  Then we have
  \begin{align*}
    n'_0 &=
           \left(\bm{n} - \bm{p}'_0\right) \cdot \left(\bm{w}'_0\right)^\perp\\
         &=\left(\bm{n} - \bm{p}_2 - (\bm{p}_2\wedge\bm{w}_0)\bm{w}_0\right) \cdot \left(\bm{w}_2^\perp+3b\bm{w}_0^\perp\right)\\
         &=n_2+3b(\bm{n}-\bm{p}_2)\cdot\bm{w}_0^\perp -
           (\bm{p}_2\wedge \bm{w}_0)\bm{w}_0\cdot\bm{w}_2^\perp.
  \end{align*}
  We also have
  \(\bm{u}\cdot \bm{v}^\perp = \bm{v}\wedge\bm{u}\), so
  \((\bm{p}_2\wedge\bm{w}_0)\bm{w}_0\cdot\bm{w}_2^\perp =
  (\bm{w}_0\wedge\bm{w}_2)\bm{p}_2\cdot\bm{w}_0^\perp =
  -3b\bm{p}_2\cdot\bm{w}_0^\perp\). Therefore
  \[n'_0 = n_2+3b\bm{n}\cdot\bm{w}_0^\perp=n_2+3bn_0.\qedhere\]
\end{proof}

\paragraph{Lemma.}\label{lma:mutation_technical}
{\em Suppose that \(\Pi\) is a Markov triangle and that \(\bm{n}\) is
  a point in \(\OP{inn}(\Pi)\). Suppose that \(n_2\geq (1-3b)n_0\)
  (respectively \(n_2\leq (1-3a)n_1\)). Then the same inequality holds
  for any \(\mu_{\bm{\bin}}\Pi\).}
\begin{proof}
  We will focus on the first inequality; the proof of the second
  is similar. By induction, it is sufficient to prove the claim
  for \(\Pi'=\Pi(a',b',c')=\mu_0\Pi\) and
  \(\Pi''=\Pi(a'',b'',c'')=\mu_1\Pi\). For \(\Pi''\), the
  quantity \(n''_2\) is positive and hence automatically
  satisfies \(n''_2\geq (1-3b'')n''_0\). For \(\Pi'\), we
  have \[n'_0 = n_2+3bn_0,\quad n'_2 = -n_0.\] Since
  \(b'\geq 1\), we have \(1-3b'\leq -2\). Therefore the
  assumption \(n_2\geq (1-3b)n_0\)
  implies
  \begin{equation}\label{eq:mutation_technical}
    (1-3b')(n_2+3bn_0)\leq
    (1-3b')n_0.\end{equation} Therefore
  \begin{align*}
    n'_2&=-n_0\\
    &\geq (1-3b')n_0&&\mbox{since }1-3b'<-1\\
    &\geq (1-3b')(n_2+3bn_0)&&\mbox{by Equation \eqref{eq:mutation_technical}}\\
    &=(1-3b')n'_0&&\mbox{as required.}\qedhere
  \end{align*}
\end{proof}

\paragraph{Corollary.}\label{cor:mutation_technical}
{\em Suppose that \(\Pi\) is a Markov triangle
and that \(\bm{n}\) is a point in \(\OP{inn}(\Pi)\). Suppose that
\begin{itemize}
\item \(n_2\geq (1-3b)n_0\) and \(n_0\geq 2\), or
\item \(n_2\leq (1-3a)n_1\) and \(n_1\leq -2\).
\end{itemize}
Then these inequalities persist for any mutation
\(\mu_{\bm{\bin}}\Pi\).}
\begin{proof}
  By Lemma \ref{lma:mutation_technical}, the inequality
  \(n_2\geq (1-3b)n_0\) persists for any mutation. This means
  that under a \(\mu_0\) mutation, \(n'_0=n_2+3bn_0\geq n_0\),
  so \(n_0\) cannot decrease. Under a \(\mu_1\) mutation,
  \(n_0\) is unchanged. Therefore \(n_0\geq 2\) for any
  combination of mutations. The proof for \(n_1\leq -2\) is
  similar.
\end{proof}

\pg\label{pg:displacement_argument} We can now prove Proposition
\ref{prp:bound_on_displacement}. Let \(\Pi=\Pi(a,b,c)\) be a
Markov triangle where \(c\) is the biggest number in the triple
and let \(Y=\HP(c)\) be the corresponding Manetti surface. As in
the proof of Lemma \ref{lma:trop_strict}, we imagine that
\(\Pi\) has a bunch of length-zero edges \(e_i\) corresponding
to curves of the minimal resolution \(\tilde{Y}\) of the
corresponding Manetti surface. The tropicalisation
\(\Trop(\tilde{Y},\pi_*^{-1}B)\) is then obtained from \(\Pi\)
by moving the edges parallely inwards along their normal
directions. As the edge \(e_{i_0}\) moves inwards along its
normal direction, it first encounters the eigenline \(W_0\) (to
which it is parallel) and then it moves over \(\OP{inn}(\Pi)\)
until it reaches the closest integral point, \(\bm{n}\). During
this final stage of its motion, the affine displacement it
suffers is \(n_0\). On the other side of the polygon, the edge
\(e_{i_1}\) similarly undergoes a displacement of at least
\(n_1\) where \(\bm{n}\) is a different closest integral point.

If \((a,b,c)=(29,2,169)\) (corresponding to the Type C triangle
\(\mu_{(0,0)}\Pi(1,2,5)\)) and \(\bm{n}=(-1,1)\) then \(n_0=2\) and
\(n_2=0\). By Lemma \ref{lma:trop_strict}, the same integral point is
the closest to \(e_{i_0}\) for this (Type C) triple and all its
mutations. Corollary \ref{cor:mutation_technical}, the inequality
\(n_0\geq 2\) persists for all subsequent mutations, so we deduce that
\(\tilde{b}_{i_0}\leq -2\) for such triples.

Similar arguments work with the triples shown in Table
\ref{tbl:displacement_cases}. Any triple which is not listed in the
statement of Proposition \ref{prp:bound_on_displacement} is either in
Table \ref{tbl:displacement_cases} or is obtained from such a triple
by mutation, and hence satisfies either \(\tilde{b}_{i_0}\leq -n_0\leq
-2\) or \(\tilde{b}_{i_1}\leq n_1\leq -2\). This completes the proof
of Proposition \ref{prp:bound_on_displacement}.

\begin{table}[htb]
  \centering
  \begin{tabular}{llllll}
    \(\bm{\bin}\) & \((a,b,c)\) & Type & \(\bm{n}\) & &\\
    \hline
    \((1,1,1)\) & \((1,34,89)\) & A & \((14,1)\) & \(n_1=-3\) & \(n_2=1\) \\
    \((1,1,0)\) & \((34,13,1325)\) & B & \((0,1)\) & \(n_0=31\) & \(n_2=1\)\\
    \((1,0)\) & \((13,5,194)\) & B & \((0,1)\) & \(n_0=12\) & \(n_2=1\)\\
    \((0,1)\) & \((5,29,433)\) & C & \((13,1)\) & \(n_1=-27\) & \(n_2=1\) \\
    \((0,0)\) & \((29,2,169)\) & C & \((-1,1)\) & \(n_0 = 2\) & \(n_2=0\)
  \end{tabular}
  \caption{Cases for the proof of Proposition
    \ref{prp:bound_on_displacement}. Each row refers to the triangle \(\mu_{\bm{\bin}}\Pi(1,2,5)\).}
  \label{tbl:displacement_cases}
\end{table}

\section{Calculating discrepancies}
\label{sct:discrep}

\paragraph{Partial resolution of the double
  cover.}\label{pg:partial_res_double_cover} Let \(\pi\colon
\tilde{Y}\to Y\) be a resolution of singularities with exceptional
divisor \(E=\sum_i E_i\). Suppose that \(n\) is an odd number such
that \(nB\) is Cartier on \(Y\) and that \(\LL\) is a line bundle on
\(Y\) with \(\LL^{\otimes 2}\cong\OO(nB)\). Then \(\pi^*\LL\) is a
line bundle on \(\tilde{Y}\) which is a square root of
\(\OO(\pi^*(nB))\); if we take the associated double cover then it
will usually fail to be normal because \(\pi^*(nB)\) has multiplicity
\(n\) along the proper transform \(\pi^{-1}_*B\) and potentially
nonzero multiplicities \(b_i\) along the exceptional curves \(E_i\):
\[\pi^*(nB)=n\pi^{-1}_*B+\sum b_iE_i.\]
We can reduce these multiplicities to \(1\) along
\(\pi^{-1}_*B\) and either \(0\) or \(1\) along each \(E_i\) by
tensoring \(\pi^*(\LL)\) with
\(\OO(\pi^{-1}_*B)^{\otimes (n-1)/2}\) and with
\(\bigotimes_i\OO(-E_i)^{\otimes \lfloor b_i/2\rfloor}\). The
associated double cover \(\tilde{f}\colon\tilde{X}\to\tilde{Y}\)
is normal and branched along
\(\pi^{-1}_*B+\sum_{b_i\text{ odd}}E_i\). If we collapse
\(\tilde{f}^{-1}_*E\) then we obtain a normal surface which is a
double cover of \(Y\) branched along\footnote{In general, there might be branching over singular points of \(Y\). This does not occur in any of the cases we encounter, so we ignore it here. To see why it doesn't occur, observe that the only Manetti surface for which an octic can be disjoint from the singular point is \(\PP(1,1,4)\) (because otherwise the moment-image projection of the singular point fails to be an integral point of the tropicalisation) and for \(\PP(1,1,4)\) we will rule out this kind of branching in \ref{pg:1_1_4_I}.} \(B\). By uniqueness of the
Alexeev--Pardini construction described in
\ref{pg:ap_double_covers}, this coincides with \(X\). Thus
\(\tilde{X}\) is a partial resolution of \(X\).

\paragraph{Necessary condition for log canonical singularities.}
\label{pg:necessary_condition} Let \(\pi\colon\tilde{Y}\to Y\)
be the minimal resolution of \(Y=\HP(c)\) with exceptional
curves \(D_1,\ldots,D_n\) and let \(B\subset Y\) be an octic
curve. If we let \(Z=\tilde{Y}\) and \(\Delta=B/2\) in
Eq. \eqref{eq:discrep} and intersect both sides with \(D_i\)
then we get
\[-2-D_i^2+\frac{1}{2}\pi^{-1}_*B\cdot D_i=\sum_j
  a(D_j,Y,B/2)D_i\cdot D_j\] where we used the adjunction
formula \(K_{\tilde{Y}}\cdot D_i=-2-D_i^2\). Let us write
\(\bm{\alpha}=(\alpha_1,\ldots,\alpha_n)\) with
\(\alpha_i=a(D_i,Y,B/2)\) and
\(\bm{\beta}=(\beta_1,\ldots,\beta_n)\) with
\(\beta_i=-2-D_i^2+\frac{1}{2}\pi^{-1}_*B\cdot D_i\), and let
\(N\) be the matrix with \(N_{ij}=D_i\cdot D_j\). Then
\begin{equation}\label{eq:alpha}\bm{\alpha}=N^{-1}\bm{\beta}.\end{equation}
We can figure out the numbers \(\beta_i\) from the mirror
tropicalisation \(\Trop(\tilde{Y},\pi^{-1}_*B)\), which then gives the
discrepancies \(\alpha_i\). If any of these discrepancies is \(<-1\)
then the pair \((Y,B/2)\) fails to be log canonical, so the double
cover \(X\) of \(Y\) branched along \(B\) also fails to be log
canonical by \ref{pg:useful_properties_lt_lc}(4).

\section{Classification}
\label{sct:classif}

\pg We now proceed case by case with the list of Manetti
surfaces given in Theorem \ref{pg:non_normal}, with the aim of
classifying octic curves in these surfaces which could yield a
log canonical octic double. The octic curves correspond to
sections of the sheaf \(\OO(B)\) spanned by the GHK theta
functions \(\theta_p\) corresponding to the integer points \(p\)
in the mirror tropicalisation, as enumerated in Section
\ref{sct:integer_points}. If \(\theta=\sum a_p\theta_p\) is such
a section then we say its {\em support}
\(\OP{supp}(\theta)\subset\ZZ\Trop(Y,B)\) is the set of integer
points \(p\) in the mirror tropicalisation with \(a_p\neq
0\). We will stratify the space of sections of this sheaf
according to their support: write \(\Gamma_P\) for the subset of
sections whose support is a specific set
\(P\subset\ZZ\Trop(Y,B)\) of integer points. If
\(\OP{supp}(\theta)\subset P\) then there is a one-parameter
deformation
\(\theta_t=\theta+t\sum_{p\in P\setminus\OP{supp}(\theta)}
a_p\theta_p\) with \(\theta_t\in \Gamma_P\) for \(t\neq 0\); we
are free to choose the coefficients \(a_p\) to ensure that
\(\theta_t\) belongs to any specified Zariski-dense set. In
\ref{pg:hp_13}--\ref{pg:hp_34}, we will prove:

\paragraph{Theorem.}\label{thm:non_lc} {\em For \(Y=\HP(13)\),
  \(Y=\HP(29)\) and \(Y=\HP(34)\), there is a Zariski-dense set
  \(U\) in the top stratum \(\Gamma_{\ZZ\Trop(Y,B)}\) such that
  \(\left(Y,\frac{1}{2}\theta^{-1}(0)\right)\) is not log
  canonical for any \(\theta\in U\).}

\paragraph{Corollary.}\label{cor:ignore_13_29_34} {\em If \(X\)
  is a normal octic double of \(Y=\HP(13)\), \(Y=\HP(29)\) or
  \(Y=\HP(34)\) then \(X\) cannot be log canonical.}
\begin{proof}
  If \(B=\theta^{-1}(0)\) is an octic curve in \(Y\) such that
  the double cover \(X\) has log canonical singularities then
  the same will be true for any small deformation of \(B\) by
  \ref{pg:useful_properties_lt_lc}(2). But we can find a small
  deformation \(\theta_t\) with \(\theta_t\in U\) for
  \(t\neq 0\), which ensures that \((Y,B_t/2)\) is not log
  canonical, so that the double cover \(X_t\) branched along
  \(B_t\) is not log canonical by
  \ref{pg:useful_properties_lt_lc}(4).
\end{proof}

\pg In a similar way, we can ignore \(\PP(5^2,2^2,29^2)\),
\(\HP(5,29)\), \(\HP(2,29)\), \(\PP(1,5^2,13^2)\), and
\(\PP(1,13^2,34^2)\) because they are degenerations of
\(\HP(29)\), \(\HP(13)\) and \(\HP(34)\). This leaves us with
the following three cases to analyse:
\[\PP(1,1,4),\quad \PP(1,4,25),\quad \HP(5).\]
We will analyse these three cases first, and then proceed to the
cases of \(\HP(13)\), \(\HP(29)\) and \(\HP(34)\) which will
consistute the proof of Theorem \ref{thm:non_lc}.

\paragraph{\(\PP(1,1,4)\).} We will work this case out in a little
more detail because it is simple and explicit enough to illustrate the
ideas in the toric setting. If \(Y=\PP(1,1,4)\) with weighted
homogeneous coordinates \([x:y:z]\) then an octic curve is a Weil
divisor defined by the vanishing of a section \(\Omega\) of the
\(\QQ\)-line bundle \(\OO(16)\), that is a polynomial in \(x,y,z\) of
weighted degree \(16\). The mirror tropicalisation is \(\Pi(1,1,2)\),
shown in Figure \ref{fig:p114}. In this figure, we have also labelled
a selection of the integer points with monomials \(\OO(16)\) according
to their affine distances to the three edges (see e.g. {\cite[Example
    5.4.5]{CoxLittleSchenck}}). The branched double cover is
\[X = \left\{[u:x:y:z]\in\PP(8,1,1,4) \,:\,
u^2=\Omega(x,y,z)\right\}.\] We consider three cases:
\begin{itemize}
\item Case I: The \(\frac{1}{4}(1,1)\) singularity of \(Y\) at
  \([0:0:1]\) does not belong to \(B\). Equivalently, \((0,4)\in
  \OP{supp}(\Omega)\).
\item Case II: \((0,4)\not\in\OP{supp}(\Omega)\) but
  \((i,3)\in\OP{supp}(\Omega)\) for some \(i\).
\item Case III: \((0,4)\not\in\OP{supp}(\Omega)\) and
  \((i,3)\not\in\OP{supp}(\Omega)\) for \(i=0,1,2,3,4\).
\end{itemize}

\begin{figure}[htb]
  \begin{tikzpicture}[scale=0.7]
    \node at (0,0) [left] {\(y^{16}\)};
    \foreach \x in {0,1,...,16} {\node at (\x,0) {\(\bullet\)};}
    \foreach \x in {0,1,...,12} {\node at (\x,1) {\(\cdot\)};}
    \foreach \x in {0,1,...,8} {\node at (\x,2) {\(\cdot\)};}
    \foreach \x in {0,1,...,4} {\node at (\x,3) {\(\cdot\)};}
    \foreach \y in {1,...,4} {\node at (0,\y) {\(\bullet\)};}
    \foreach \y in {0,1,...,4} {\node at ({16-4*\y},\y) {\(\bullet\)};}
    \draw (0,0) -- (16,0) node[right] {\(x^{16}\)} -- (0,4) node [left] {\(z^4\)} -- cycle;
    \node at (4,3) [above right] {\(x^4z^3\)};
    \node at (8,2) [above right] {\(x^8z^2\)};
    \node at (12,1) [above right] {\(x^{12}z\)};
  \end{tikzpicture}
  \caption{The mirror tropicalisation \(\Pi(1,1,2)\) of \(\PP(1,1,4)\)
    with some integral points labelled by monomial sections from
    \(\OO(16)\).}
  \label{fig:p114}
\end{figure}

\paragraph{\(\PP(1,1,4)\), Case I: \((0,4)\in\OP{supp}(\Omega)\).}\label{pg:1_1_4_I}
This case was already considered by Anthes {\cite[Example
  5.5]{Anthes}}. If \([0:0:1]\not\in B\) then our branched
double cover has precisely two non-Gorenstein singularities,
\([\pm \sqrt{\Omega(0,0,1)}:0:0:1]\in X\subset\PP(8,1,1,4)\)
each of type \(\frac{1}{4}(1,1)\). Note that these are distinct:
if they were the same point in \(\PP(8,1,1,4)\) then we would
have
\[[\lambda^8\sqrt{\Omega(0,0,1)}:0:0:\lambda^4]=[-\sqrt{\Omega(0,0,1)}:0:0:1]\]
so \(\lambda^4=1\) and \(\lambda^8=-1\), which is impossible. In
particular, it is not possible to have branching over the
singular point of \(\PP(1,1,4)\); this kind of branching could
only occur if the weight of the \(u\) coordinate were congruent
to \(2\) modulo \(4\), which yields {\em degenerate covers} in
the terminology of \cite{Manetti2}.

\paragraph{\(\PP(1,1,4)\), Cases II and III:}\label{pg:1_1_4_II}
In the remaining cases, the octic \(B\) passes through the
\(\frac{1}{4}(1,1)\) singularity of \(Y\). Let \(\pi\colon\tilde{Y}\to
Y\) be the minimal resolution; this is toric: its fan contains one new
ray corresponding to the exceptional curve \(C_1\) of square
\(-4\). Since the corresponding edge in
\(\Trop(\tilde{Y},\pi^*B)=\Trop(Y,B)\) has length zero, \(C_1\cdot
\pi^*B=0\), but since \(B\) passes through the singularity, this means
that \(\pi^*B\) contains \(C_1\) as an irreducible component. Case II
is when \(\pi^*B=\pi^{-1}_*B+C_1\) and Case III is when
\(\pi^*B=\pi^{-1}_*B+nC_1\) for \(n\geq 2\). This means that in Case
II, \(\Trop(\tilde{Y},\pi^{-1}_*B)\) is the polygon shown in Figure
\ref{fig:p114_II}.

\begin{figure}[htb]
  \center
  \begin{tikzpicture}[scale=0.7]
    \foreach \x in {0,1,...,16} {\node at (\x,0) {\(\bullet\)};}
    \foreach \x in {0,1,...,12} {\node at (\x,1) {\(\cdot\)};}
    \foreach \x in {0,1,...,8} {\node at (\x,2) {\(\cdot\)};}
    \foreach \x in {0,1,...,4} {\node at (\x,3) {\(\bullet\)};}
    \foreach \y in {1,...,3} {\node at (0,\y) {\(\bullet\)};}
    \foreach \y in {0,1,...,3} {\node at ({16-4*\y},\y) {\(\bullet\)};}
    \draw (0,0) -- (16,0) -- (4,3) -- (0,3) -- cycle;
  \end{tikzpicture}
  \caption{\(\Trop(\tilde{Y},\pi^{-1}_*(B))\) in Case II.}
  \label{fig:p114_II}
\end{figure}

In Case III, the edge corresponding to \(C_1\) will have moved
even further downwards and have affine length at least \(8\); by
\ref{pg:geom_mirror_trop}(i), this means
\(C_1\cdot \pi^{-1}_*B\geq 8\). We can therefore rule out Case
III using the necessary condition for log canonical
singularities from \ref{pg:necessary_condition}: we have
\(N=(-4)\) and \(\bm{\beta}=(2+C_1\cdot\pi^{-1}_*B)\), so the
necessary condition for \((Y,B/2)\) to be log canonical is
\(C_1\cdot\pi^{-1}_*B\leq 4\).

Therefore we can focus on Case II. We distinguish five subcases
according to how \(B\) and \(C_1\) intersect:
\begin{itemize}
\item[(a)] four transverse intersections,
\item[(b)] one intersection with multiplicity \(2\) and two transverse intersections,
\item[(c)] two intersections with multiplicity \(2\),
\item[(d)] one intersection with multiplicity \(3\) and one transverse intersection,
\item[(e)] one intersection with multiplicity \(4\).
\end{itemize}
We will find that in Subcases (d) and (e) the double cover must
have a singularity which is not log canonical. Subcases (a-c)
will yield singularities whose dual graphs are shown in Figures
\ref{fig:graphs}(a-c) respectively; these are
\(\ZZ/2\)-quotients of simple elliptic (in Subcase (a)) or cusp
(in Subcases (b-c)) singularities.\footnote{In an earlier
  version of this paper, some of these cases were missing,
  specifically those in which \(B\) is not smooth along
  \(B\cap C_1\). We are grateful to H. Asaike, M. Enokizono,
  M. Hattori and Y. Koto for drawing our attention to this gap
  and refer the reader to their paper \cite{AEHK}
  for a full classification of normal degenerations of Horikawa
  surfaces with any geometric genus.}

We first focus on the generic situation where \(B\) is smooth
along \(B\cap C_1\). In each subcase, we can find the partial
resolution \(\tilde{X}\) of the double cover \(X\) following the
prescription outlined in \ref{pg:partial_res_double_cover}. In
fact, we can go further by taking a more refined resolution
\(\dbtilde{Y}\to \tilde{Y}\) which is a log resolution of the
branch curve. In Case II(a), this means blowing up the four
intersections \(\pi^{-1}_*B\cap C_1\) to obtain \(-1\)-curves
\(E_1,E_2,E_3,E_4\). The proper transform of \(C_1\) becomes a
curve of square \(-8\); when we pass to the double cover, it
becomes a curve \(\overline{C}_1\) with square \(-4\); the
curves \(E_i\) become curves \(\overline{E}_i\) of square
\(-2\). See Figure \ref{fig:min_res_114} for a pictorial summary
of the argument.

\begin{figure}[htb]
  \begin{center}
    \begin{tikzpicture}[scale=0.8]
      \node (y) at (-1,0) {\(Y\)};
      \node (y2) at (-1,3) {\(\tilde{Y}\)};
      \draw[->] (y2) -- (y);
      \node (y3) at (-1,6) {\(\dbtilde{Y}\)};
      \draw[->] (y3) -- (y2);
      \node (y3) at (-1,6) {\(\dbtilde{Y}\)};
      \draw[->] (y3) -- (y2);
      \node (x) at (-3,0) {\(X\)};
      \node (x2) at (-3,3) {\(\tilde{X}\)};
      \draw[->] (x2) -- (x);
      \node (x3) at (-3,6) {\(\dbtilde{X}\)};
      \draw[->] (x3) -- (x2);
      \draw[->] (x3) -- (y3) node (ff1) [midway,above] {\(\dbtilde{f}\)};
      \node at (ff1) [below=1em] {2:1};
      \draw[->] (x2) -- (y2) node (ff2) [midway,above] {\(\tilde{f}\)};
      \node at (ff2) [below=1em] {2:1};
      \draw[->] (x) -- (y) node (ff3) [midway,above] {\(f\)};
      \node at (ff3) [below=1em] {2:1};
      \begin{scope}[shift={(0.5,0)}]
        \node (sing) at (2,0) {\(\bullet\)};
        \node (singlabel) at (4,0) {\([0:0:1]\)};
        \draw[->] (singlabel) -- (sing);
        \draw[dashed,line width=0.6mm] (0,-1) to[out=60,in=210] (2,0) to[out=30,in=-30] (3,1) to[out=150,in=60] (2,0) to[out=-120,in=180] (2,-1) to[out=0,in=-60] (2,0) to[out=120,in=30] (1,1) to[out=-150,in=150] (2,0) to[out=-30,in=120] (4,-1) node [right] {\(B\)};

        \begin{scope}[shift={(0,3)}]
          \draw (0,0) -- (4,0) node [right] {\(C_1\)};
          \node at (2,0) [above] {\(-4\)};
          \draw[dashed,line width=0.6mm] (0,-1) to[out=60,in=180] (1,1) to[out=0,in=180] (2,-1) to[out=0,in=180] (3,1) to[out=0,in=120] (4,-1) node [right] {\(B\)};
        \end{scope}

        \begin{scope}[shift={(0,6)}]
          \draw (0,1) -- (4,1) node [right] {\(C_1\)};
          \node at (2,1) [above] {\(-8\)};
          \draw[dashed,line width=0.6mm] (0,-1) -- (4,-1) node [right] {\(B\)};
          \draw (0.5,1.2) -- (0.5,-1.2) node [midway,right] {\(E_1\)};
          \draw (1.5,1.2) -- (1.5,-1.2) node [midway,right] {\(E_2\)};
          \draw (2.5,1.2) -- (2.5,-1.2) node [midway,right] {\(E_3\)};
          \draw (3.5,1.2) -- (3.5,-1.2) node [midway,right] {\(E_4\)};
        \end{scope}
      \end{scope}

      \begin{scope}[shift={(-1,0)}]
        \begin{scope}[shift={(-8,0)}]
          \node at (2,0) {\(\bullet\)};
          \node at (2,0) [below] {singularity};
        \end{scope}

        \begin{scope}[shift={(-8,3)}]
          \draw (0,0) -- (4,0);
          \node at (2,0) [below] {four nodes};
          \node at (0.5,0) {\(\bullet\)};
          \node at (1.5,0) {\(\bullet\)};
          \node at (2.5,0) {\(\bullet\)};
          \node at (3.5,0) {\(\bullet\)};
        \end{scope}

        \begin{scope}[shift={(-8,6)}]
          \draw (4,0) -- (0,0) node [left] {\(-4\)};
          \node at (4,0) [right] {\(\overline{C}_1\)};
          \draw (0.5,-0.5) -- (0.5,0.5) node [above] {\(-2\)};
          \node at (0.5,-0.5) [below] {\(\overline{E}_1\)};
          \draw (1.5,-0.5) -- (1.5,0.5) node [above] {\(-2\)};
          \node at (1.5,-0.5) [below] {\(\overline{E}_2\)};
          \draw (2.5,-0.5) -- (2.5,0.5) node [above] {\(-2\)};
          \node at (2.5,-0.5) [below] {\(\overline{E}_3\)};
          \draw (3.5,-0.5) -- (3.5,0.5) node [above] {\(-2\)};
          \node at (3.5,-0.5) [below] {\(\overline{E}_4\)};
        \end{scope}
      \end{scope}
    \end{tikzpicture}
  \caption{The process for finding the minimal resolution of the
    double cover of \(Y=\PP(1,1,4)\) branched along a curve \(B\) in
    Case II(a).}
  \label{fig:min_res_114}
\end{center}
\end{figure}

The result is that in Case II(a), \(X\) has a
singularity whose minimal resolution has the following dual
graph:

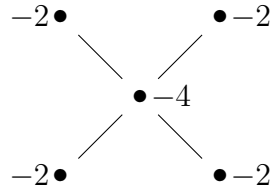
\begin{center}
  \begin{tikzpicture}
    \node (a) at (0,0) {\(\circ\)};
    \node at (a) [right] {\(-4\)};
    \node (b) at (45:1.5) {\(\bullet\)};
    \node at (b) [right] {\(-2\)};
    \node (c) at (135:1.5) {\(\bullet\)};
    \node at (c) [left] {\(-2\)};
    \node (d) at (-135:1.5) {\(\bullet\)};
    \node at (d) [left] {\(-2\)};
    \node (e) at (-45:1.5) {\(\bullet\)};
    \node at (e) [right] {\(-2\)};
    \draw (a) -- (b);
    \draw (a) -- (c);
    \draw (a) -- (d);
    \draw (a) -- (e);
  \end{tikzpicture}
\end{center}

This is a \(\ZZ/2\)-quotient of a simple elliptic
singularity. We summarise the analogous results in Subcases
II(b--e) in Figure \ref{fig:114_more} when \(B\) is assumed to
be smooth along \(B\cap C_1\). We see that II(b) and II(c) yield
\(\ZZ/2\)-quotients of cusp singularities (which are log
canonical) whilst II(d) and II(e) yield singularities that are
not log canonical.

\begin{figure}[htb]
  \centering
  \begin{tikzpicture}[scale=0.7]
    \draw[line width=0.6mm] (0,0) -- (4,0) node [right] {\(C_1\)};
    \draw[dashed,line width=0.6mm] (0,-1) to[out=45,in=180] (1,0) to[out=0,in=180] (2,-1) to[out=0,in=180] (3,1) to[out=0,in=135] (4,-1) node [right] {\(B\)};
    \begin{scope}[shift={(0,2.5)}]
      \draw[line width=0.6mm] (0,1) -- (4,1) node [pos=0.4,above] {\(-8\)};
      \node at (4,1) [right] {\(C_1\)};
      \draw[dashed,line width=0.6mm] (0,0) -- (4,0) node [right] {\(B\)};
      \draw (1,1.2) -- (1,-1.2) node [pos=0.28,left] {\(-1\)};
      \draw (0,-1) -- (2,-1) node [pos=0.2,below] {};
      \draw (2.5,1.2) -- (2.5,-0.2) node [midway,left] {\(-1\)};
      \draw (3.5,1.2) -- (3.5,-0.2) node [midway,left] {\(-1\)};
    \end{scope}
    \begin{scope}[shift={(0,6)}]
      \draw (0,1) -- (4,1) node[pos=0.4,above] {\(-4\)};
      \draw (2.5,1.2) -- (2.5,-0.2) node [midway,left] {};
      \draw (3.5,1.2) -- (3.5,-0.2) node [midway,left] {};
      \draw (1,1.2) -- (1,-1.2) node [pos=0.28,left] {};
      \draw (0,-1) -- (2,-1) node [pos=0.2,below] {};
      \draw (0,-0.5) -- (2,-0.5) node [pos=0.2,above] {};
    \end{scope}
    \begin{scope}[shift={(2,9)}]
      \node (a) at (-1,0) {\(\bullet\)};
      \node (b) at (1,0) {\(\circ\)};
      \node (c) at (2,1) {\(\bullet\)};
      \node (d) at (2,-1) {\(\bullet\)};
      \node (e) at (-2,1) {\(\bullet\)};
      \node (f) at (-2,-1) {\(\bullet\)};
      \draw (a) -- (b) -- (c);
      \draw (b) -- (d);
      \draw (a) -- (e);
      \draw (a) -- (f);
    \end{scope}
    \begin{scope}[shift={(0,11)}]
      \node at (0,0) {(b)}; 
    \end{scope}
    \begin{scope}[shift={(5.5,0)}]
      \draw[line width=0.6mm] (0,0) -- (4,0) node [right] {\(C_1\)};
      \draw[dashed,line width=0.6mm] (0,-1) to[out=45,in=180] (1,0) to[out=0,in=180] (2,-1) to[out=0,in=180] (3,0) to[out=0,in=135] (4,-1) node[right] {\(B\)};
      \begin{scope}[shift={(0,2.5)}]
        \draw[line width=0.6mm](0,1) -- (4,1) node [midway,above] {\(-8\)};
        \node at (4,1) [right] {\(C_1\)};
        \draw[dashed,line width=0.6mm] (0,0) -- (4,0) node [right] {\(B\)};
        \draw (1,1.2) -- (1,-1.2) node [pos=0.28,left] {\(-1\)};
        \draw (0,-1) -- (1.5,-1) node [pos=0.2,below] {};
        \draw (3,1.2) -- (3,-1.2) node [pos=0.28,right] {\(-1\)};
        \draw (2.5,-1) -- (4,-1) node [pos=0.8,below] {};
      \end{scope}
      \begin{scope}[shift={(0,6)}]
        \draw (0,1) -- (4,1) node[midway,above] {\(-4\)};
        \draw (1,1.2) -- (1,-1.2) node [pos=0.28,left] {};
        \draw (0,-1) -- (1.5,-1) node [pos=0.2,below] {};
        \draw (0,-0.5) -- (1.5,-0.5) node [pos=0.2,above] {};
        \draw (3,1.2) -- (3,-1.2) node [pos=0.28,left] {};
        \draw (2.5,-1) -- (4,-1) node [pos=0.8,below] {};
        \draw (2.5,-0.5) -- (4,-0.5) node [pos=0.8,above] {};
      \end{scope}
      \begin{scope}[shift={(2,9)}]
        \node (a) at (-1,0) {\(\bullet\)};
        \node (g) at (0,0) {\(\circ\)};
        \node (b) at (1,0) {\(\bullet\)};
        \node (c) at (2,1) {\(\bullet\)};
        \node (d) at (2,-1) {\(\bullet\)};
        \node (e) at (-2,1) {\(\bullet\)};
        \node (f) at (-2,-1) {\(\bullet\)};
        \draw (a) -- (g) -- (b) -- (c);
        \draw (b) -- (d);
        \draw (a) -- (e);
        \draw (a) -- (f);
      \end{scope}
      \begin{scope}[shift={(0,11)}]
        \node at (0,0) {(c)};
      \end{scope}
    \end{scope}
    \begin{scope}[shift={(11,0)}]
      \draw[line width=0.6mm] (0,0) -- (4,0) node [right] {\(C_1\)};
      \draw[dashed,line width=0.6mm] (0,-1) to[out=45,in=180] (1.5,0) to[out=0,in=180] (3,1) to[out=0,in=135] (4,-1) node [right] {\(B\)};
      \begin{scope}[shift={(0,2.5)}]
        \draw[line width=0.6mm] (0,1) -- (4,1) node [pos=0.65,above] {\(-8\)};
        \node at (4,1) [right] {\(C_1\)};
        \draw[dashed,line width=0.6mm] (0,0) -- (4,0) node [right] {\(B\)};
        \draw (1.5,1.2) -- (1.5,-1.2) node [pos=0.28,left] {\(-1\)};
        \draw (0,-1) -- (3,-1) node [midway,above] {};
        \draw (3.5,1.2) -- (3.5,-0.2) node [midway,left] {\(-1\)};
        \draw (0.5,-1.2) -- (0.5,-0.5) node [left] {};
      \end{scope}
      \begin{scope}[shift={(0,6)}]
        \draw (0,1) -- (4,1) node[pos=0.65,above] {\(-4\)};
        \draw (1.5,1.2) -- (1.5,-1.2) node [pos=0.7,left] {};
        \draw (0,-1) -- (3,-1) node [pos=0.2,below] {};
        \draw (0,0) -- (3,0) node [pos=0.2,above] {};
        \draw (3.5,1.2) -- (3.5,-0.2) node [pos=0.6,right] {};
        \draw (2.5,-0.3) -- (2.5,0.3) node [above] {};
        \draw (2.5,-0.7) -- (2.5,-1.3) node [below] {};
      \end{scope}
      \begin{scope}[shift={(2,10)}]
        \node (a2) at (0,0) {\(\bullet\)};
        \node (b2) at (1,0) {\(\bullet\)};
        \node (c2) at (2,0) {\(\bullet\)};
        \node (d2) at (-1,0) {\(\bullet\)};
        \node (e2) at (-2,0) {\(\bullet\)};
        \node (f2) at (0,-1) {\(\circ\)};
        \node (g2) at (0,-2) {\(\bullet\)};
        \draw (e2) -- (d2) -- (a2) -- (b2) -- (c2);
        \draw (a2) -- (f2) -- (g2);
      \end{scope}
      \begin{scope}[shift={(0,11)}]
        \node at (0,0) {(d)};
      \end{scope}
    \end{scope}
    \begin{scope}[shift={(16.5,0)}]
      \draw[line width=0.6mm] (0,0) -- (4,0) node [right] {\(C_1\)};
      \draw[dashed,line width=0.6mm] (0,-1) to[out=45,in=180] (2,0) to[out=0,in=135] (4,-1) node [right] {\(B\)};
      \begin{scope}[shift={(0,2.5)}]
        \draw[line width=0.6mm] (0,1) -- (4,1) node [pos=0.25,above] {\(-8\)};
        \node at (4,1) [right] {\(C_1\)};
        \draw[dashed,line width=0.6mm] (0,0) -- (4,0) node [right] {\(B\)};
        \draw (2,1.2) -- (2,-1.2) node [pos=0.28,left] {\(-1\)};
        \draw (0,-1) -- (4,-1) node [pos=0.8,above] {};
        \draw (0.5,-1.2) -- (0.5,-0.5) node [left] {};
        \draw (0,-0.6) -- (1.5,-0.6);
      \end{scope}
      \begin{scope}[shift={(0,6)}]
        \draw (0,1) -- (4,1) node[pos=0.25,above] {\(-4\)};
        \draw (2,1.2) -- (2,-1.2) node [pos=0.28,left] {};
        \draw (0,-1) -- (4,-1) node [pos=0.8,above] {};
        \draw (0.5,-1.2) -- (0.5,-0.5) node [left] {};
        \draw (0,-0.6) -- (1.5,-0.6);
        \draw (0,0) -- (4,0) node [pos=0.8,above] {};
        \draw (3.5,-0.2) -- (3.5,0.5) node [left] {};
        \draw (2.5,0.4) -- (4,0.4);
      \end{scope}
      \begin{scope}[shift={(2,9)}]
        \node (a1) at (0,0) {\(\bullet\)};
        \node (b1) at (1,0) {\(\bullet\)};
        \node (c1) at (2,0) {\(\bullet\)};
        \node (d1) at (-1,0) {\(\bullet\)};
        \node (e1) at (-2,0) {\(\bullet\)};
        \node (f1) at (0,-1) {\(\circ\)};
        \node (g1) at (-2,-1) {\(\bullet\)};
        \node (h1) at (2,-1) {\(\bullet\)};          
        \draw (g1) --(e1) -- (d1) -- (a1) -- (b1) -- (c1) -- (h1);
        \draw (a1) -- (f1);
      \end{scope}
      \begin{scope}[shift={(0,11)}]
        \node at (0,0) {(e)};
      \end{scope}
    \end{scope}
  \end{tikzpicture}
  \caption{The singularities of the branched double cover in
    Cases II(b--e) assuming that \(B\) is smooth along
    \(B\cap C_1\). From bottom: the intersection between \(B\)
    and \(C_1\) in \(\tilde{Y}\); the log resolution; the double
    cover of the log resolution; the dual graph of the
    exceptional locus of the singularity of the branched double
    cover. Any unmarked curves other than \(B\) have square
    \(-2\). In the resolution graphs, \(\bullet\) indicates a
    \(-2\)-curve whilst \(\circ\) indicates a \(-4\)-curve.}
  \label{fig:114_more}
\end{figure}

If \(B\) is not smooth along \(B\cap C_1\) then we must be in
one of Subcases (b-e).

\paragraph{\(\PP(1,1,4)\), Case II(d-e), \(B\) not smooth.}
In Subcases (d) and (e), if \(B\) is not
smooth then we can perturb its equation locally analytically in
a neighbourhood of \(C_1\) so that it is locally modelled on a
smooth curve intersecting \(C_1\) with multiplicity \(3\) or
\(4\), which means that the singularity of the branched double
cover is adjacent to the non-log canonical singularities in
Figure \ref{fig:114_more}(d-e). This means these subcases can
never yield log canonical singularities.

\paragraph{\(\PP(1,1,4)\), Case II(b), \(B\) not smooth.} In
Subcase (b), let \(p\in B\cap C_1\) be the point
with\footnote{Here, \(i_p(B,C_1)\) denotes the local
  intersection number.} \(i_p(B,C_1)=2\). If \(B\) is not
smooth at \(p\) then it has multiplicity \(2\). There are two possibilities:
\begin{itemize}
\item[(i)] The germ of \(B\) at \(p\) has two smooth branches
  \(B_1\) and \(B_2\) which are both transverse to \(C_1\) and
  intersect one another with multiplicity \(k\geq 1\).
\item[(ii)] The germ of \(B\) at \(p\) is irreducible and it has
  a single tangent plane \(T_pB\) which is distinct from
  \(T_pC_1\).
\end{itemize}
In Case II(b.i), we can blow-up the two transverse intersections
of \(C_1\) with \(B\) (writing \(B_3\) and \(B_4\) for the germs
of \(B\) at these points) and also \(k\) times infinitesimally
close to \(p\) to obtain the following configuration:

\begin{center}
  \begin{tikzpicture}
    \node (w) at (-2,1) {\(B_3\)};
    \node (x) at (-1,1) {\(\bullet\)};
    \node (y) at (-2,-1) {\(B_4\)};
    \node (z) at (-1,-1) {\(\bullet\)};
    \node at (x) [above] {\(-1\)};
    \node at (z) [above] {\(-1\)};
    \node (a) at (0,0) {\(C_1\)};
    \node (b) at (1,0) {\(\bullet\)};
    \node at (b) [above] {\(-2\)};
    \node (c) at (2,0) {\(\cdots\)};
    \node (d) at (3,0) {\(\bullet\)};
    \node at (d) [above] {\(-2\)};
    \node (e) at (4,0) {\(\bullet\)};
    \node at (e) [above] {\(-1\)};
    \node (f) at (5,1) {\(B_1\)};
    \node (g) at (5,-1) {\(B_2\)};
    \draw (w) -- (x);
    \draw (x) -- (a);
    \draw (y) -- (z);
    \draw (z) -- (a);
    \draw (a) -- (b);
    \draw (b) -- (c);
    \draw (c) -- (d);
    \draw (d) -- (e);
    \draw (e) -- (f);
    \draw (e) -- (g);
    \draw [decorate,decoration={brace,amplitude=5pt,mirror}] (1,-0.5) -- (4,-0.5) node[midway,yshift=-1em]{\(k\)};
  \end{tikzpicture}
\end{center}

The \(k+2\) curves introduced by this blow-up are all part of the branch locus, so we further
blow-up all the intersection points to obtain

\begin{center}
  \begin{tikzpicture}
    \node (w) at (-2,1) {\(B_3\)};
    \node (x) at (-1,1) {\(\bullet\)};
    \node (y) at (-2,-1) {\(B_4\)};
    \node (z) at (-1,-1) {\(\bullet\)};
    \node at (x) [above] {\(-1\)};
    \node at (z) [above] {\(-1\)};
    \node (a) at (0,0) {\(C_1\)};
    \node (b) at (1,0) {\(\bullet\)};
    \node at (b) [above] {\(-1\)};
    \node (c) at (2,0) {\(\bullet\)};
    \node at (c) [above] {\(-4\)};
    \node (d) at (3,0) {\(\bullet\)};
    \node at (d) [above] {\(-1\)};
    \node (e) at (4,0) {\(\cdots\)};
    \node (f) at (5,0) {\(\bullet\)};
    \node at (f) [above] {\(-4\)};
    \node (g) at (6,1) {\(\bullet\)};
    \node at (g) [above] {\(-1\)};
    \node (h) at (7,1) {\(B_1\)};
    \node (i) at (6,-1) {\(\bullet\)};
    \node at (i) [above] {\(-1\)};
    \node (j) at (7,-1) {\(B_2\)};
    \draw (w) -- (x);
    \draw (x) -- (a);
    \draw (y) -- (z);
    \draw (z) -- (a);
    \draw (a) -- (b);
    \draw (b) -- (c);
    \draw (c) -- (d);
    \draw (d) -- (e);
    \draw (e) -- (f);
    \draw (f) -- (g);
    \draw (g) -- (h);
    \draw (f) -- (i);
    \draw (i) -- (j);
    \draw [decorate,decoration={brace,amplitude=5pt,mirror}] (1,-0.5) -- (5,-0.5) node[midway,yshift=-1em]{\(2k\)};
  \end{tikzpicture}
\end{center}

Here, \(C_1\) is a \(-8\)-curve and everything except the
\(-1\)-curves are part of the branch locus. When we take the
branched double cover, we obtain:

\begin{center}
  \begin{tikzpicture}
    \node (x) at (-1,1) {\(\bullet\)};
    \node (z) at (-1,-1) {\(\bullet\)};
    \node at (x) [above] {\(-2\)};
    \node at (z) [above] {\(-2\)};
    \node (a) at (0,0) {\(\bullet\)};
    \node at (a) [above] {\(-4\)};
    \node (b) at (1,0) {\(\bullet\)};
    \node at (b) [above] {\(-2\)};
    \node (c) at (2,0) {\(\bullet\)};
    \node at (c) [above] {\(-2\)};
    \node (d) at (3,0) {\(\bullet\)};
    \node at (d) [above] {\(-2\)};
    \node (e) at (4,0) {\(\cdots\)};
    \node (f) at (5,0) {\(\bullet\)};
    \node at (f) [above] {\(-2\)};
    \node (g) at (6,1) {\(\bullet\)};
    \node at (g) [above] {\(-2\)};
    \node (i) at (6,-1) {\(\bullet\)};
    \node at (i) [above] {\(-2\)};
    \draw (x) -- (a);
    \draw (z) -- (a);
    \draw (a) -- (b);
    \draw (b) -- (c);
    \draw (c) -- (d);
    \draw (d) -- (e);
    \draw (e) -- (f);
    \draw (f) -- (g);
    \draw (f) -- (i);
    \draw [decorate,decoration={brace,amplitude=5pt,mirror}] (1,-0.5) -- (5,-0.5) node[midway,yshift=-1em]{\(2k\)};
  \end{tikzpicture}
\end{center}

In Case II(b.ii), the picture is similar. Since the multiplicity
of \(B\) at \(p\) is equal to \(2\), the Puiseux expansion can
have at most one term, and the germ of \(B\) at \(p\) is locally
analytically equivalent to \(y^2=x^{2k+1}\). This yields:

\begin{center}
  \begin{tikzpicture}
    \node (w) at (-2,1) {\(B\)};
    \node (x) at (-1,1) {\(\bullet\)};
    \node (y) at (-2,-1) {\(B\)};
    \node (z) at (-1,-1) {\(\bullet\)};
    \node at (x) [above] {\(-1\)};
    \node at (z) [above] {\(-1\)};
    \node (a) at (0,0) {\(C_1\)};
    \node (b) at (1,0) {\(\bullet\)};
    \node at (b) [above] {\(-1\)};
    \node (c) at (2,0) {\(\bullet\)};
    \node at (c) [above] {\(-4\)};
    \node (d) at (3,0) {\(\bullet\)};
    \node at (d) [above] {\(-1\)};
    \node (e) at (4,0) {\(\cdots\)};
    \node (f) at (5,0) {\(\bullet\)};
    \node at (f) [above] {\(-4\)};
    \node (g) at (6,0) {\(\bullet\)};
    \node at (g) [above] {\(-1\)};
    \node (h) at (7,1) {\(B\)};
    \node (i) at (7,-1) {\(\bullet\)};
    \node at (i) [above] {\(-2\)};
    \draw (w) -- (x);
    \draw (x) -- (a);
    \draw (y) -- (z);
    \draw (z) -- (a);
    \draw (a) -- (b);
    \draw (b) -- (c);
    \draw (c) -- (d);
    \draw (d) -- (e);
    \draw (e) -- (f);
    \draw (f) -- (g);
    \draw (g) -- (h);
    \draw (g) -- (i);
    \draw [decorate,decoration={brace,amplitude=5pt,mirror}] (1,-0.5) -- (5,-0.5) node[midway,yshift=-1em]{\(2k\)};
  \end{tikzpicture}
\end{center}

and in the double cover:

\begin{center}
  \begin{tikzpicture}
    \node (x) at (-1,1) {\(\bullet\)};
    \node (z) at (-1,-1) {\(\bullet\)};
    \node at (x) [above] {\(-2\)};
    \node at (z) [above] {\(-2\)};
    \node (a) at (0,0) {\(\bullet\)};
    \node at (a) [above] {\(-4\)};
    \node (b) at (1,0) {\(\bullet\)};
    \node at (b) [above] {\(-2\)};
    \node (c) at (2,0) {\(\bullet\)};
    \node at (c) [above] {\(-2\)};
    \node (d) at (3,0) {\(\bullet\)};
    \node at (d) [above] {\(-2\)};
    \node (e) at (4,0) {\(\cdots\)};
    \node (f) at (5,0) {\(\bullet\)};
    \node at (f) [above] {\(-2\)};
    \node (g) at (6,1) {\(\bullet\)};
    \node at (g) [above] {\(-2\)};
    \node (i) at (6,-1) {\(\bullet\)};
    \node at (i) [above] {\(-2\)};
    \draw (x) -- (a);
    \draw (z) -- (a);
    \draw (a) -- (b);
    \draw (b) -- (c);
    \draw (c) -- (d);
    \draw (d) -- (e);
    \draw (e) -- (f);
    \draw (f) -- (g);
    \draw (f) -- (i);
    \draw [decorate,decoration={brace,amplitude=5pt,mirror}] (1,-0.5) -- (5,-0.5) node[midway,yshift=-1em]{\(2k+1\)};
  \end{tikzpicture}
\end{center}

\paragraph{\(\PP(1,1,4)\), Case II(c), \(B\) not smooth.} The
local analysis is the same as Case II(b) except that now we end
up with two forked chains of \(-2\)-curves coming out of the
\(-4\)-curve living over \(C_1\).

\paragraph{\(\PP(1,1,4)\), Case II(b-c): further
  bounds.}\label{pg:bound_on_ts} It
remains to bound the number of vertices of the resolution
graph. It helps to think of \(\tilde{Y}\cong\FF_4\) as a GIT quotient
\(\CC^4/(\CC^*)^2\), where the action is written in terms of Cox
variables as
\[[x:y:z:w]=[\lambda x:\lambda y:\mu \lambda^4z:\mu w].\] Now
our branch curve \(B\cup C_1\) is given by an equation of the
form \(w\Omega(x,y,z,w)\) where \(\Omega\) transforms like
\(\lambda^{16}\mu^3\) (if we think of \(\tilde{Y}\) as \(\FF_4\)
then \(B\) lives in the same linear system as \(3C_1+16F\) where
\(F=\{x=0\}\) is a ruling). Suppose that \(p\in B\cap C_1\) is a
singular point of \(B\). Using an automorphism of the Hirzebruch
surface \(\tilde{Y}\), we can assume that
\(p=[0:1:1:0]\). Working in the affine chart \(y=z=1\), the
polynomial \(\Omega(x,w)\) is a sum of monomials \(x^aw^{3-b}\)
indexed by the integral points
\((a,b)\in \ZZ\Trop\left(\tilde{Y},\pi^{-1}_*B\right)\). In Case
II(b-c), since the multiplicity of the singular point is \(2\),
the quadratic term in \(\Omega(x,w)\) is nonzero, so the corank
of this quadratic form is at most \(1\) and, by Arnold's
determinator of singularities {\cite[Section 16.2]{AGZV}}, the
germ of \(B\) at \(p\) is locally analytically equivalent to an
\(A_q\) singularity for some \(q\), that is \(\xi^2=\eta^{q+1}\)
for some analytic functions \(\xi(x,w)\) and \(\eta(x,w)\). The
highest order term possible in \(\Omega\) is \(x^{16}w^3\), so
the multiplicity of this \(A_q\) singularity is at most
\(18\). This is realised by, for example
\[\Omega(x,w) = (w-x)^2 - (w-x)^3- x^{16}w^3,\]
which belongs to Case II(b), giving the singularity from Figure
\ref{fig:graphs}(b) with \(t=19\). In Case II(c), the same
argument tells us that the multiplicities of each singularity
separately are bounded above by \(18\), so that graph (c) from
Figure \ref{fig:graphs} satisfies \(t_1+t_2\leq 38\).

However, we conjecture that \(t_1+t_2\leq 21\), which is
realised by
\[\Omega(x,w) = (w-x)^2 - x^{16}w^3,\]
which has an \(A_{18}\) singularity at \([0:1:1:0]\) and an
\(A_2\) singularity at \([1:0:1:0]\), giving the singularity
from Figure \ref{fig:graphs}(c) with \(t_1=19\) and
\(t_2=3\). The reason for our conjecture is that in order for
\(\Omega\) to have an \(A_{18}\) singularity at \([0:1:1:0]\) we
need it to have a term \(x^{16}w^3\), but in the affine chart
\([1:y:1:w]\) this term is cubic (\(w^3\)), so one would have to
choose the remaining terms very carefully in order to be able to
absorb this into the quadratic part of \(\Omega(y,w)\) by an
analytic change of coordinates. These extra terms would also
have to be chosen in such a way that they do not reduce the
multiplicity of the singularity at \([0:1:1:0]\). This seems
like too many constraints.

\paragraph{\(\PP(1,4,25)\).} In this case,
\(\OP{Trop}(Y,B)=\Pi(1,2,5)\) and the integer points correspond
to ``octic'' monomials, i.e. weighted degree \(80=8abc\) in the
weighted homogeneous coordinates \([x:y:z]\). We have marked
with a \(\star\)
(or \begin{tikzpicture}[baseline=-1mm]\node[circle,fill=none,draw=black,inner
  sep=0pt,minimum size=1pt] at (0,0)
  {\(\star\)};\end{tikzpicture} when on the boundary) the most
important integer points
\[\bm{q}_1=(0,3),\quad \bm{q}_2=(1,3),\quad \bm{q}_3=(7,2),\quad
  \bm{q}_4=(13,1),\quad \bm{q}_5=(19,0),\quad\bm{q}_6=(20,0);\] the rest we have drawn with a
\(\bullet\) if they lie on the boundary of the polygon and a
\(\cdot\) if they are on the interior.

\begin{tikzpicture}[scale=0.6]
  \foreach \x in {0,1,2,...,18} {\node at (\x,0) {\(\bullet\)};}
  \foreach \x in {1,...,13} {\node at (\x,1) {\(\cdot\)};}
  \foreach \x in {1,...,6} {\node at (\x,2) {\(\cdot\)};}
  \node at (1,3) {\(\cdot\)};
  \foreach \y in {1,2} {\node at (0,\y) {\(\bullet\)};}
  \node[circle,fill=none,draw=black,inner sep=0pt,minimum size=1pt] at (0,3) {\(\star\)};
  \node at (1,3) {\(\star\)};
  \node at (7,2) {\(\star\)};
  \node at (13,1) {\(\star\)};
  \node[circle,fill=none,draw=black,inner sep=0pt,minimum size=1pt] at (20,0) {\(\star\)};
  \draw (0,0) -- (20,0) -- (0,16/5) node [above] {\((0,16/5)\)} -- cycle;
  \node at (0,3) [left] {\(\bm{q}_1\)};
  \node at (1,3) [below right] {\(\bm{q}_2\)};
  \node at (7,2) [above] {\(\bm{q}_3\)};
  \node at (13,1) [above] {\(\bm{q}_4\)};
  \node at (20.1,0) [right] {\(\bm{q}_6\)};
  \node[circle,fill=none,draw=black,inner sep=0pt,minimum size=1pt] (q5) at (19,0) {\(\star\)};
  \node at (q5) [below] {\(\bm{q}_5\)};
  \node at (9.5,1.7) [above] {\(e_x\)};
  \node at (0,1.5) [left] {\(e_y\)};
  \node at (9.5,0) [below] {\(e_z\)};
\end{tikzpicture}

The inward normal rays to the edges \(e_x\),
\(e_y\) and \(e_z\) are \(\rho_x=(-4,-25)\), \(\rho_y=(1,0)\)
and \(\rho_z=(0,1)\) corresponding to the Cox variables
\(x,y,z\) of weights \(1,4,25\) respectively. We write \(D_x\),
\(D_y\) and \(D_z\) for the corresponding divisors. In our
tables below, we will keep track of the edge \(e_x\) with normal
\(\rho_x\) but \(e_y\) and \(e_z\) will play no role.

The singularities of \(\PP(1,4,25)\) comprise:
\begin{itemize}
\item a \(\frac{1}{4}(1,1)\) singularity living over \((20,0)\),
\item a \(\frac{1}{25}(1,4)\) singularity living over \((0,16/5)\).
\end{itemize}
Let \(\pi\colon \tilde{Y}\to Y\) be the minimal resolution of
\(Y=\PP(1,4,25)\). The exceptional locus \(\OP{Exc}(\pi)\) comprises
Wahl chain \(C_1,C_2,C_3,C_4\) with self-intersections \(-7,-2,-2,-2\)
and a \(-4\)-curve \(C_5\). This is still a toric variety: the mirror
tropicalisation \(\Trop(\tilde{Y},\pi^*B)\) is still \(\Pi(1,2,5)\)
but the fan acquires new rays
\[\rho_1=(0,-1),\ \rho_2=(-1,-7),\ \rho_3=(-2,-13),\ \rho_4=(-3,-19),\ \rho_5=(-1,-6)\]
which we think of as the normals to zero-length edges
\(e_1,\ldots,e_4\) at the vertex \((0,16/5)\) and \(e_5\) at
\((20,0)\). Since \((0,16/5)\) is not an integer point, all of our
octics must pass through the \(\frac{1}{25}(1,4)\) singularity. Note
that \(\bm{q}_6\in\OP{supp}(\Omega)\) if and only if
\(B=\Omega^{-1}(0)\) does not pass through the \(\frac{1}{4}(1,1)\)
singularity. We consider the following exhaustive list of cases:
\begin{itemize}
\item[I.] \(\bm{q}_2,\bm{q}_6\in\OP{supp}(\Omega)\).
\item[II.] \(\bm{q}_6\not\in\OP{supp}(\Omega)\) but at least one of
  \(\bm{q}_2,\ldots,\bm{q}_5\) is in \(\OP{supp}(\Omega)\).
\item[III.] \(\bm{q}_6\not\in\OP{supp}(\Omega)\) and
  \(\bm{q}_2,\ldots,\bm{q}_5\not\in\OP{supp}(\Omega)\).
\item[IV.] \(\bm{q}_6\in\OP{supp}(\Omega)\) and
  \(\bm{q}_2\not\in\OP{supp}(\Omega)\).
\end{itemize}

\paragraph{\(\PP(1,4,25)\), Case I.}
Let \(\pi\colon \tilde{Y}\to Y\) be the minimal resolution. Lemma
\ref{lma:trop_strict} implies that
\(\Trop(\tilde{Y},\pi^{-1}_*B)\subset\Pi_5\). For a moment, let's
assume that at least one point of \(e_y\) is in
\(\OP{supp}(\Omega)\). By \ref{pg:moving_edges}, our assumptions on
the integer points in \(\OP{supp}(\Omega)\) imply that
\(\Trop(\tilde{Y},\pi^{-1}_*B)=\Pi_5\).

\begin{center}
  \begin{tikzpicture}[scale=0.6]
    \foreach \x in {0,1,2,...,19} {\node at (\x,0) {\(\bullet\)};}
    \foreach \x in {1,...,12} {\node at (\x,1) {\(\cdot\)};}
    \foreach \x in {1,...,6} {\node at (\x,2) {\(\cdot\)};}
    \foreach \x in {0,1} {\node at (\x,3) {\(\bullet\)};}
    \node at (7,2) {\(\cdot\)};
    \node at (13,1) {\(\cdot\)};
    \node at (20,0) {\(\bullet\)};
    \foreach \y in {0,1,2,3} {\node at (0,\y) {\(\bullet\)};}
    \draw (0,0) -- (20,0) node[below] {\(\bm{q}_6\)} -- (1,3) node [above] {\(\bm{q}_2\)} -- (0,3) node[above] {\(\bm{q}_1\)} -- cycle;
    \node at (0.6,2.6) {\(e_1\)};
    \node at (9,2) {\(e_4\)};
    \node at (9.5,0) [below] {\(e_z\)};
    \node at (0,1.5) [left] {\(e_y\)};
  \end{tikzpicture}
\end{center}

We label only those edges whose affine length is nonzero;
\(e_2\) and \(e_3\) are concentrated at the point \(\bm{q}_2\)
and \(e_x\) and \(e_5\) are concentrated at the point
\(\bm{q}_5\). Recall from \ref{pg:displacement_argument} that
the tropicalisation \(\Trop(\tilde{Y},\pi_*^{-1}B)\) is obtained
from the Markov triangle \(\Pi(1,2,5)\) by displacing the edges
parallely in their inward normal directions. The following table
summarises the affine lengths of the new edges together with how
far (in affine distance) the edges were displaced in passing
from the Markov triangle to the tropicalisation.

\begin{center}
  \begin{tabular}{ccccccc}
    Edge & \(e_1\) & \(e_2\) & \(e_3\) & \(e_4\) & \(e_x\) & \(e_5\) \\
    \hline Affine length & \(1\) & \(0\) & \(0\) & \(1\) & \(0\) & \(0\) \\
    Affine displacement & \(1/5\) & \(2/5\) & \(3/5\) & \(4/5\) & \(0\) & \(0\)
  \end{tabular}
\end{center}

This means that \(\pi^{-1}_*B\) intersects each of \(C_1\) and \(C_4\)
transversely at one point, is disjoint from \(C_2\), \(C_3\) and
\(C_5\) and that the branch curve specified in
\ref{pg:partial_res_double_cover} is \(\pi^{-1}_*B+C_1+C_3\). We take
a log resolution \(\dbtilde{Y}\to\tilde{Y}\) of this branch locus by
blowing up \(\pi^{-1}_*B\cap C_1\) to get a \(-1\)-curve \(E_1\). When
we take the branched cover
\(\dbtilde{f}\colon\dbtilde{X}\to\dbtilde{Y}\), the chain
\(E_1,C_1,\ldots,C_4\) turns into a chain of curves with
self-intersections \(-2\), \(-4\), \(-4\), \(-1\), \(-4\), which is a
resolution of the \(\frac{1}{50}(1,29)\) Wahl singularity. The curve
\(C_5\) turns into two \(-4\) curves, each of which yield
\(\frac{1}{4}(1,1)\) singularities in the branched double cover \(X\to
Y\). See Figure \ref{fig:p_1_4_25_case_I}.

If \(\OP{supp}(\Omega)\) does not contain a point of \(e_y\) then
\(e_y\) may move to the right by an affine distance of \(1\). The only
way this affects the result is that \(D_y\) becomes an irreducible
component of \(\pi^{-1}_*B\), but the analysis remains unaffected.

\begin{figure}[htb]
  \centering
  \begin{tikzpicture}[every text node part/.style={align=center}]
    \draw[line width=0.6mm] (-0.2,-0.2) -- (2.2,0.2) node [midway,sloped] {\(C_1\) \\ \(-7\)};
    \draw (1.8,0.2) -- (4.2,-0.2) node [midway,sloped] {\(C_2\) \\ \(-2\)};
    \draw[line width=0.6mm] (3.8,-0.2) -- (6.2,0.2) node [midway,sloped] {\(C_3\) \\ \(-2\)};
    \draw (5.8,0.2) -- (8.2,-0.2) node [midway,sloped] {\(C_4\) \\ \(-2\)};
    \draw[dashed] (7.8,-0.2) -- (10.2,0.2) node [midway,sloped] {\(D_x\) \\ \(-1\)};;
    \draw (9.8,0.2) -- (12.2,-0.2) node [midway,sloped] {\(C_5\) \\ \(-4\)};
    \draw[dashed,line width=0.6mm] (0,-0.5) to[out=90,in=180] (0.6,0.8) -- (6,0.8) to[out=0,in=90] (6.6,-0.5);
    \node at (0,0.8) {\(B\)};
    \node (ytilde) at (-1.5,0) {\(\tilde{Y}\)};
    \begin{scope}[shift={(0,2.5)}]
      \draw[line width=0.6mm] (-0.2,-0.2) -- (2.2,0.2) node [midway,sloped] {\(C_1\) \\ \(-8\)};
      \draw (1.8,0.2) -- (4.2,-0.2) node [midway,sloped] {\(C_2\) \\ \(-2\)};
      \draw[line width=0.6mm] (3.8,-0.2) -- (6.2,0.2) node [midway,sloped] {\(C_3\) \\ \(-2\)};
      \draw (5.8,0.2) -- (8.2,-0.2) node [midway,sloped] {\(C_4\) \\ \(-2\)};
      \draw[dashed] (7.8,-0.2) -- (10.2,0.2) node [midway,sloped] {\(D_x\) \\ \(-1\)};;
      \draw (9.8,0.2) -- (12.2,-0.2) node [midway,sloped] {\(C_5\) \\ \(-4\)};
      \draw[dashed,line width=0.6mm] (-0.2,0.8) -- (6,0.8) to[out=0,in=90] (6.6,-0.5);
      \draw (0,-0.4) -- (0,1) node [midway] {\(-1\,\, E_1\,\,\)};
      \node at (-0.2,0.8) [above] {\(B\)};
      \node (ytildetilde) at (-1.5,0) {\(\dbtilde{Y}\)};
      \draw[->] (ytildetilde) -- (ytilde);
    \end{scope}
    \begin{scope}[shift={(0,5)}]
      \draw (-0.2,-0.2) -- (2.2,0.2) node [midway,sloped,below] {\(-4\)};
      \draw (1.8,0.2) -- (4.2,-0.2) node [midway,sloped,below] {\(-4\)};
      \draw (3.8,-0.2) -- (6.2,0.2) node [midway,sloped,below] {\(-1\)};
      \draw (5.8,0.2) -- (8.2,-0.2) node [midway,sloped,below] {\(-4\)};
      \draw[dashed] (7.8,-0.2) -- (10.2,0.2) node [midway,sloped,below] {\(-1\)};;
      \draw (9.8,0.2) -- (12.2,-0.2) node [midway,sloped,below] {\(-4\)};
      \draw[dashed] (6.8,0) -- (9.2,0.6) node [midway,sloped,above] {\(-1\)};;
      \draw (8.7,0.5) -- (11.2, 0.5) node [midway,sloped,above] {\(-4\)};
      \draw (0,-0.4) -- (0,1) node [midway,left] {\(-2\)};
      \node (xtildetilde) at (-1.5,0) {\(\dbtilde{X}\)};
      \draw[->] (xtildetilde) -- (ytildetilde) node (ff1) [midway,right] {2:1};
      \node at (ff1) [left=1em] {\(\dbtilde{f}\)};
    \end{scope}
    \begin{scope}[shift={(0,7)}]
      \node (a) at (0,0) {\(\bullet\)};
      \node (b) at (2,0) {\(\bullet\)};
      \node (c) at (4,0) {\(\bullet\)};
      \node (d) at (6,0) {\(\bullet\)};
      \node (e) at (8,0) {\(\bullet\)};
      \node (f) at (10,0) {\(\bullet\)};
      \draw (a) -- (b);
      \draw (b) -- (c);
      \draw (c) -- (d);
      \node at (a) [above] {\(-2\)};
      \node at (b) [above] {\(-4\)};
      \node at (c) [above] {\(-3\)};
      \node at (d) [above] {\(-3\)};
      \node at (e) [above] {\(-4\)};
      \node at (f) [above] {\(-4\)};
    \end{scope}
  \end{tikzpicture}
  \caption{\(\PP(1,4,25)\), Case I. The minimal resolution
    \(\tilde{Y}\), its blow-up \(\dbtilde{Y}\), and the double cover
    \(\dbtilde{X}\) of \(\dbtilde{Y}\). Thick curves are the branch
    locus. Any non-dashed curves in \(\dbtilde{X}\) are contracted to
    get \(X\). The top row shows the Wahl chains associated to the
    three singularities of \(X\).}
  \label{fig:p_1_4_25_case_I}
\end{figure}

\paragraph{\(\PP(1,4,25)\), Case II.} The only difference between this
case and the previous one is that, with \(\bm{q}_6\) gone, the edges
\(e_x\) and \(e_5\) move normally inwards and yield the following data:

\begin{center}
  \begin{tikzpicture}[scale=0.6,every text node part/.style={align=center}]
    \foreach \x in {1,2,...,19} {\node at (\x,0) {\(\bullet\)};}
    \foreach \x in {1,...,12} {\node at (\x,1) {\(\cdot\)};}
    \foreach \x in {1,...,6} {\node at (\x,2) {\(\cdot\)};}
    \foreach \x in {0,1} {\node at (\x,3) {\(\bullet\)};}
    \foreach \y in {0,1,2,3} {\node at (0,\y) {\(\bullet\)};}
    \node at (7,2) {\(\bullet\)};
    \node at (13,1) {\(\bullet\)};
    \node at (20,0) {\(\circ\)};
    \draw (0,0) -- (19,0) node [below] {\(\bm{q}_5\)} -- (1,3) node [above] {\(\bm{q}_2\)} -- (0,3) node[above] {\(\bm{q}_1\)} -- cycle;
    \node at (19,0) [above] {\((19,0)\)};
    \node at (19,0) {\(\bullet\)};
    \node at (7,2) [above] {\(\bm{q}_3\)};
    \node at (13,1) [above] {\(\bm{q}_4\)};
    \node at (0.4,3) [below] {\(e_1\)};
    \node at (10,2) {\(e_5\)};
  \end{tikzpicture}
\end{center}

\begin{center}
  \begin{tabular}{ccccccc}
    Edge &  \(e_1\) & \(e_2\) & \(e_3\) & \(e_4\) & \(e_x\) & \(e_5\) \\
    \hline
    Affine length & \(1\) & \(0\) & \(0\) & \(0\) & \(0\) & \(3\) \\
    Affine displacement & \(1/5\) & \(2/5\) & \(3/5\) & \(4/5\) &
    \(1\) & \(1\)
  \end{tabular}
\end{center}

Since \(e_x\) moved an affine distance \(1\), the curve \(D_x\) is now
part of the branch locus, and the singularity above the
\(\frac{1}{25}(1,4)\) is exactly as in Case I. The intersection of
\(\pi^{-1}_*B\) with \(C_5\) now consists of four points (counted with
multiplicity): \(C_5\cap D_x\) and the three zeros of the cubic
\(\Omega\) restricted to \(C_5\). We distinguish various cases
according to how these roots collide, but these are precisely the
cases we analysed for \(\PP(1,1,4)\), Case II.

\paragraph{\(\PP(1,4,25)\), Case III.} In this case, the edges
\(e_2,e_3,e_4,e_x,e_5\) all move normally inwards until they pass
through \(\bm{q}_1\). In particular, \(e_x\) has moved a total affine
distance of \(5\) from its original position in \(\Pi(1,2,5)\), so
that \(D_x\) appears with multiplicity \(5\) in the branch locus. This
means that the double cover is not normal.

\paragraph{\(\PP(1,4,25)\), Case IV.} For a generic \(\Omega\) with
\(\bm{q}_6\in\OP{supp}(\Omega)\) and
\(\bm{q}_2\not\in\OP{supp}(\Omega)\), the data will be:

\begin{center}
  \begin{tikzpicture}[scale=0.6]
    \foreach \x in {0,1,2,...,19} {\node at (\x,0) {\(\bullet\)};}
    \foreach \x in {0,1,...,12} {\node at (\x,1) {\(\cdot\)};}
    \foreach \x in {0,1,...,6} {\node at (\x,2) {\(\cdot\)};}
    \foreach \x in {0} {\node at (\x,3) {\(\bullet\)};}
    \node at (1,3) {\(\circ\)};
    \node at (7,2) {\(\bullet\)};
    \node at (13,1) {\(\cdot\)};
    \node at (20,0) {\(\bullet\)};
    \foreach \y in {0,1,2,3} {\node at (0,\y) {\(\bullet\)};}
    \draw (0,0) -- (20,0) node[right] {\(\bm{q}_6\)} -- (7,2) node [above] {\(\bm{q}_3\)} -- (0,3) node[above] {\(\bm{q}_1\)} -- cycle;
    \node at (3,3) {\(e_2\)};
    \node at (13,1.5) {\(e_3\)};
  \end{tikzpicture}
\end{center}

\begin{center}
  \begin{tabular}{ccccccc}
    Edge &  \(e_1\) & \(e_2\) & \(e_3\) & \(e_4\) & \(e_x\) & \(e_5\) \\
    \hline
    Affine length & \(0\) & \(1\) & \(1\) & \(0\) & \(0\) & \(0\) \\
    Affine displacement & \(1/5\) & \(7/5\) & \(8/5\) & \(4/5\) & \(0\) & \(0\)
  \end{tabular}
\end{center}

Using Eq. \eqref{eq:alpha} for the pair \((Y,B)\) at the
\(\frac{1}{25}(1,4)\) singularity, we get
\[\bm{\alpha}=(-9/10,-13/10,-6/5,-3/5)^T,\] and since there are
discrepancies which are \(<-1\), the singularity of the double cover
is not log canonical by \ref{pg:necessary_condition}. Since a
non-generic \(\Omega\) is adjacent to something which is not log
canonical, it must itself have non-log canonical singularities.

\paragraph{\(\HP(5)\).} The Manetti surface \(\HP(5)\) is
obtained from \(\PP(1,4,25)\) by smoothing the
\(\frac{1}{4}(1,1)\) singularity and keeping the
\(\frac{1}{25}(1,4)\) singularity. This is almost toric, with
mirror tropicalisation:

\begin{center}
  \begin{tikzpicture}[scale=0.6]
    \foreach \x in {0,1,2,...,19} {\node at (\x,0) {\(\bullet\)};}
    \foreach \x in {1,...,13} {\node at (\x,1) {\(\cdot\)};}
    \foreach \x in {1,...,6} {\node at (\x,2) {\(\cdot\)};}
    \foreach \x in {0} {\node at (\x,3) {\(\bullet\)};}
    \node at (1,3) {\(\cdot\)};
    \foreach \y in {1,2} {\node at (0,\y) {\(\bullet\)};}
    \node at (7,2) {\(\ast\)};
    \node at (20,0) {\(\ast\)};
    \draw (0,0) -- (20,0) node[below] {\((20,0)\)} -- (0,16/5) node [left] {\((0,16/5)\)} -- cycle;
    \draw[dashed] (20,0) -- (8/3,8/3) node [midway,above,sloped] {\tiny branch cut \((13,-2)\)};
    \node at (8/3,8/3) {\(\times\)};
    \node at (0,1.5) [left] {\(e_y\)};
    \node at (9.5,0) [below] {\(e_z\)};
    \node at (2,2.8) [above] {\(e_z\)};
  \end{tikzpicture}
\end{center}

The affine monodromy across the branch cut is \(\begin{pmatrix} -25 &
-169\\ 4 & 27\end{pmatrix}\). We draw integer points as \(\ast\) if
they lie on the branch cut. The end of the branch cut (denoted
\(\times\)) is at \((8/3,8/3)\), and the edge \(e_z\) crosses the
branch cut so appears broken (we have labelled both pieces). Note that
there are now no curves \(C_5\), \(D_x\): instead, \(D_x\) has merged
with \(D_z\) to become a single curve which we will continue to denote
by \(D_z\).

As for \(\PP(1,4,25)\), the key integer points are
\[\bm{q}_1=(0,3),\quad \bm{q}_2=(1,3),\quad \bm{q}_3=(7,2),\quad
\bm{q}_4=(13,1),\quad \bm{q}_6=(20,0).\]
The generic case (when the support of \(\Omega\) is maximal) yields
the following data:

\begin{center}
  \begin{tikzpicture}[scale=0.6]
    \foreach \x in {0,1,2,...,19} {\node at (\x,0) {\(\bullet\)};}
    \foreach \x in {1,...,13} {\node at (\x,1) {\(\cdot\)};}
    \foreach \x in {1,...,6} {\node at (\x,2) {\(\cdot\)};}
    \foreach \x in {0} {\node at (\x,3) {\(\bullet\)};}
    \node at (1,3) {\(\bullet\)};
    \foreach \y in {1,2} {\node at (0,\y) {\(\bullet\)};}
    \node at (7,2) {\(\ast\)};
    \node at (20,0) {\(\ast\)};
    \draw (0,0) -- (20,0) node[below] {\((20,0)\)} -- (1,3) -- (0,3) -- cycle;
    \draw[dashed] (20,0) -- (8/3,8/3) node [midway,above,sloped] {\tiny branch cut \((13,-2)\)};
    \node at (8/3,8/3) {\(\times\)};
    \node at (0.4,3) [above] {\(e_1\)};
    \node at (1.8,2.8) [above] {\(e_4\)};
    \node at (0,1.5) [left] {\(e_y\)};
    \node at (9.5,0) [below] {\(e_z\)};
  \end{tikzpicture}
\end{center}

\begin{center}
  \begin{tabular}{ccccccc}
    Edge &  \(e_1\) & \(e_2\) & \(e_3\) & \(e_4\)  \\
    \hline
    Affine length & \(1\) & \(0\) & \(0\) & \(1\)  \\
    Affine displacement & \(1/5\) & \(2/5\) & \(3/5\) & \(4/5\) 
  \end{tabular}
\end{center}

Note: the edge \(e_z\) now terminates at \(\bm{q}_6\). Just as in
\(\PP(1,4,25)\), Case I, \(\pi^{-1}_*B\) hits \(C_1\) and \(C_4\) once
transversely and is disjoint from \(C_2\) and \(C_3\), so the same
argument as in that case yields a \(\frac{1}{50}(1,29)\) singularity
on \(X\). If we start to lose integer points from
\(\OP{supp}(\Omega)\) then nothing changes about this analysis unless
\(\bm{q}_2\not\in\OP{supp}(\Omega)\). When
\(\bm{q}_2\not\in\OP{supp}(\Omega)\), the data becomes:

\begin{center}
  \begin{tikzpicture}[scale=0.6]
    \node (q2) at (-1.5,2) {\(\bm{q}_2\)};
    \draw[->] (q2) -- (0.8,3);
    \foreach \x in {0,1,2,...,19} {\node at (\x,0) {\(\bullet\)};}
    \foreach \x in {1,...,13} {\node at (\x,1) {\(\cdot\)};}
    \foreach \x in {1,...,6} {\node at (\x,2) {\(\cdot\)};}
    \foreach \x in {0} {\node at (\x,3) {\(\bullet\)};}
    \foreach \y in {1,2} {\node at (0,\y) {\(\bullet\)};}
    \node at (7,2) {\(\ast\)};
    \node at (20,0) {\(\ast\)};
    \node at (1,3) {\(\circ\)};
    \draw (0,0) -- (20,0) node[below] {\((20,0)\)} -- (7,2) -- (0,3) -- cycle;
    \draw[dashed] (20,0) -- (8/3,8/3) node [midway,above,sloped] {\(e_3\), \small parallel to branch cut};
    \node at (8/3,8/3) {\(\times\)};
    \node at (2,2.8) [above] {\(e_2\)};
    \node at (9.5,0) [below] {\(e_z\)};
    \node at (0,1.5) [left] {\(e_y\)};
  \end{tikzpicture}
\end{center}

\begin{center}
  \begin{tabular}{ccccccc}
    Edge &  \(e_1\) & \(e_2\) & \(e_3\) & \(e_4\)  \\
    \hline
    Affine length & \(0\) & \(1\) & \(1\) & \(0\)  \\
    Affine displacement & \(1/5\) & \(7/5\) & \(8/5\) & \(4/5\) 
  \end{tabular}
\end{center}

Here we encounter some new phenomena in the picture: the singular
point of our integral affine manifold \(\times\) lies outside the
positive polygon and the edge \(e_3\) runs along the branch cut. The
branch curve hits \(C_2\) and \(C_3\) each once transversely and is
disjoint from \(C_1\) and \(C_4\), and we obtain the same non-log
canonical singularity as in \(\PP(1,4,25)\), Case IV.

\paragraph{\(\HP(13)\).}\label{pg:hp_13} The Manetti surface
\(\HP(13)\) is obtained
from \(\PP(1,25,169)\) by smoothing the \(\frac{1}{25}(1,4)\)
singularity and keeping the \(\frac{1}{169}(1,25)\)
singularity. Note that this singularity has index \(13\), which
is odd, so that \(13\) times an octic is Cartier and the
description of the double cover in
\ref{pg:partial_res_double_cover} applies. The triangle
\(\Pi(1,5,13)\) is Type A, so its integer points coincide with
\(\ZZ\Pi_A\). The mirror tropicalisation of \(\HP(13)\) has
vertices at \(\bm{p}_0=(0,0)\), \(\bm{p}_1=(104/5,0)\) and
\(\bm{p}_2=(0,40/13)\), and a branch cut emanating from
\(\bm{p}_1\) in the \((-34,5)\)-direction with monodromy matrix
\(M_{\bm{p}_1}=\begin{pmatrix} -169 & -1156\\ 25 &
  171\end{pmatrix}\). Below we show \(\Pi_A\) (vertices at
\((0,3)\), \((14,1)\) and \((20,0)\)) together with the ambient
Markov triangle in grey. Note that, even if it is not
immediately clear from the picture, the branch cut lies above
the convex hull of the integer points, and intersects it only at
\((14,1)\).

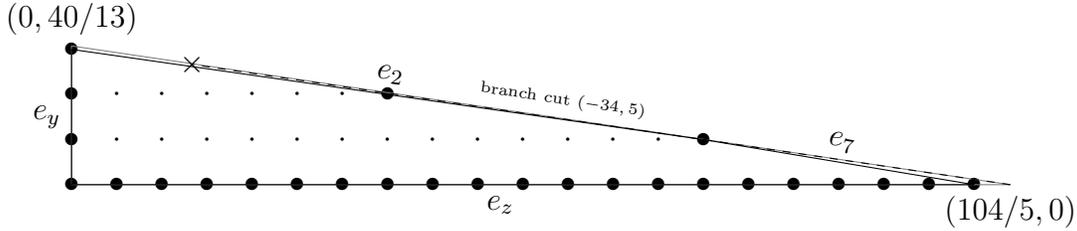
\begin{center}
  \begin{tikzpicture}[scale=0.6]
    \foreach \x in {0,1,2,...,20} {\node at (\x,0) {\(\bullet\)};}
    \foreach \x in {1,...,13} {\node at (\x,1) {\(\cdot\)};}
    \foreach \x in {1,...,6} {\node at (\x,2) {\(\cdot\)};}
    \foreach \x in {0} {\node at (\x,3) {\(\bullet\)};}
    \foreach \y in {1,2} {\node at (0,\y) {\(\bullet\)};}
    \node at (7,2) {\(\bullet\)};
    \node at (14,1) {\(\bullet\)};
    \node at (0,40/13) [above] {\((0,40/13)\)};
    \node at (104/5,0) [below] {\((104/5,0)\)};
    \draw[gray,opacity=0.8] (0,0) -- (104/5,0) -- (0,40/13) -- cycle;
    \draw (0,0) -- (20,0) -- (14,1) node [midway,above,sloped] {\(e_7\)} -- (0,3) node [midway,above,sloped] {\(e_2\)} -- cycle;
    \draw[dashed] (104/5,0) -- (8/3,8/3) node [pos=0.55,above,sloped] {\tiny branch cut \((-34,5)\)};
    \node at (8/3,8/3) {\(\times\)};
    \node at (9.5,0) [below] {\(e_z\)};
    \node at (0,1.5) [left] {\(e_y\)};
  \end{tikzpicture}
\end{center}

We resolve the \(\frac{1}{169}(1,25)\) singularity by adding
edges \(e_1,\ldots,e_7\) with inward normals
\begin{gather*}
  \rho_1=(0,-1),\, \rho_2=(-1,-7),\, \rho_3=(-5,-34),\,
  \rho_4=(-9,-61),\\
  \rho_5=(-13,-88),\, \rho_6=(-17,-115),\,
  \rho_7=(-21,-142)\end{gather*} corresponding to curves \(C_i\)
with
\[C_1^2=-7,\quad C_2^2=-5,\quad C_3^2=\cdots=C_7^2=-2.\] Note
that the edge corresponding to \(C_3\) is parallel to the branch
cut: this means that, as they move inwards, the edges
\(e_4,\ldots,e_7\) cross the branch cut, and appear as
\(M_{\bm{p}_1}^{-1}e_i\). The zero-length edges are \(e_1\)
(concentrated at \((0,3)\)), and \(e_3,\ldots,e_6\) concentrated
at \((14,1)\). The table of affine lengths and displacements
becomes:

\begin{center}
  \begin{tabular}{cccccccc}
    Edge &  \(e_1\) & \(e_2\) & \(e_3\) & \(e_4\) & \(e_5\) & \(e_6\)& \(e_7\) \\
    \hline
    Affine length & \(0\) & \(2\) & \(0\) & \(0\) & \(0\) & \(0\) & 1\\
    Affine displacement & \(1/13\) & \(7/13\) & \(8/13\) & \(9/13\)
    &  \(10/13\) & \(11/13\) & \(12/13\)\\
  \end{tabular}
\end{center}

Using Eq. \eqref{eq:alpha} to calculate the discrepancies of
\((Y,B/2)\), we find that the discrepancies of \(C_2\) and
\(C_3\) are respectively \(-31/26\) and \(-14/13\), so this pair
is not log canonical.

\paragraph{\(\HP(29)\).}\label{pg:hp_29} The Manetti surface \(\HP(29)\) has as its
mirror tropicalisation the polygon with vertices at
\(\bm{p}_0=(-16/5,0)\), \(\bm{p}_1=(20,0)\) and
\(\bm{p}_2=(80/29,80/29)\) and branch cuts emanating from
\(\bm{p}_0\) in the \((11,5)\)-direction and from \(\bm{p}_1\)
in the \((-13,2)\)-direction with monodromy matrices
\(M_{\bm{p}_0}=\begin{pmatrix}56 & -121 \\ 25 &
  -54\end{pmatrix}\) and
\(M_{\bm{p}_1}=\begin{pmatrix}-25 & -169 \\ 4 &
  27\end{pmatrix}\). All the integral points are contained in
the polygon \(\Pi_C\), whose vertices are at \((-3,0)\),
\((-1,1)\), \((2,2)\), \((7,2)\) and \((20,0)\).

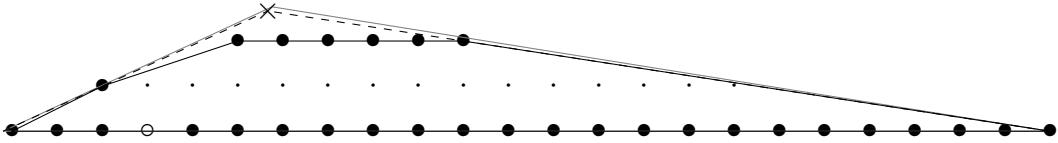
\begin{center}
  \begin{tikzpicture}[scale=0.6]
    \foreach \x in {-3,-2,-1} {\node at (\x,0) {\(\bullet\)};}
    \node at (0,0) {\(\circ\)};
    \foreach \x in {1,2,...,20} {\node at (\x,0) {\(\bullet\)};}
    \foreach \x in {-1,0,...,13} {\node at (\x,1) {\(\cdot\)};}
    \foreach \x in {2,3,...,7} {\node at (\x,2) {\(\bullet\)};}
    \node at (-1,1) {\(\bullet\)};
    \draw[gray,opacity=0.5] (-16/5,0) -- (20,0) -- (80/29,80/29) -- (-16/5,0);
    \draw[dashed] (-16/5,0) -- (8/3,8/3) node {\(\times\)};
    \draw[dashed] (20,0) -- (8/3,8/3);
    \draw (-3,0) -- (-1,1) -- (2,2) -- (7,2) -- (20,0) -- cycle;
  \end{tikzpicture}
\end{center}

The resolution of the \(\frac{1}{841}(1,637)\) singularity introduces
exceptional curves, \(C_i\), \(i=1,\ldots,10\) with
\[C_1^2=-5,\quad C_2^2=\cdots=C_6^2=-2,\quad C_7^2=-10,\quad
C_8^2=C_9^2=C_{10}^2=-2.\] These correspond to zero-length edges
\(e_1,\ldots,e_{10}\) in the tropicalisation sitting at \(\bm{p}_2\), with
inward normal vectors:
\begin{gather*}
  \rho_1 = (6,-13),\quad \rho_2=(5,-11),\quad \rho_3=(4,-9),\quad
  \rho_4=(3,-7),\quad
  \rho_5=(2,-5),\\ \rho_6=(1,-3),\quad\rho_7=(0,-1),\quad
  \rho_8=(-1,-7),\quad \rho_9=(-2,-13),\quad\rho_{10}=(-3,-19).
\end{gather*}
The edge \(e_2\) is parallel to the branch cut through \(\bm{p}_0\) and the
edge \(e_9\) is parallel to the branch cut through \(\bm{p}_1\). By the
time it has reached the inner triangle, the edge \(e_1\) has passed
through the branch cut so appears as \(M_{\bm{p}_0}e_1\), i.e. parallel to
\((2,1)\), having affine length \(1\) in the picture (note that if an
edge transforms according to \(M\) then its normal transforms
according to \((M^{-1})^T\)). Similarly, the edge \(e_{10}\) appears
as \(M_{\bm{p}_1}^{-1}e_{10}\), parallel to \((6,-1)\), having affine
length zero in the picture. The table of affine lengths and
displacements becomes:

\begin{center}
  \begin{tabular}{cccccc}
    Edge &  \(e_1\) & \(e_2\) & \(e_3\) & \(e_4\) & \(e_5\)\\
    \hline
    Affine length & \(1\) & \(0\) & \(0\) & \(0\) & \(0\)\\
    Affine displacement & \(9/29\) & \(16/29\) & \(23/29\) & \(30/29\)
    & \(37/29\) \\
    \hline
    Edge & \(e_6\)
     & \(e_7\) & \(e_8\) & \(e_9\) & \(e_{10}\) \\
    Affine length & \(1\) &
    \(5\) & \(0\) & \(1\) & \(0\)  \\
    Affine displacement & \(44/29\) & \(22/29\) & \(31/29\) &
    \(40/29\) & \(20/29\)\\
  \end{tabular}
\end{center}

Using Eq. \eqref{eq:alpha} to calculate the discrepancies of
\((Y,B/2)\), we find that many of the discrepancies are \(<-1\) (for example \(C_2\) has discrepancy \(-31/29\)) so this pair is not log canonical.

\paragraph{\(\HP(34)\).}\label{pg:hp_34} The Manetti surface
\(\HP(34)\) is obtained from
\(\PP(1,169,1156)\) by smoothing the \(\frac{1}{169}(1,25)\)
singularity and keeping the \(\frac{1}{1156}(1,169)\) singularity.
The mirror tropicalisation of \(\HP(34)\) has vertices at
\(\bm{p}_0=(0,0)\), \(\bm{p}_1=(272/13,0)\) and
\(\bm{p}_2=(0,104/34)\), and a branch cut emanating from \(\bm{p}_1\)
in the \((-89,13)\)-direction with monodromy matrix
\(M_{\bm{p}_1}=\begin{pmatrix} -1156 & -7921\\ 169 &
1158\end{pmatrix}\). Below we show the convex hull \(\Pi_A\) of
the integral points (vertices at \((0,3)\), \((14,1)\) and \((20,0)\))
together with the ambient Markov triangle in grey. Note that the
branch cut lies strictly above the convex hull of the integer points,
and never intersects it.

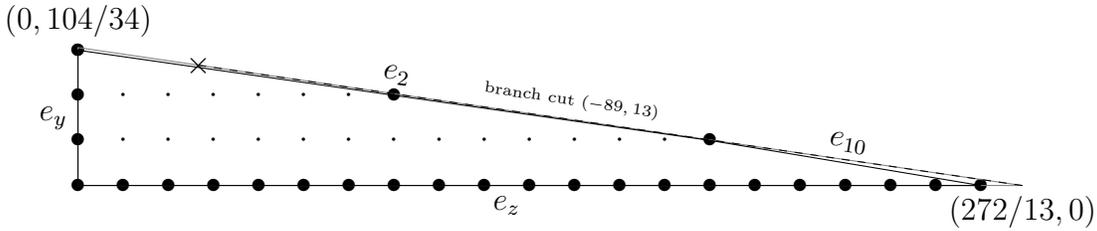
\begin{center}
  \begin{tikzpicture}[scale=0.6]
    \foreach \x in {0,1,2,...,20} {\node at (\x,0) {\(\bullet\)};}
    \foreach \x in {1,...,13} {\node at (\x,1) {\(\cdot\)};}
    \foreach \x in {1,...,6} {\node at (\x,2) {\(\cdot\)};}
    \foreach \x in {0} {\node at (\x,3) {\(\bullet\)};}
    \foreach \y in {1,2} {\node at (0,\y) {\(\bullet\)};}
    \node at (0,3) {\(\bullet\)};
    \node at (7,2) {\(\bullet\)};
    \node at (14,1) {\(\bullet\)};
    \node at (272/13,0) [below] {\((272/13,0)\)};
    \node at (0,104/34) [above] {\((0,104/34)\)};
    \draw[gray,opacity=0.8] (0,0) -- (272/13,0) -- (0,104/34) -- cycle;
    \draw (0,0) -- (20,0) -- (14,1) node [midway,above,sloped] {\(e_{10}\)} -- (0,3) node[midway,above,sloped] {\(e_2\)} -- cycle;
    \draw[dashed] (272/13,0) -- (8/3,8/3) node [pos=0.55,above,sloped] {\tiny branch cut \((-89,13)\)};
    \node at (8/3,8/3) {\(\times\)};
    \node at (0,1.5) [left] {\(e_y\)};
    \node at (9.5,0) [below] {\(e_z\)};
  \end{tikzpicture}
\end{center}

The resolution of the \(\frac{1}{1156}(1,169)\) singularity
introduces edges \(e_1,\ldots,e_{10}\) with inward normals
\begin{gather*}
  \rho_1 = (0,-1),\quad \rho_2=(-1,-7),\quad \rho_3=(-7,-48),\quad
  \rho_4=(-13,-89),\quad \\
  \rho_5=(-19,-130),\quad \rho_6=(-44,-301),\quad 
\rho_7=(-69,-472),\quad \\
  \rho_8=(-94,-643),\quad \rho_9=(-119,-814),\quad\rho_{10}=(-144,-985).
\end{gather*}
corresponding to exceptional curves, \(C_i\), \(i=1,\ldots,10\)
with
\[C_1^2=-7,\quad C_2^2=-7, \quad C_3^2=C_4^2=-2, \quad C_5^2=-3,
  \quad C_6^2=\cdots=C_{10}^2=-2\] Note that the edge \(e_4\) is
parallel to the branch cut: this means that, as they move
inwards, the edges \(e_5,\ldots,e_{10}\) cross the branch cut,
and appear as \(M_{\bm{p}_1}^{-1}e_i\). Moving the edges
normally inwards to \(\Pi_A\), edge \(e_1\) ends up with
zero-length concentrated at \((0,3)\), and \(e_3,\ldots,e_9\)
end up with zero-length concentrated at \((20,0)\). The table of
affine lengths and displacements becomes:
\begin{center}
  \begin{tabular}{ccccccc}
    Edge &  \(e_1\) & \(e_2\) & \(e_3\) & \(e_4\) & \(e_5\) & \(e_6\)  \\
    \hline
    Affine length & \(0\) & \(2\) & \(0\) & \(0\) & \(0\)  & \(0\)  \\
    Affine displacement & \(1/17\) & \(7/17\) & \(14/17\) & \(21/17\) & \(11/17\) & \(12/17\) \\
    \hline
    Edge &  \(e_7\) & \(e_8\) & \(e_9\) & \(e_{10}\)  \\
    Affine length & \(1\) & \(0\) & \(0\) & \(1\)  \\
    Affine displacement & \(13/17\) & \(14/17\) & \(15/17\) & \(16/17\)\\
  \end{tabular}
\end{center}

Using Eq. \eqref{eq:alpha} to calculate the discrepancies of
\((Y,B/2)\), we find that many discrepancies \(<-1\) (for example \(C_2\) has discrepancy \(-2657/2312\)) so this pair is not log canonical.

There are some important features in this example which did not appear
in the other examples. First, the index of the Wahl singularity is
\(34\) which is even, so it is not immediately clear why there is an
odd \(n\) such that \(nB\) is Cartier. However, since \(B\) is octic,
its homology class is divisible by \(2\), so \(17B\) is
Cartier. Second, the curve \(\eE_1\) from the non-toric blow-up in the
construction of \(\tilde{Y}\) appears as a part of the branch locus,
with multiplicity \(1\). To see why, observe that the polygon
\(\Pi_A\) has a zero-length edge concentrated at \((14,1)\) defined by
\[\tilde{b}_1\coloneqq \langle
  u,\rho_4\rangle=\left(\begin{pmatrix}8/3 \\ 8/3\end{pmatrix} -
    \begin{pmatrix} 14 \\ 1\end{pmatrix}\right)
  \cdot \begin{pmatrix} -13 \\ -89\end{pmatrix} = -1.\]
By \ref{pg:geom_mirror_trop}(ii), this means that \(\eE_1\cdot
\pi^{-1}_*B=-1\), so \(\pi^{-1}_*B\) contains \(\eE_1\) as a component
with multiplicity \(1\).

\section{Conclusion}
\label{sct:conclusion}

We conclude with some assorted remarks about the KSBA-stable
surfaces we have found.

\pg First, in all cases, the minimal resolution is an elliptic
surface of geometric genus \(3\). One can be quite explicit
about the singular fibres; the interested reader can look at
Section 12 of the arXiv version 1 of this paper for details.

\pg Not every branched double cover of \(Y=\PP(1,1,4)\),
\(\PP(1,4,25)\) or \(\HP(5)\) with a basket of non-Gorenstein
singularities described by Theorem \ref{thm:main_thm} needs to
be stable; for example the branch curve \(B\) could develop some
awful non log canonical singularity away from the singularities
of \(Y\). However, the other condition for KSBA-stability
(ampleness of \(K_X\)) always holds. To see why, note that
\(K_Y+\frac{1}{2}B\) is ample on \(Y\) since the Picard rank of
\(Y\) is \(1\), so any line bundle of positive degree is ample
(in the class group of \(Y\) we have
\(K_Y+\frac{1}{2}B=(-3+8/2)t=t\) where \(t\) is the product of
Markov numbers corresponding to the singularities of \(Y\)). Now
\(K_X=f^*\left(K_Y+\frac{1}{2}B\right)\), and since the pullback
of an ample divisor by a finite map is ample (see {\cite[Lemma
  30.17.2, tag 0B5V]{Stacks}}), \(K_X\) is ample. Therefore if
\((Y,B/2)\) is log canonical then \(X\) is KSBA-stable. We
verified this in each of the cases studied in Section
\ref{sct:classif} assuming that \(B\) generic amongst octics
satsfying the constraints of the case.

\pg\label{pg:smoothable} We will now analyse how the strata we found sit in the KSBA
boundary. First, note that {\em any octic double Manetti surface
  is smoothable.}
\begin{proof}
  Suppose that \(Y_0\) is a Manetti surface and let \(n\) be the
  product of the Markov numbers associated to the singularities
  of \(Y_0\), that is:
  \[n=\begin{cases}
      abc&\mbox{ if }Y_0=\PP(a^2,b^2,c^2)\\
      ab&\mbox{ if }Y_0=\HP(a,b)\\
      c&\mbox{ if }Y_0=\HP(c).
    \end{cases}\] Let \(B_0\subset Y_0\) an octic curve and
  \(f\colon X\to Y_0\) the double cover of \(Y_0\) branched over
  \(B_0\). Let \(\mathcal{Y}\to \Delta\) be a \(\QQ\)-Gorenstein
  smoothing of \(Y_0\) with general fibre \(\PP^2\). By
  {\cite[Theorem 1.6(4)]{DeVlemingStapleton}}, after possibly
  making a base change, we can find a divisor
  \(\mathcal{B}\subset\mathcal{Y}\) extending \(B_0\subset Y_0\)
  if and only if the Weil divisor class \([B_0]\) is divisible
  by \(n\). Since \(B_0\) is an octic curve, \([B_0]=8n\), so
  \(\mathcal{B}\) exists, and we can use the Alexeev--Pardini
  construction outlined in \ref{pg:ap_double_covers} to take the
  double cover of \(\mathcal{Y}\) branched over
  \(\mathcal{B}\). This yields a 3-fold
  \(\mathcal{X}\to\Delta\); it is not completely obvious that
  the central fibre \(X_0\) is the same as the \(S_2\) double
  cover of \(Y_0\) branched along \(B_0\). Note that the
  Alexeev--Pardini cover is obtained by taking
  \(\OP{Spec}_{\mathcal{Y}}\left(\OO\oplus\OO(-\mathcal{B})\right)\),
  so to show that it behaves well under restriction, we need to
  show that
  \(\OO_{\mathcal{X}}\left(-\mathcal{B}\right)|_{X_0} =
  \OO_{X_0}\left(-B_0\right)\). This follows from
  {\cite[Proposition 2.79]{KollarFamilies}} since
  \(\mathcal{B}\) is \(\QQ\)-Cartier. Therefore \(\XX\) is a
  deformation of the octic double Manetti surface. If \(B_0\) is
  chosen generically then \(\mathcal{B}|_{Y_t}\) will be a
  generic octic and hence smooth, so \(X_t\) will be smooth.
\end{proof}

\pg The moduli space of smooth octic double planes is
\(36\)-dimensional: the projective space of octics is the
\(44\)-dimensional space \(\PP(\OP{Sym}^8\left(V^*)\right)\)
where \(V\) is the standard representation of the
\(9\)-dimensional group \(GL(3)\), so
\(\dim(\PP\left(\OP{Sym}^8(V^*)\right)/PGL(3)=44-8=36\). For
each basket of singularities in Theorem \ref{thm:main_thm} we
obtain a stratum of the KSBA boundary: we write \(D^Y_C\) for
this stratum where \(Y\in\{\PP(1,1,4),\PP(1,4,25),\HP(5)\}\) is
the base of the branched cover and \(C\in\{I,IIa,IIb,IIc\}\) is
the basket (\(Y=\HP(5)\) only has one case which we label
\(I\)). When \(C=IIa,b,c\) and \(Y=\PP(1,1,4)\) (respectively
\(Y=\PP(1,4,25)\)), we only consider the generic case where the
branch curve is smooth at the points where it intersects \(C_1\)
(respectively \(C_5\)): the deeper strata seem harder to
control. Of course these strata will themselves be stratified
according to any additional Gorenstein singularities that may
appear in the basket; see Anthes \cite{Anthes} for how the big
open stratum with no non-Gorenstein surfaces is stratified. We
do not discuss this here.

\paragraph{Proposition.} {\em The strata have the following
  dimensions:}

\begin{align*}
  \dim\left(D^{\PP(1,1,4)}_I\right)&=35 &
                                          \dim\left(D^{\PP(1,1,4)}_{IIa}\right)&=34
  &   \dim\left(D^{\PP(1,1,4)}_{IIb}\right)&=33 &
                                                  \dim\left(D^{\PP(1,1,4)}_{IIc}\right)&=32
  \\
  \dim\left(D^{\PP(1,4,25)}_I\right)&=34 &
                                          \dim\left(D^{\PP(1,4,25)}_{IIa}\right)&=33
  &   \dim\left(D^{\PP(1,4,25)}_{IIb}\right)&=32 &
                                                  \dim\left(D^{\PP(1,4,25)}_{IIc}\right)&=31
\end{align*}
\[\dim\left(D^{\HP(5)}_I\right)=35\]

{\em Moreover:
  \begin{gather*}
    D^Y_{IIc}\subset \overline{D^Y_{IIb}},\quad
    D^Y_{IIb}\subset \overline{D^Y_{IIa}},\quad D^Y_{IIa}\subset
    \overline{D^Y_{I}}\quad\mbox{for}\quad Y\in\{\PP(1,1,4),\PP(1,4,25)\},\\
    D^{\PP(1,4,25)}_C\subset\overline{D^{\PP(1,1,4)}_C}\quad\mbox{for}\quad
    C\in\{I,IIa,IIb,IIc\},\\
    \mbox{and}\quad
  D^{\PP(1,4,25)}_I\subset
  \overline{D^{\PP(1,1,4)}_I}\cap\overline{D^{\HP(5)}_I}.
  \end{gather*}
}
\begin{proof}
  For any Manetti surface, the space of generic ``octics'' is
  \(45\)-dimensional, which becomes \(44\) after
  projectivising. This is because a (GHK)-basis is in bijection
  with the integral points of the corresponding Markov triangle;
  this is easily seen to be \(45\) for \(\Pi(1,1,1)\), and the
  number of integral points is invariant under mutation.
  
  For a fixed surface \(Y=\PP(1,1,4)\) or \(\PP(1,4,25)\), Cases
  II(a-c) arise by imposing one, two or three codimension \(1\)
  conditions; for \(Y=\PP(1,1,4)\) where we write the equation
  of the branch curve as a polynomial in Cox variables
  \(\Omega(x,y,z,w)\), these conditions are (a) the absence of
  the monomial \(z^4\) from \(\Omega\), (b) both the absence of
  \(z^4\) and the vanishing of the discriminant of the
  homogeneous quartic \(\Omega(x,y,1,0)\), and (c) both the
  other conditions and the vanishing of the second subresultant
  between \(\Omega(x,y,1,0)\) and its derivative, which detects
  the presence of two common roots. The conditions for
  \(\PP(1,4,25)\) are similar but harder to write out because we
  need more Cox variables.

  Since the surface \(\PP(1,4,25)\) is a common degeneration of
  both \(\HP(5)\) and \(\PP(1,1,4)\), and since for a fixed
  \(Y\) the codimension \(1\) conditions for the three cases
  II(a-c) are nested, most of the adjacencies between the strata
  are clear. One can also write down the degeneration from
  \(\PP(1,1,4)\) to \(\PP(1,4,25)\) explicitly in coordinates in
  \(\PP(1,1,4,5)\) as in {\cite[Example
    6.3]{HackingExceptional}} and see that these conditions for
  II(a-c) specialise, so if \((\PP(1,4,25),B)\) arises as a KSBA
  limit of a sequence of pairs \((\PP(1,1,4),B_k)\) of type IIa
  (respectively IIb, IIc), \(B\) will be in the closure of
  \(D^{\PP(1,1,4)}_{IIa}\) (respectively
  \(D^{\PP(1,1,4)}_{IIb}\), \(D^{\PP(1,1,4)}_{IIc}\)).

  As we move through the cases for a fixed \(Y\), the dimension
  drops by \(1\) because we are imposing a single additional
  constraint on the coefficients of \(\Omega\) and in the
  generic situation these constraints are
  independent.\footnote{Recall that we are not delving further
    into the situation where the branch curve can become
    singular along \(C_1\) or \(C_5\): it seems harder to
    capture this behaviour in terms of the coefficients of
    \(\Omega\).} It remains to compute the dimensions for the
  top strata \(D^Y_I\), and for that it is sufficient to show:
  \[\dim(\OP{Aut}(\PP(1,1,4))=\dim(\OP{Aut}(\HP(5)))=9,\qquad
    \dim(\OP{Aut}(\PP(1,4,25)))=10.\] The general automorphism
  of \(\PP(1,1,4)\) is:
  \[[x:y:z]\mapsto
    [a_1x+a_2y:a_3x+a_4y:a_5z+a_6x^4+a_7x^3y+\cdots+a_{10}y^5]\]
  (we thank S\"{o}nke Rollenske for pointing the five terms we
  originally missed) which has \(9\) free parameters up to an
  overall scale factor.

  The general automorphism
  of \(\PP(1,4,25)\) is:
  \[[x:y:z]\mapsto
    [a_1x:a_2y+a_3x^4:a_4z+a_5x^{25}+a_6x^{21}y+\cdots+a_{11}xy^6]\]
  which has \(10\) free parameters up to scale.

  To understand the automorphism group of \(\HP(5)\), recall
  that its minimal resolution is obtained from \(\FF_7\) by a
  sequence of two toric and one non-toric blow-ups (see
  \ref{pg:hp5_minimal_model}--we use the notation from that
  paragraph). Let \(Y\) denote the result of performing the two
  toric blow-ups to \(\FF_7\). An automorphism of \(Y\) which
  fixes the \(-1\)-curve denoted \(F\) in Figure \ref{fig:hp5}
  lifts to a unique automorphism of \(\HP(5)\) and every
  automorphism of \(\HP(5)\) arises this way, because the three
  points \(F\cap B\), \(F\cap E\) and \(F\cap G\) must be fixed,
  so any automorphism of \(\HP(5)\) fixes \(F\) pointwise and
  can be blown-down.

  \begin{figure}[htb]
    \begin{center}
      \begin{tikzpicture}
        \node (a4) at (0,0) {\(\bullet\)};
        \node (b4) at (8,0) {\(\bullet\)};
        \node (c4) at (1,1) {\(\bullet\)};
        \node (d4) at (0,1) {\(\bullet\)};
        \node (e4) at (6,0.5) {\(\bullet\)};
        \node (f4) at (3.5,0.8) {\(\bullet\)};
        \draw (a4.center) -- (b4.center) node [midway,below] {\(z\)} -- (e4.center) node [midway,above right] {\(w\)} -- (f4.center) node [midway,above] {\(v\)} -- (c4.center) node [midway,above right] {\(y\)} -- (d4.center) node [midway,above] {\(u\)} -- (a4.center) node [midway,left] {\(x\)};
        \node (g1) at ((5,0.3) {};
        \node (g2) at (6,1.4) {};
      \end{tikzpicture}
      \caption{The moment polygon of \(Y\).}
      \label{fig:aut_calc}
    \end{center}
  \end{figure}
  
  It therefore suffices to find the group of automorphisms of
  \(Y\) fixing \(F\) pointwise. Since \(Y\) is toric, with
  hexagonal moment polygon, there is a GIT model for \(Y\) as
  \(\CC^6\sslash\GG_m^4\). Let \(x,y,z,u,v,w\) be the Cox
  variables associated to the edges of the moment hexagon as
  shown in Figure \ref{fig:aut_calc}. The inward normals
  associated to these edges are
  \[\rho_x=(0,1),\,\rho_y=(-1,-7),\,\rho_z=(0,1),\,
    \rho_u=(0,-1),\,\rho_v=(-2,-13),\,\rho_w=(-1,-6),\]
  from which we can read off the weights of the Cox variables
  under the \(\GG_m^4\)-action:

  \begin{center}
    \begin{tabular}{cccccc}
      \(x\) & \(y\) & \(z\) & \(u\) & \(v\) & \(w\)\\
      \(1\) & \(1\) & \(7\) & \(0\) & \(0\) & \(0\)\\
      \(0\) & \(0\) & \(1\) & \(1\) & \(0\) & \(0\)\\
      \(2\) & \(0\) & \(13\) & \(0\) & \(1\) & \(0\)\\
      \(1\) & \(0\) & \(6\) & \(0\) & \(0\) & \(1\)\\
    \end{tabular}
  \end{center}
  
  The unstable locus is \(\{vwyz=vwxz=uwxz=uxyz=uvxy=uvwy=0\}\).
  An automorphism has the form
  \[[x:y:z:u:v:w]\mapsto \left[ax+byv^2w:cy:dz+\sum_{i=0}^6
      e_iuv^{13-2i}w^{6-i}x^iy^{7-i}:fu:gv:hw\right]\] for some
  \(a,b,c,d,e_0,\ldots,e_6,f,g,h\). Invertibility implies that
  \(c,f,g,h\) are all nonzero, and using the \(\GG_m^4\)-action
  we can assume they are all equal to \(1\). This leaves \(10\)
  free parameters \(a,b,d,e_0,\ldots,e_6\). The points
  \(F\cap B=[1:0:1:1:0:1]\) and \(F\cap E=[1:1:1:1:0:0]\) are
  automatically fixed, so it suffices to fix another point on
  \(F\), say \([1:1:1:1:0:1]\). This reduces to the condition
  that
  \[[a:1:d:1:0:1]=[1:1:1:1:0:1].\] If we pick a square root of
  \(a\) then act using \((1,1,1/\sqrt{a},1)\in\GG_m^4\) we get
  \[[1:1:d/\sqrt{a}:1:0:1]=[1:1:1:1:0:1]\] and have used up all
  of the freedom in the group action, so the relevant subgroup
  is defined by the condition that \(d^2=a\). This leaves \(9\)
  free parameters, as required.
\end{proof}

\appendix
\section{Example: Manetti surfaces as almost toric surfaces}
\label{sct:almost_toric}

To help the reader, we include a detailed example of the first
really nontrivial almost toric Manetti surface, \(\HP(5)\).

\pg\label{pg:hp5_minimal_model} Take \(\PP(1,4,25)\), and smooth
the \(\frac{1}{4}(1,1)\) singularity. This gives the Manetti
surface \(\HP(5)\) with a single \(\frac{1}{25}(1,4)\)
singularity. We obtain \(\HP(5)\) as follows. Take the
Hirzebruch surface \(\FF_7\) and let \(A+B+C+D\) be its toric
boundary divisor, where \(A^2=-7\), \(B^2=D^2=0\) and
\(C^2=7\). Blow up at \(B\cap C\), label the exceptional curve
\(E\) and label all strict transforms with the names of their
blown-down selves. Blow up at \(B\cap E\) and call the new
\(-1\)-curve \(F\). Finally blow up an interior
point\footnote{i.e. a point of \(F\) which is not a toric fixed
  point.} of \(F\) and call the final \(-1\)-curve \(G\). We get
the configuration of curves shown in Figure \ref{fig:hp5}.

\begin{figure}[htb]
  \begin{center}
    \begin{tikzpicture}[scale=0.8]
      \begin{scope}
        \node (a1) at (0,0) {\(\bullet\)};
        \node (b1) at (8,0) {\(\bullet\)};
        \node (c1) at (1,1) {\(\bullet\)};
        \node (d1) at (0,1) {\(\bullet\)};
        \draw (a1.center) -- (b1.center) node [midway,below] {\(C\)} -- (c1.center) node [midway,above right] {\(B\)} -- (d1.center) node [midway,above] {\(A\)} -- (a1.center) node [midway,left] {\(D\)};
        \node at (0.5,0.5) {\(\FF_7\)};
      \end{scope}
      \begin{scope}[shift={(0,2)}]
        \node (a2) at (0,0) {\(\bullet\)};
        \node (b2) at (8,0) {\(\bullet\)};
        \node (c2) at (1,1) {\(\bullet\)};
        \node (d2) at (0,1) {\(\bullet\)};
        \node (e2) at (6,0.5) {\(\bullet\)};
        \draw (a2.center) -- (b2.center) -- (e2.center) node [midway,above right] {\(E\)} -- (c2.center) -- (d2.center) -- (a2.center);
      \end{scope}
      \begin{scope}[shift={(0,4)}]
        \node (a3) at (0,0) {\(\bullet\)};
        \node (b3) at (8,0) {\(\bullet\)};
        \node (c3) at (1,1) {\(\bullet\)};
        \node (d3) at (0,1) {\(\bullet\)};
        \node (e3) at (6,0.5) {\(\bullet\)};
        \node (f3) at (3.5,0.8) {\(\bullet\)};
        \draw (a3.center) -- (b3.center) -- (e3.center) -- (f3.center) node [pos=0.7,above] {\(F\)} -- (c3.center) -- (d3.center) -- (a3.center);
        \node (h) at (5,0.6) {};
      \end{scope}
      \begin{scope}[shift={(0,6)}]
        \node (a4) at (0,0) {\(\bullet\)};
        \node (b4) at (8,0) {\(\bullet\)};
        \node (c4) at (1,1) {\(\bullet\)};
        \node (d4) at (0,1) {\(\bullet\)};
        \node (e4) at (6,0.5) {\(\bullet\)};
        \node (f4) at (3.5,0.8) {\(\bullet\)};
        \draw (a4.center) -- (b4.center) node [midway,below] {\(6\)} -- (e4.center) node [midway,above right] {\(-2\)} -- (f4.center) node [pos=0.7,above] {\(-2\)} -- (c4.center) node [midway,above right] {\(-2\)} -- (d4.center) node [midway,above] {\(-7\)} -- (a4.center) node [midway,left] {\(0\)};
        \node (g1) at ((5,0.3) {};
        \node (g2) at (6,1.4) {};
        \draw (g1.center) -- (g2.center) node [right] {\(G\)};
        \node at (6,1.4) [left] {\(-1\)};
      \end{scope}
      \draw[dotted,->-] (e2.center) -- (b1.center);
      \draw[dotted,->-] (b2.center) -- (b1.center);
      \draw[dotted,->-] (f3.center) -- (e2.center);
      \draw[dotted,->-] (e3.center) -- (e2.center);
      \draw[dotted,->-] (g1.center) -- (h.center);
      \draw[dotted,->-] (g2.center) -- (h.center);
    \end{tikzpicture}
    \caption{The minimal resolution of \(\HP(5)\) built from
      \(\FF_7\) by two toric blow-ups and one non-toric
      blow-up. Curves and their self-intersections are
      labelled.}\label{fig:hp5}
  \end{center}
\end{figure}

The surface \(\HP(5)\) is obtained by contracting the divisor
\(A+B+F+E\) (a Wahl chain with self-intersections
\(-7,-2,-2,-2\)). Note that \(\HP(5)\) is uniquely determined up to
isomorphism by this recipe: although we had to choose the interior
point of \(F\) for the final, non-toric, blow-up, any two choices are
related by the torus action.

\pg There is a (Symington-style) almost toric picture of
\(\HP(5)\) (Figures \ref{fig:hp5_mutations} and
\ref{fig:hp5_resolution}). Start with the moment polygon of
\(\PP(1,4,25)\): a triangle with vertices at \((0,0)\),
\((25,0)\) and \((0,4)\). Make a generalised nodal trade at the
corner \((25,0)\) by introducing a focus-focus singularity
connected to this corner along a branch cut pointing in the
\((-13,2)\) direction (you find this direction by adding the two
primitive edge vectors \((-1,0)\) and \((-25,4)\) emanating from
this vertex: \((-1,0)+(-25,4)=2(-13,2)\)). This is an almost
toric base diagram for the orbifold \(\HP(5)\) (Figure
\ref{fig:hp5_mutations}(f))--it can be obtained by a sequence of
mutations from the standard moment triangle for \(\PP^2\)
(Figure \ref{fig:hp5_mutations} (a-f)).

If we take the minimal resolution then the base diagram changes
(Figure \ref{fig:hp5_resolution}(b)): we make the top corner Delzant
by chopping off pieces, introducing new edges pointing in the
directions:
\[(1,0),\quad (7,-1),\quad (13,-2),\quad (19,-3).\]
These correspond to the \(-7,-2,-2,-2\) curves \(A,B,F,E\). Note that
the edge \(F\) is parallel to the branch cut. If we perform a suitable
change of branch cut (Figure \ref{fig:hp5_resolution}(c)) and a
non-toric blow-down (Figure \ref{fig:hp5_resolution}(d)), the result
is a toric surface which can be blown-down to \(\FF_7\) (Figure
\ref{fig:hp5_resolution}(e)).

\begin{figure}[htb]
  \begin{center}
    \begin{tikzpicture}[scale=0.4]
      \begin{scope}[shift={(12,6)}]
        \node at (-2,6) {(a)};
        \draw (0,0) -- (10,0) node {\tiny \(\bullet\)} -- (0,10)  node {\tiny \(\bullet\)} -- (0,0)  node {\tiny \(\bullet\)};
        \draw[blue,dotted,->] (8,0.5) -- (8,-1.5) -- (-4,-1.5) node[pos=0.5,above] {nodal trade} -- (-4,0.5);
        \node at (10,0) [right] {\tiny \((10,0)\)};
        \node at (0,10) [left] {\tiny \((0,10)\)};
        \node at (0,0) [above right] {\tiny \((0,0)\)};
      \end{scope}
      \begin{scope}[shift={(0,6)}]
        \node at (-2,6) {(b)};
        \draw (0,0) -- (10,0) -- (0,10) node {\tiny \(\bullet\)} -- (0,0) node {\tiny \(\bullet\)};
        \draw[dashed,blue] (10,0) -- (2,4) node [pos=0.5,below,sloped] {\tiny \((-2,1)\)};
        \node[blue] at (2,4) {\(\times\)};
        \draw[blue,dotted,->] (3,3.5) arc [radius=1,start angle=-30,end angle=150] to[out=-120,in=90] (1.5,-1.9);
        \node[blue] at (3.2,-0.8) {mutate};
      \end{scope}
      \node at (-2,3) {(c)};
      \draw (0,0) -- (20,0) node {\tiny \(\bullet\)} -- (0,5) -- (0,0) node {\tiny \(\bullet\)};
      \node at (0,0) [above right] {\tiny \((0,0)\)};
      \node at (20,0) [right] {\tiny \((20,0)\)};
      \node at (0,5) [left] {\tiny \((0,5)\)};
      \draw[dashed,blue] (0,5) -- (2,4) node [pos=0.75,below,sloped] {\tiny \((2,-1)\)};
      \node[blue] at (2,4) {\(\times\)};
      \draw[red,dotted,->] (19,-0.5) -- (19,-5.5) node [sloped,midway,above] {nodal trade};
      \begin{scope}[shift={(0,-6)}]
        \node at (-2,3) {(d)};
        \draw (0,0) node {\tiny \(\bullet\)} -- (20,0) -- (0,5) -- cycle;
        \draw[dashed,blue] (0,5) -- (2,4) node {\(\times\)};
        \draw[dashed,red] (20,0) -- (5,3) node {\(\times\)};
        \draw[red,dotted,->] (5.5,2.8) arc [radius=0.5,start angle=-30,end angle=150] to[out=-120,in=90] (4.5,-2.9);
        \node[red] at (6.3,-1.5) {mutate};
      \end{scope}
      \begin{scope}[shift={(0,-12)}]
        \node at (-2,3) {(e)};
        \draw (0,0) node {\tiny \(\bullet\)} -- (25,0) -- (0,4) -- cycle;
        \node at (0,0) [above right] {\tiny \((0,0)\)};
        \node at (0,4) [above] {\tiny \((0,4)\)};
        \node at (25,0) [above] {\tiny \((25,0)\)};
        \draw[dashed,red] (0,4) -- (5,3) node {\(\times\)};
        \draw[dashed,blue] (25,0) -- (8.1,2.6) node [pos=0.75,below] {\tiny \((-13,2)\)};;
        \node[blue] at (8.1,2.6) {\(\times\)};
        \draw[red,dotted,->] (3,3.2) -- (3,-2) node [pos=0.75,right] {nodal trade};
      \end{scope}
      \begin{scope}[shift={(0,-18)}]
        \node at (-2,3) {(f)};
        \draw (0,0) node {\tiny \(\bullet\)} -- (25,0) -- (0,4) node {\(\star\)} -- cycle;
        \draw[dashed,blue] (25,0) -- (8.1,2.6) node [pos=0.75,below] {\tiny \((-13,2)\)};;
        \node[blue] at (8.1,2.6) {\(\times\)};
        \node at (15,3) {\(\HP(5)\)};
      \end{scope}
    \end{tikzpicture}
    \caption{A sequence of mutations to get the almost toric base
      diagram for \(\HP(5)\). The point marked \(\star\) is a
      \(\frac{1}{25}(1,4)\) singularity. The points marked \(\bullet\)
      are Delzant vertices. The points marked \(\times\) are
      focus-focus fibres.}
    \label{fig:hp5_mutations}
  \end{center}
\end{figure}

\begin{figure}[htb]
  \begin{center}
    \begin{tikzpicture}[scale=0.4]
      \node at (-2,3) {(a)};
      \draw (0,0) node {\tiny \(\bullet\)} -- (25,0) -- (0,4) node {\(\star\)} -- cycle;
      \draw[dashed] (25,0) -- (8.1,2.6) node [pos=0.75,below] {\tiny \((-13,2)\)};;
      \node at (8.1,2.6) {\(\times\)};
      \draw[->,dotted] (0,4) -- (0.25,-2) node [pos=0.85,left] {\tiny \((1,0)\)};
      \draw[->,dotted] (0,4) -- (3,-2.5) node [pos=0.85] {\tiny \((7,-1)\)};
      \draw[->,dotted] (0,4) -- (7,-3.2) node [pos=0.85,right] {\tiny \((13,-2)\)};
      \draw[->,dotted] (0,4) -- (15,-4.4) node [pos=0.85,right] {\tiny \((19,-3)\)};
      \node at (16,-1.5) {resolve by cutting in};
      \node at (18,-2.4) {directions indicated};
      \node at (25.2,1.5) [right] {\(\HP(5)\)};
      \begin{scope}[shift={(0,-6)}]
        \node at (-2,3) {(b)};
        \draw (0,0) node {\tiny \(\bullet\)} -- (0,3.8)  node {\tiny \(\bullet\)} -- ++ (0.5,0) node {\tiny \(\bullet\)} -- ++
        (3.5,-0.5) node {\tiny \(\bullet\)}  -- ++ (6.5,-1) node (a) {\tiny \(\bullet\)}
        -- ++ (9.5,-1.5) node (b) {\tiny \(\bullet\)}  -- (25,0) -- cycle;
        \draw[dashed] (25,0) -- (8.1,2.6) node [pos=0.75,below] {\tiny \((-13,2)\)};;
        \node at (8.1,2.6) {\(\times\)};
        \node at (25.2,1.5) [right] {\(\widetilde{\HP}(5)\)};
        \draw[->,dotted] (8.1,2.5) -- (8.1,-3) node [pos=0.65,left] {change branch cut};
      \end{scope}
      \begin{scope}[shift={(0,-12)}]
        \node at (-2,3) {(c)};
        \draw (0,0) node {\tiny \(\bullet\)} -- (0,3.8) node {\tiny \(\bullet\)} -- ++ (0.5,0) node {\tiny \(\bullet\)} -- ++ (3.5,-0.5) node {\tiny \(\bullet\)}  -- ++ (3.25,-1/2) ++ (3.2,-1/2) -- ++ (3.25,-1/2) node (a1) {\tiny \(\bullet\)} -- (25.2,0) node (b1) {\tiny \(\bullet\)} -- (0,0);
        \draw[dotted,->] (a) -- (a1);
        \draw[dotted,->] (b) -- (b1);
        \draw[dashed] (7.25,2.8) -- (8.1,2.6) node {\(\times\)} -- (10.5,2.3);
        \draw[dotted,->] (8.1,1.5) -- (8.1,-3) node [pos=0.7,right] {non-toric blow-down};
        \node at (25.2,1.5) [right] {\(\widetilde{\HP}(5)\)};
        \node at (17,2) {\tiny \((6,-1)\)};
      \end{scope}
      \begin{scope}[shift={(0,-18)}]
        \node at (-2,3) {(d)};
        \draw (0,0) node {\tiny \(\bullet\)} -- (0,3.8) node {\tiny \(\bullet\)} -- ++ (0.5,0) node {\tiny \(\bullet\)} -- ++ (3.5,-0.5) node (x) {\tiny \(\bullet\)}  -- ++ (3.25,-1/2) -- ++ (3.2,-1/2) -- ++ (3.25,-1/2) node (y) {\tiny \(\bullet\)} -- (25.2,0) node (z) {\tiny \(\bullet\)} -- (0,0);
        \node at (0.25,3.8) [above] {\tiny \((1,0)\)};
        \node at (2.2,3.5) [above] {\tiny \((7,-1)\)};
        \node at (10,3) {\tiny \((13,-2)\)};
        \node at (17,2) {\tiny \((6,-1)\)};
        \node at (0.25,3.8) [below] {\tiny \(A\)};
        \node at (2.2,3.6) [below] {\tiny \(B\)};
        \node at (10,2.5) [below] {\tiny \(F\)};
        \node at (17,1.5) [below] {\tiny \(E\)};
        \node at (7,0) [above] {\tiny \(C\)};
        \node at (0,1.5) [right] {\tiny \(D\)};
      \end{scope}
      \begin{scope}[shift={(0,-24)}]
        \node at (-2,3) {(e)};
        \draw (0,0) node {\tiny \(\bullet\)} -- (0,3.8) node {\tiny \(\bullet\)} -- (0.5,3.8) node {\tiny \(\bullet\)} -- (25.2,0) node (w) {\tiny \(\bullet\)} -- (0,0);
        \node at (17,2) {\tiny \((7,-1)\)};
        \draw[->,dotted] (x) -- (w);
        \draw[->,dotted] (y) -- (w);
        \draw[->,dotted] (z) -- (w) node[pos=0.4,left] {contract};
        \node at (25.2,1.5) [right] {\(\FF_7\)};
      \end{scope}
    \end{tikzpicture}
    \caption{The almost toric base diagram for the minimal resolution
      \(\widetilde{\HP}(5)\) and its blow-down to \(\FF_7\). Since
      this is drawn to scale, the branch cuts are hard to see. In (c),
      there are two branch cuts related by the affine monodromy around
      the focus-focus singularity. These form two legs of a triangular
      ``hole'' whose base is parallel to \((13,-2)\). Non-toric
      blow-down means filling in this hole. The labels in (d) refer to
      Figure \ref{fig:hp5} (with \(G\) blown-down).}
    \label{fig:hp5_resolution}
  \end{center}
\end{figure}

\clearpage

\bibliographystyle{plain}
\bibliography{double_octics}

\end{document}